\documentclass[11pt, reqno]{amsart}
\usepackage{amssymb}
\usepackage{amsmath}
\usepackage{amsthm}
\usepackage{braket}
\usepackage{pdfpages}
\usepackage{graphicx}
\usepackage{caption}
\usepackage{subcaption}
\usepackage{mathrsfs} 
\usepackage{tikz-cd} 
\usepackage{bm}
\usepackage{enumitem}
\usepackage{mathtools}
\usepackage{pgfplots}
\usepackage{import}
\usepackage{xifthen}
\usepackage{transparent}
\usepackage{mathtools}
\usepackage{tabularx}
\usepackage{hhline}

\newcolumntype{L}[1]{>{\raggedright\arraybackslash}p{#1}}
\newcolumntype{C}[1]{>{\centering\arraybackslash}p{#1}}
\newcolumntype{R}[1]{>{\raggedleft\arraybackslash}p{#1}}

\usepackage[utf8]{inputenc}
\usepackage[english]{babel}

\usepackage{hyperref}

\usepackage{lmodern,babel,adjustbox,booktabs,multirow}

\numberwithin{equation}{section}
\theoremstyle{plain}
\newtheorem{theorem}{Theorem}[section]

\theoremstyle{theorem}
\newtheorem{prop}[theorem]{Proposition}
\newtheorem{lem}[theorem]{Lemma}
\newtheorem{cor}[theorem]{Corollary}

\newtheorem*{question*}{Question}

\theoremstyle{definition}
\newtheorem{defn}[theorem]{Definition}

\newtheorem{assumption}[theorem]{Assumption}

\newcommand{\R}{\mathbb{R}}
\newcommand{\C}{\mathbb{C}}

\newcommand{\Proj}{\mathbb{P}}
\newcommand{\Q}{\mathbb{Q}}
\newcommand{\Z}{\mathbb{Z}}
\newcommand{\Hyp}{\mathbb{H}}

\newcommand{\N}{\mathbb{N}}
\newcommand{\D}{\mathbb{D}}
\newcommand{\B}{\mathbb{B}}

\newcommand{\T}{\mathcal{T}}
\newcommand{\RT}{\mathscr{T}}
\newcommand{\RV}{\mathscr{V}}
\newcommand{\RP}{\Pi}
\newcommand{\p}{p}
\newcommand{\fm}{\text{fm}}
\newcommand{\fc}{\text{fc}}
\newcommand{\zc}{\text{zc}}
\newcommand{\bp}{\mathcal{F}}
\newcommand{\rl}{\mathcal{R}}
\newcommand{\U}{\mathcal{U}}

\DeclareMathOperator{\Rat}{Rat}

\DeclareMathOperator{\QH}{\mathcal{QB}}
\DeclareMathOperator{\BP}{\mathcal{B}}
\DeclareMathOperator{\MP}{\mathcal{P}}

\DeclareMathOperator{\PSL}{PSL}

\DeclareMathOperator{\Isom}{Isom}

\DeclareMathOperator{\Res}{Res}

\DeclareMathOperator{\Int}{Int}
\DeclareMathOperator{\chull}{Cvx\, Hull}

\DeclareMathOperator{\proj}{proj}

\numberwithin{figure}{section}

\title[On geometrically finite degenerations II]{On geometrically finite degenerations II: convergence and divergence}
\author{Yusheng Luo}
\address{Dept. of Mathematics \& University of Michigan, Ann Arbor, MI 48109 USA}
\email{yusheng.s.luo@gmail.com}
\date{\today}
\begin{document}

\begin{abstract}
In this paper, we study quasi post-critically finite degenerations for rational maps.
We construct limits for such degenerations as geometrically finite rational maps on a finite tree of Riemann spheres.
We prove the boundedness for such degenerations of hyperbolic rational maps with Sierpinski carpet Julia set and give criteria for the convergence for quasi-Blaschke products $\QH_d$, making progress towards the analogues of Thurston's compactness theorem for acylindrical $3$-manifold and the double limit theorem for quasi-Fuchsian groups in complex dynamics.
In the appendix, we apply such convergence results to show the existence of certain polynomial matings.
\end{abstract}

\maketitle

\setcounter{tocdepth}{1}
\tableofcontents

\section{Introduction}\label{sec:intro}
The study of iterations of rational maps on Riemann sphere $\hat\C$ has been a central topic in dynamics.
Classically, hyperbolic rational maps are easy to analyze.
They form an open and conjecturally dense subset in the moduli space, and a connected component $\mathcal{H}$ is called a {\em hyperbolic component}.
The problem of how the hyperbolic components are positioned has been studied broadly in the literature, especially in low degrees \cite{DH85, BH88, Rees90}.

Since a hyperbolic component $\mathcal{H}$ is not bounded in general, it is interesting and useful to understand 
\begin{enumerate}[label=\arabic*)]
\item when a sequence $[f_n] \in \mathcal{H}$ is bounded; and 
\item how a sequence $[f_n] \in \mathcal{H}$ diverges.
\end{enumerate}
In this paper, we answer these two questions for {\em quasi post-critically finite} degenerations:

\begin{itemize}
\item (Theorem \ref{thm:crm}: Limit as rational maps on a tree of Riemann spheres.) We show a quasi post-critically finite degeneration $[f_n] \in \mathcal{H}$ with connected Julia set has a subsequence converging to a geometrically finite rational map on a finite tree of Riemann spheres; 
\item (Theorem \ref{thm:schcb}: Boundedness of Sierpinski carpet maps.) We show if $\mathcal{H}$ has Sierpinski carpet Julia set, then any quasi post-critically finite degeneration $[f_n] \in \mathcal{H}$ has a convergent subsequence in the moduli space of rational maps;
\item (Theorem \ref{thm:dlm}: A double limit theorem.) We give criteria for convergence of quasi post-critically finite degenerations for quasi-Blaschke products $\QH_d$, i.e., the hyperbolic component of rational maps containing $z^d$.
\end{itemize}

In the prequel \cite{Luo21}, we proved a realization theorem for quasi post-critically finite degenerations.
We used it to classify geometrically finite polynomials on the boundary of the main hyperbolic component $\mathcal{H}_d$. 
The geometrically finite polynomials are conjecturally dense on the boundary.
In some loose sense, quasi post-critically finite degenerations provide a dense set of `directions' for degenerations in a hyperbolic component.

Together with the realization theorem in the prequel, the convergence results in this paper can be applied to show the existence of polynomial matings, or more generally, tunings of rational maps.
In the appendix, we give one such application (see Theorem \ref{thm:mating}).

Our theory is motivated by the parallel developments for Kleinian groups and fits into the well-known Sullivan's dictionary.
\begin{itemize}
\item Degenerating Kleinian groups give isometric group actions on $\R$-trees \cite{MorganShalen84, Bestvina88, Paulin88}.
Limiting branched coverings on $\R$-trees are constructed for arbitrary degenerations of rational maps in \cite{L19p, L19b} (see also \cite{McM09, Kiwi15}).
The limiting dynamics on $\R$-trees are typically quite complicated.
Under our imposed geometric condition, Theorem \ref{thm:crm} asserts that the limiting dynamics is tame.
This allows us to fully analyze the convergence/divergence problem.

\item A compact Kleinian manifold is acylindrical if and only if its limit set is a Sierpinski carpet.
Thurston showed that deformation space of acylindrical $3$-manifold is bounded \cite{Thurston86}.
Base on this, McMullen conjectured that a hyperbolic component with Sierpinski carpet Julia set is bounded \cite{McM95}.
Theorem \ref{thm:schcb} verifies this for quasi post-critically finite degenerations, and is one step towards the conjecture.

\item Thurston's double limit theorem gives a sufficient condition for a sequence of quasi-Fuchsian groups to have a convergent subsequence \cite{Thurston98}.
Theorem \ref{thm:dlm} gives an analogue in complex dynamics for quasi post-critically finite degenerations.
\end{itemize}


\subsection*{Statement of results}
Let $f: \hat\C \longrightarrow \hat\C$ be a rational map of degree $d\geq 2$.
It is said to be {\em hyperbolic} if all critical points converge to attracting periodic cycles under iteration.
For our purposes, it is convenient to mark the rational map by an ordered list of $d+1$ (not necessarily distinct) fixed points (see \cite{Milnor12}).
We denote $\Rat_{d, \fm}$ as the space of all fixedpoint marked rational maps of degree $d$.
The group of M\"obius transformation $\PSL_2(\C)$ acts naturally on $\Rat_{d,\fm}$ and we define the {\em moduli space} 
$$
\mathcal{M}_{d,\fm} = \Rat_{d,\fm}/\PSL_2(\C).
$$

Let $\mathcal{H} \subseteq \mathcal{M}_{d,\fm}$ be a hyperbolic component with connected Julia set.
Let $\mathcal{U} = \bigcup_{s \in |\mathcal{S}|} \U_s$ be the union of critical and post-critical Fatou components for $[f] \in \mathcal{H}$, where $|\mathcal{S}|$ is the index set.
We define a `metric' $d_{\U}$ on $\mathcal{U}$ by
\begin{itemize}
\item $d_{\U}(x,y) = \infty$ if $x, y$ are in different components of $\U$;
\item $d_{\U}(x,y) = d_{\U_s}(x, y)$ if $x, y$ are in the same component $\U_s \subseteq \U$,
\end{itemize}
where $d_{\U_s}$ is the hyperbolic metric on $\U_s$.

A sequence $[f_n] \in \mathcal{H}$ is said to be {\em ($K$-)quasi post-critically finite} if we can label the critical points as $c_{1,n},..., c_{2d-2, n}$, and for any sequence $c_{i,n}$, there exist quasi pre-period $l_i$ and quasi period $q_i$ such that
$$
d_{\U_n}(f_n^{l_i}(c_{i,n}), f_n^{l_i+q_i}(c_{i,n})) \leq K,
$$
where $\U_n= \bigcup_{s \in |\mathcal{S}|} \U_{s,n}$ is the corresponding set for $[f_n]$.

We say a sequence $[f_n] \in \mathcal{H}$ is a {\em ($K$-)quasi post-critically finite degeneration} if $[f_n]$ is ($K$-)quasi post-critically finite and it leaves every compact subset of $\mathcal{H}$.
We also say $[f_n]$ is {\em diverging in $\mathcal{M}_{d, \fm}$}, or simply {\em diverging} if $[f_n]$ leaves every compact subset of $\mathcal{M}_{d, \fm}$.
Note that a degenerating sequence $[f_n] \in \mathcal{H}$ may not be diverging.

We define a {\em tree of Riemann spheres} $(\RT, \hat\C^\RV)$ as a finite tree $\RT$ with vertex set $\RV$, a disjoint union of Riemann spheres $\hat \C^\RV := \bigcup_{a\in \RV}\hat \C_{a}$, together with markings $\xi_a: T_a\RT \xhookrightarrow{} \hat\C_a$ for $a\in \RV$.

A {\em rational map} $(F, R)$ on a tree of Riemann spheres is a map 
$$
F: (\RT, \RV) \longrightarrow (\RT, \RV)
$$ 
that is injective on each edge and a union of maps $R:= \bigcup_{a\in \RV} R_a$ so that
\begin{itemize}
\item $R_a: \hat\C_a \longrightarrow \hat \C_{F(a)}$ is a rational map of degree $\geq 1$;
\item $R_a \circ \xi_a = \xi_{F(a)} \circ DF_a$,
\end{itemize}
where $DF_a: T_a\RT \longrightarrow T_{F(a)}\RT$ is the tangent map.

Define $\Xi_a := \xi_a(T_a\RT)$ and $\Xi := \cup_{a\in \RV} \Xi_a$ as the {\em singular set}.
The rational map $(F, R)$ is said to have degree $d$ if $R$ has $2d-2$ critical points in $\hat\C^\RV-\Xi$.
We also define the local degree of $a$ by $\deg(a):=\deg(R_a)$.

A sequence $f_n$ of rational maps is said to {\em converge} to $(F, R)$ on $(\RT, \hat\C^\RV)$ if there exist rescalings $A_{a,n} \in \PSL_2(\C)$ representing $a \in \RV$ such that 
$$
A_{F(a),n}^{-1} \circ f_n \circ A_{a,n}(z) \to R_a(z)
$$ 
compactly on $\hat \C_a- \Xi_a$ (See Definition \ref{defn:trs}).

The Fatou set $\Omega(R) \subseteq \hat\C^\RV$ is defined naturally as the largest open set for which the iterations $\{R^n|_{\Omega(R)}: n \geq 1\}$ form a normal family. The Julia set $\mathcal{J}(R)$ is the complement of the Fatou set.
We say $(F, R)$ is {\em geometrically finite} if every critical point of $R$ in the Julia set has a finite orbit.

Our first result shows
\begin{theorem}\label{thm:crm}
Let $\mathcal{H} \subseteq \mathcal{M}_{d, \fm}$ be a hyperbolic component with connected Julia set.
Let $[f_n] \in \mathcal{H}$ be a quasi post-critically finite degeneration.
Then after passing to a subsequence, $[f_n]$ converges to a degree $d$ geometrically finite rational map $(F, R)$ on a tree of Riemann spheres $(\RT, \hat\C^\RV)$.

Moreover, $[f_n]$ converges in $\mathcal{M}_{d, \fm}$ if and only if there exists a fixed point $a\in \RV$ of local degree $d$.
\end{theorem}

The proof is based on the idea that a quasi post-critically finite degeneration $[f_n]$ gives abundant {\em rescaling limits} sketched as below (cf. \cite{Kiwi15}).

In \S \ref{sec:DBS}, we show that after passing to a subsequence, a quasi post-critically finite degeneration $[f_n]$ gives a finite union of {\em quasi-invariant trees} $\T_{s,n} \subseteq \U_{s,n} \subseteq \U_n$ marked by
$$
\psi_n: \mathcal{T} := \bigcup_{s\in |\mathcal{S}|} \mathcal{T}_s \longrightarrow \mathcal{T}_n := \bigcup_{s\in |\mathcal{S}|} \mathcal{T}_{s,n}
$$
with vertex set $\mathcal{V} = \bigcup_{s\in |\mathcal{S}|} \mathcal{V}_s$ and induced dynamics
$$
\widetilde F: \mathcal{T} \longrightarrow \mathcal{T}.
$$
Each vertex $v\in \mathcal{V}_s \subseteq \mathcal{V}$ thus corresponds to a sequence of pointed disks $(\U_{s,n}, v_n)$.
After appropriate {\em rescaling} by $A_{v,n} \in \PSL_2(\C)$\footnote{$A_{v,n}$ is well-defined up to pre-copmosing with a bounded sequence in $\PSL_2(\C)$.}, the pointed disk $A_{v,n}^{-1}((\U_{s,n}, v_n))$ converges in Carath\'edoroy topology to a pointed disk $(\U_v, 0)$.
Since $\mathcal{T}_{n}$ is quasi-invariant under $f_n$, the sequence
$$
A_{\widetilde F(v),n}^{-1} \circ f_n \circ A_{v,n}
$$
converges compactly to a proper map between $\U_v$ and $\U_{\widetilde F(v)}$, and thus converges algebraically\footnote{See \S \ref{sec:DBS} for definition of algebraic convergence.} to a rational map of degree $\geq 1$.

Two vertices $v,w$ are said to be equivalent if $A_{v,n}^{-1} \circ A_{w,n}$ is a bounded sequence in $\PSL_2(\C)$.
Denote the quotient of this equivalence relation by $\RP = \mathcal{V} /\sim$ with induced map $F: \RP \rightarrow \RP$.

By picking a representative rescaling $A_{a,n} \in \PSL_2(\C)$ for each $a\in \RP$, we get a sequence $x_{a,n} = A_{a,n}({\bf 0}) \in \Hyp^3$, where we have identified $\Hyp^3$ as the unit ball in $\R^3$ with the origin ${\bf 0}$ and $\PSL_2(\C) \cong \Isom(\Hyp^3)$.

The finite tree $(\RT, \RV)$ is constructed as the `spine' for the degenerating hyperbolic polyhedra in $\Hyp^3$
$$
\chull(\{x_{a,n}: a\in \RP\}) \subseteq \Hyp^3.
$$
The map $F: \RP \rightarrow \RP$ induces a {\em rescaling tree map} $F: (\RT, \RV) \longrightarrow (\RT, \RV)$.
Theorem \ref{thm:crm} is proved by analyzing the map $F$ and {\em rescaling rational maps} 
$$
R_a:=\lim_{n\to\infty} A_{F(a),n}^{-1} \circ f_n \circ A_{a,n}.
$$ 

\begin{figure}[ht]
  \centering
  \resizebox{0.7\linewidth}{!}{
    \def\svgwidth{\columnwidth}
    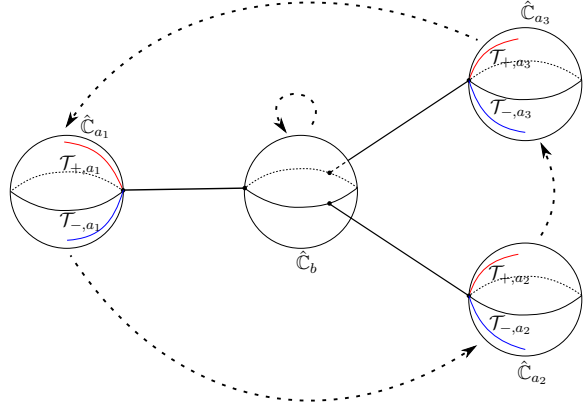

  }
  \caption{An illustration of a quasi post-critical finite degeneration in $\QH_2$ where $\bp_\pm$ degenerates to the fat rabbit and fat co-rabbit radially and at the same rate.
  The quasi-invariant trees are tripods (see Figure \ref{fig:R} upper-left).}
  \label{fig:DB}
\end{figure}

\begin{figure}[ht]
  \centering
  \resizebox{0.8\linewidth}{!}{
    \def\svgwidth{\columnwidth}
\begingroup%
  \makeatletter%
  \providecommand\color[2][]{%
    \errmessage{(Inkscape) Color is used for the text in Inkscape, but the package 'color.sty' is not loaded}%
    \renewcommand\color[2][]{}%
  }%
  \providecommand\transparent[1]{%
    \errmessage{(Inkscape) Transparency is used (non-zero) for the text in Inkscape, but the package 'transparent.sty' is not loaded}%
    \renewcommand\transparent[1]{}%
  }%
  \providecommand\rotatebox[2]{#2}%
  \newcommand*\fsize{\dimexpr\f@size pt\relax}%
  \newcommand*\lineheight[1]{\fontsize{\fsize}{#1\fsize}\selectfont}%
  \ifx\svgwidth\undefined%
    \setlength{\unitlength}{388.5bp}%
    \ifx\svgscale\undefined%
      \relax%
    \else%
      \setlength{\unitlength}{\unitlength * \real{\svgscale}}%
    \fi%
  \else%
    \setlength{\unitlength}{\svgwidth}%
  \fi%
  \global\let\svgwidth\undefined%
  \global\let\svgscale\undefined%
  \makeatother%
  \begin{picture}(1,0.95559846)%
    \lineheight{1}%
    \setlength\tabcolsep{0pt}%
    \put(0,0){\includegraphics[width=\unitlength,page=1]{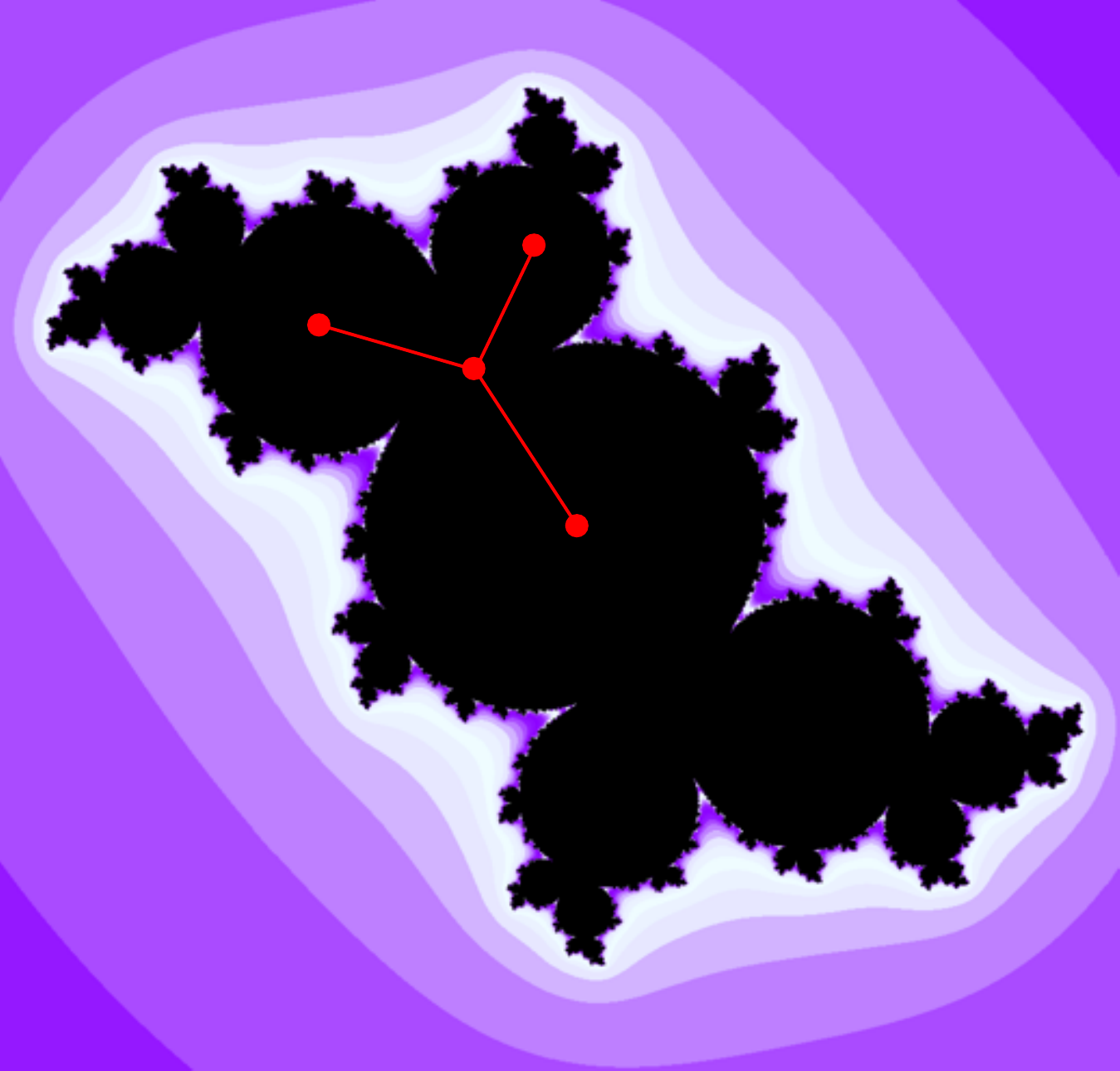}}%
    \put(0.53489906,0.48224011){\color[rgb]{1,0,0}\makebox(0,0)[lt]{\lineheight{1.25}\smash{\begin{tabular}[t]{l}$c_n$\end{tabular}}}}%
    \put(0.44939229,0.62855636){\color[rgb]{1,0,0}\makebox(0,0)[lt]{\lineheight{1.25}\smash{\begin{tabular}[t]{l}$p_n$\end{tabular}}}}%
    \put(0,0){\includegraphics[width=\unitlength,page=2]{R.pdf}}%
    \put(0.78243238,0.06479166){\color[rgb]{0,0,0}\makebox(0,0)[lt]{\lineheight{1.25}\smash{\begin{tabular}[t]{l}$\deg(c_n)=2$\end{tabular}}}}%
  \end{picture}%
\endgroup%

    \def\svgwidth{\columnwidth}
\begingroup%
  \makeatletter%
  \providecommand\color[2][]{%
    \errmessage{(Inkscape) Color is used for the text in Inkscape, but the package 'color.sty' is not loaded}%
    \renewcommand\color[2][]{}%
  }%
  \providecommand\transparent[1]{%
    \errmessage{(Inkscape) Transparency is used (non-zero) for the text in Inkscape, but the package 'transparent.sty' is not loaded}%
    \renewcommand\transparent[1]{}%
  }%
  \providecommand\rotatebox[2]{#2}%
  \newcommand*\fsize{\dimexpr\f@size pt\relax}%
  \newcommand*\lineheight[1]{\fontsize{\fsize}{#1\fsize}\selectfont}%
  \ifx\svgwidth\undefined%
    \setlength{\unitlength}{388.5bp}%
    \ifx\svgscale\undefined%
      \relax%
    \else%
      \setlength{\unitlength}{\unitlength * \real{\svgscale}}%
    \fi%
  \else%
    \setlength{\unitlength}{\svgwidth}%
  \fi%
  \global\let\svgwidth\undefined%
  \global\let\svgscale\undefined%
  \makeatother%
  \begin{picture}(1,0.95559846)%
    \lineheight{1}%
    \setlength\tabcolsep{0pt}%
    \put(0,0){\includegraphics[width=\unitlength,page=1]{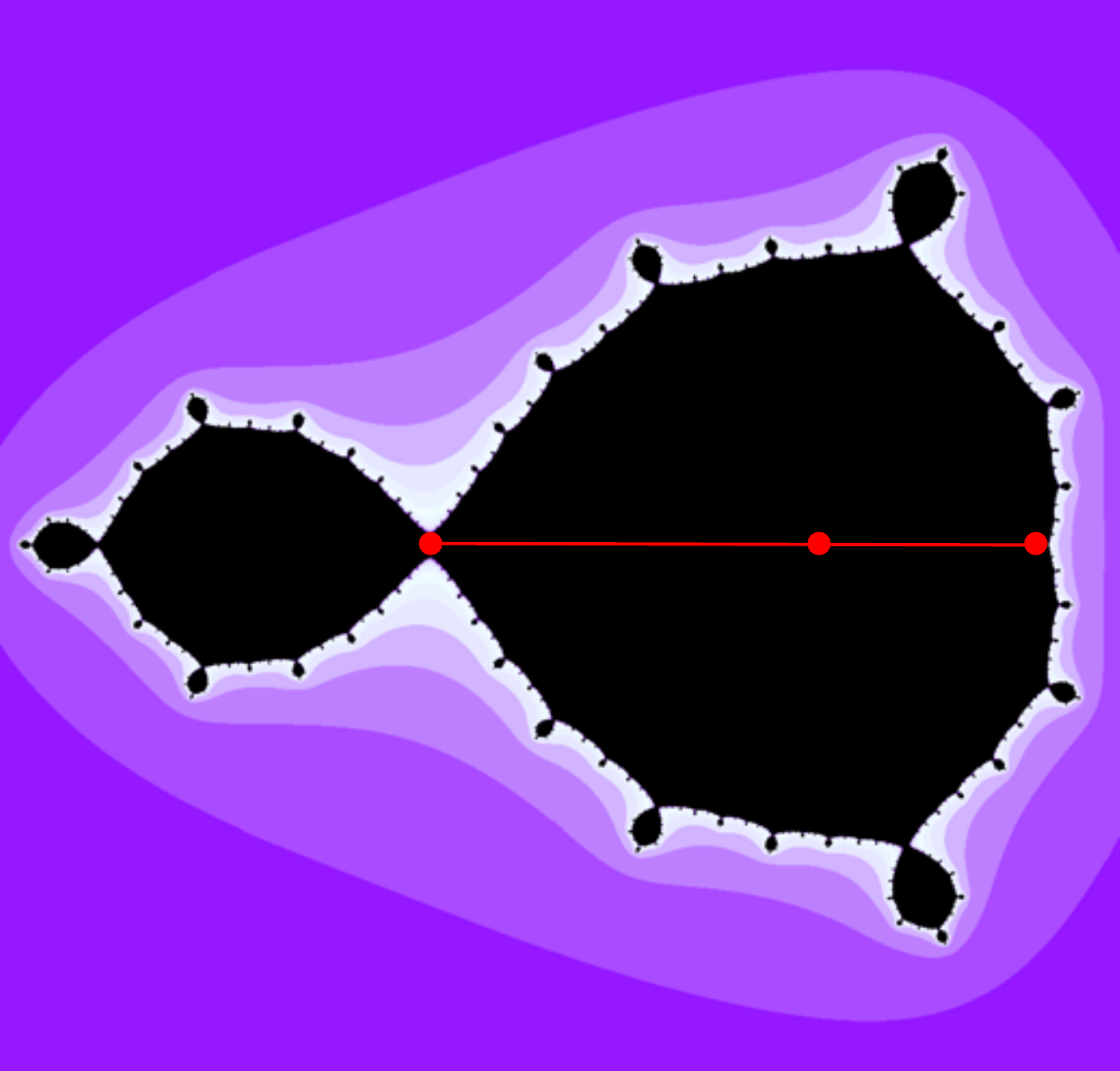}}%
    \put(0.72661393,0.43626924){\color[rgb]{1,0,0}\makebox(0,0)[lt]{\lineheight{1.25}\smash{\begin{tabular}[t]{l}$p_n$\end{tabular}}}}%
    \put(0.40145289,0.44004794){\color[rgb]{1,0,0}\makebox(0,0)[lt]{\lineheight{1.25}\smash{\begin{tabular}[t]{l}$c_n$\end{tabular}}}}%
    \put(0,0){\includegraphics[width=\unitlength,page=2]{G.pdf}}%
    \put(0.6932368,0.06260605){\color[rgb]{0,0,0}\makebox(0,0)[lt]{\lineheight{1.25}\smash{\begin{tabular}[t]{l}$\deg(p_n) = \deg(c_n)=2$\end{tabular}}}}%
  \end{picture}%
\endgroup%

  }
  \resizebox{0.8\linewidth}{!}{
    \def\svgwidth{\columnwidth}
\begingroup%
  \makeatletter%
  \providecommand\color[2][]{%
    \errmessage{(Inkscape) Color is used for the text in Inkscape, but the package 'color.sty' is not loaded}%
    \renewcommand\color[2][]{}%
  }%
  \providecommand\transparent[1]{%
    \errmessage{(Inkscape) Transparency is used (non-zero) for the text in Inkscape, but the package 'transparent.sty' is not loaded}%
    \renewcommand\transparent[1]{}%
  }%
  \providecommand\rotatebox[2]{#2}%
  \newcommand*\fsize{\dimexpr\f@size pt\relax}%
  \newcommand*\lineheight[1]{\fontsize{\fsize}{#1\fsize}\selectfont}%
  \ifx\svgwidth\undefined%
    \setlength{\unitlength}{388.5bp}%
    \ifx\svgscale\undefined%
      \relax%
    \else%
      \setlength{\unitlength}{\unitlength * \real{\svgscale}}%
    \fi%
  \else%
    \setlength{\unitlength}{\svgwidth}%
  \fi%
  \global\let\svgwidth\undefined%
  \global\let\svgscale\undefined%
  \makeatother%
  \begin{picture}(1,0.95559846)%
    \lineheight{1}%
    \setlength\tabcolsep{0pt}%
    \put(0,0){\includegraphics[width=\unitlength,page=1]{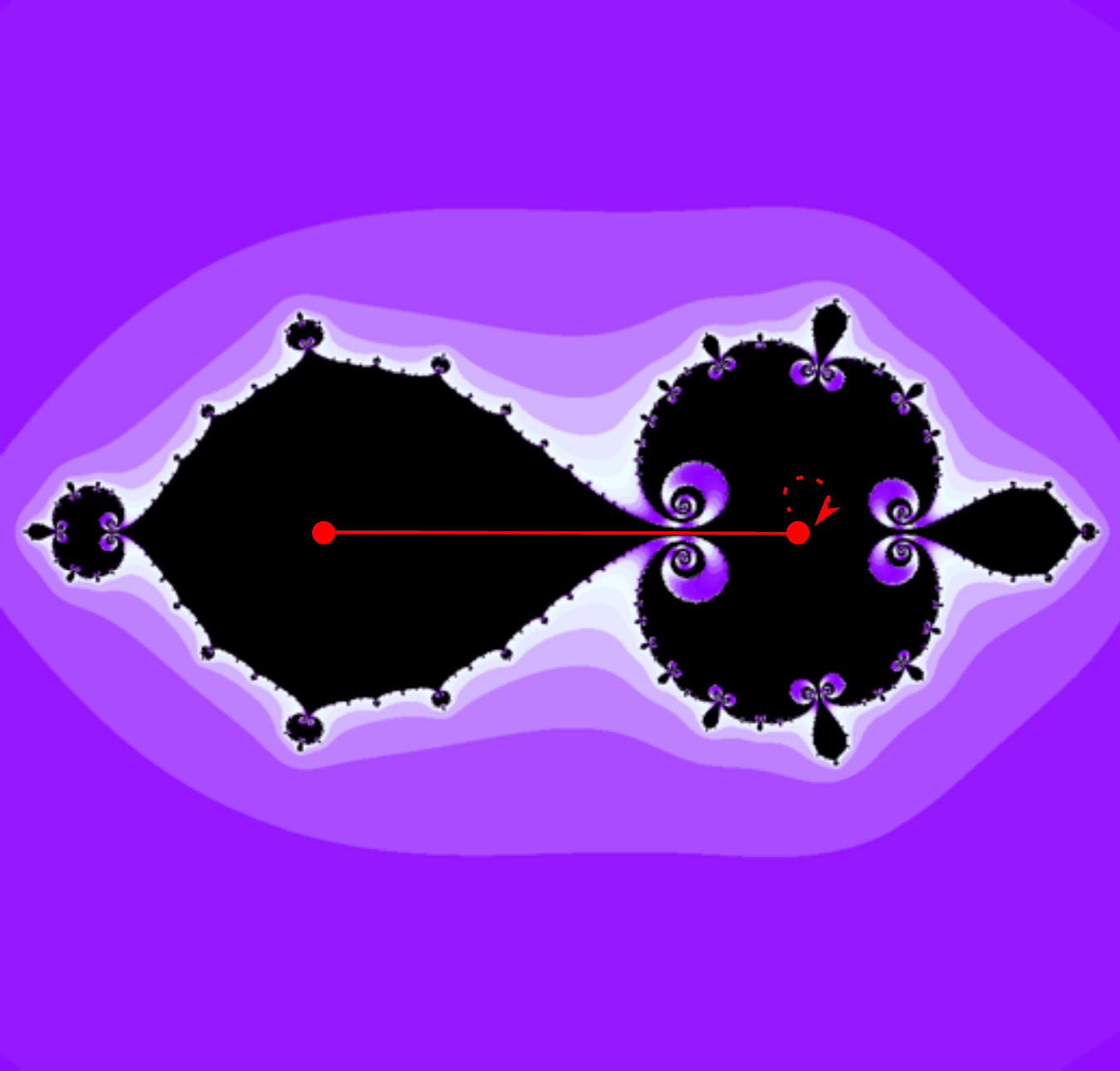}}%
    \put(0.28836901,0.44570684){\color[rgb]{1,0,0}\makebox(0,0)[lt]{\lineheight{1.25}\smash{\begin{tabular}[t]{l}$p_n$\end{tabular}}}}%
    \put(0.68801972,0.08010952){\color[rgb]{0,0,0}\makebox(0,0)[lt]{\lineheight{1.25}\smash{\begin{tabular}[t]{l}$\deg(p_n) = \deg(c_n)=2$\end{tabular}}}}%
    \put(0.70269265,0.44289106){\color[rgb]{1,0,0}\makebox(0,0)[lt]{\lineheight{1.25}\smash{\begin{tabular}[t]{l}$c_n$\end{tabular}}}}%
  \end{picture}%
\endgroup%

    \def\svgwidth{\columnwidth}
\begingroup%
  \makeatletter%
  \providecommand\color[2][]{%
    \errmessage{(Inkscape) Color is used for the text in Inkscape, but the package 'color.sty' is not loaded}%
    \renewcommand\color[2][]{}%
  }%
  \providecommand\transparent[1]{%
    \errmessage{(Inkscape) Transparency is used (non-zero) for the text in Inkscape, but the package 'transparent.sty' is not loaded}%
    \renewcommand\transparent[1]{}%
  }%
  \providecommand\rotatebox[2]{#2}%
  \newcommand*\fsize{\dimexpr\f@size pt\relax}%
  \newcommand*\lineheight[1]{\fontsize{\fsize}{#1\fsize}\selectfont}%
  \ifx\svgwidth\undefined%
    \setlength{\unitlength}{388.5bp}%
    \ifx\svgscale\undefined%
      \relax%
    \else%
      \setlength{\unitlength}{\unitlength * \real{\svgscale}}%
    \fi%
  \else%
    \setlength{\unitlength}{\svgwidth}%
  \fi%
  \global\let\svgwidth\undefined%
  \global\let\svgscale\undefined%
  \makeatother%
  \begin{picture}(1,0.95559846)%
    \lineheight{1}%
    \setlength\tabcolsep{0pt}%
    \put(0,0){\includegraphics[width=\unitlength,page=1]{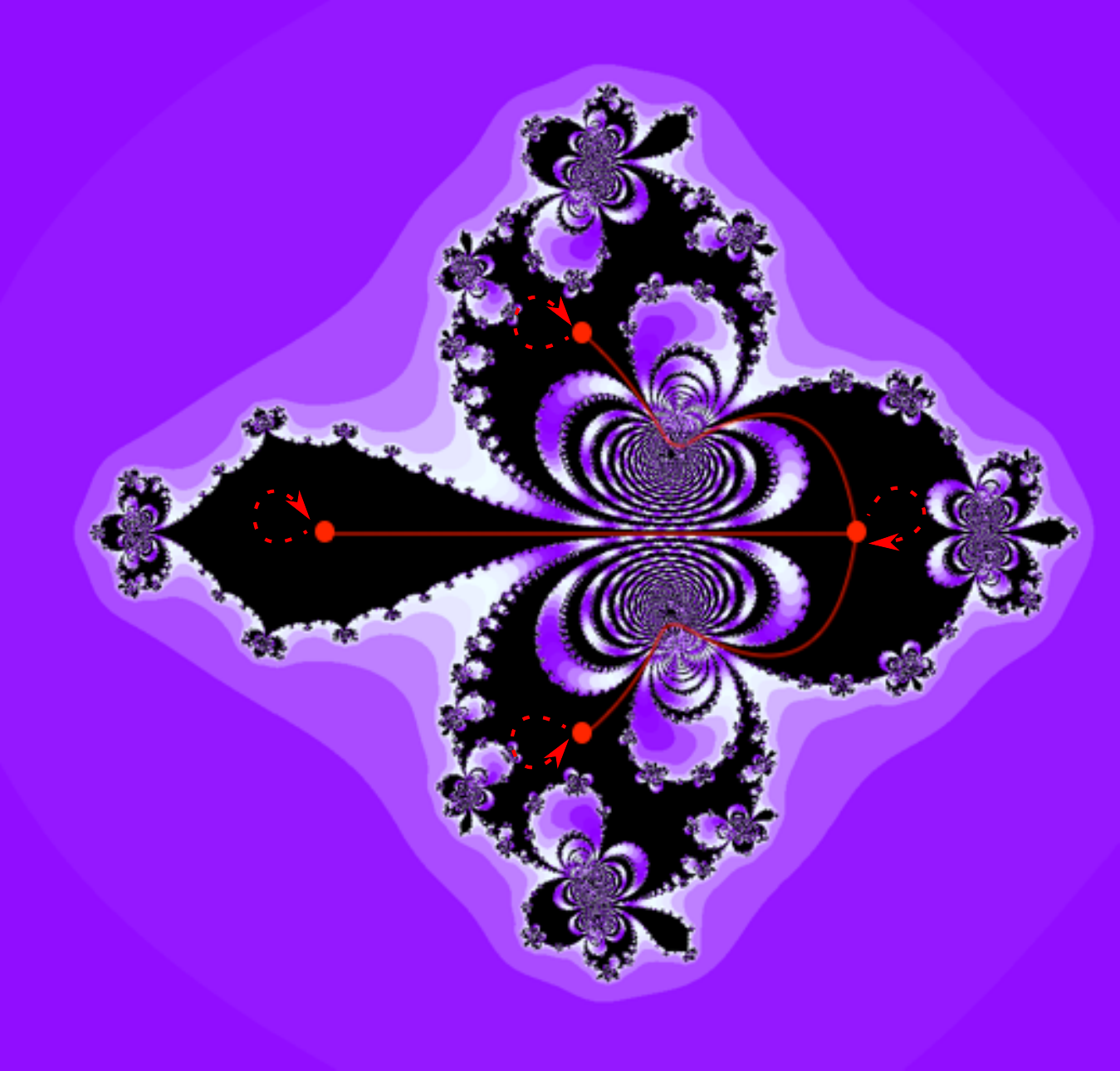}}%
    \put(0.29128956,0.45348029){\color[rgb]{1,0,0}\makebox(0,0)[lt]{\lineheight{1.25}\smash{\begin{tabular}[t]{l}$p_n$\end{tabular}}}}%
    \put(0.68138877,0.07994368){\color[rgb]{0,0,0}\makebox(0,0)[lt]{\lineheight{1.25}\smash{\begin{tabular}[t]{l}$\deg = 2$ for all vertices\end{tabular}}}}%
  \end{picture}%
\endgroup%

  }
  \caption{Illustrations of quasi-invariant trees for quasi post-critically finite degenerations.}
  \label{fig:R}
\end{figure}

\subsection*{Sierpinski carpet rational maps}
A closed set $K \subseteq \hat\C$ is called a {\em Sierpinski carpet} if it is obtained from the sphere by removing a countable dense set of open disks, bounded by disjoint Jordan curves whose diameters tend to zero.
Applying Theorem \ref{thm:crm}, we prove


\begin{figure}[ht]
  \centering
  \includegraphics[width=0.6\textwidth]{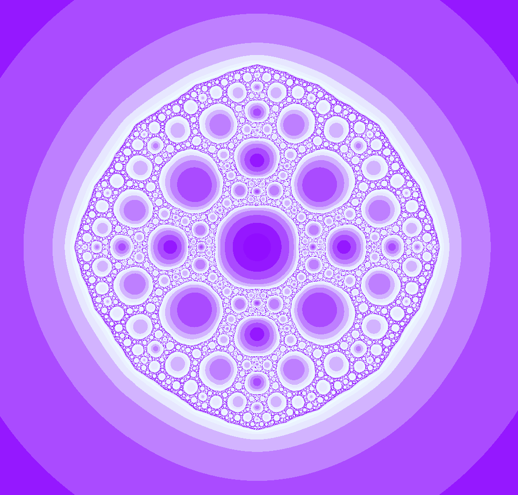}
  \caption{The rational map $f(z) = z^2+\frac{1}{16z^2}$ has a Sierpinski carpet Julia set.}
  \label{fig:SC}
\end{figure}

\begin{theorem}\label{thm:schcb}
Let $\mathcal{H}\subseteq \mathcal{M}_{d, \fm}$ be a hyperbolic component with Sierpinski carpet Julia set.
Then any quasi post-critically finite degeneration $[f_n] \in \mathcal{H}$ has a convergent subsequence in $\mathcal{M}_{d, \fm}$.
\end{theorem}

The proof is by contradiction.
By Theorem \ref{thm:crm}, $[f_n]$ converges to a rational map on $(\RT, \hat\C^\RV)$.
We first show the limits of quasi-invariant trees $\mathcal{T}_{s,a} = \lim A_{a, n}^{-1} (\T_{s,n})$ are disjoint in $\hat\C^\RV$.
This allows us to construct a curve system $\Sigma_n$ in the complement of critical and post-critical Fatou components for sufficiently large $n$, where curves in $\Sigma_n$ are in correspondence with edges of an invariant subtree $\RT^r \subseteq \RT$.

If $[f_n]$ is diverging, then $\RT^r$ is not trivial.
We can record the dynamics $F$ on the edges $E_1,..., E_k$ of $\RT^r$ by two matrices
\begin{itemize}
\item (Markov matrix): $M_{i,j} = \begin{cases} 1 &\mbox{if } E_i \subseteq F(E_j) \\ 
0 & \mbox{otherwise } \end{cases}$
\item (Degree matrix): $D_{i,j} = \begin{cases} \delta(E_i) &\mbox{if } i = j \\ 
0 & \mbox{otherwise } \end{cases}$
\end{itemize}
where $\delta(E_i)$ is the local degree at an edge.
We show that there exists a non-negative vector $\vec{v} \neq \vec{0}$ with $M\vec{v} = D\vec{v}$.

As matrices with non-negative entries, we show $D^{-1}M$ is no bigger than the Thurston's matrix associated to $\Sigma_n$ for the hyperbolic rational map $f_n$. 
So the spectral radius of the Thurston's matrix is greater or equal to $1$, giving a contradiction.

\subsection*{Quasi-Blaschke products}
A hyperbolic rational map $f$ is said to be a {\em quasi-Blaschke product} if $J(f)$ is a Jordan curve and $f$ fixes the two Fatou components.
Let $\QH_d \subseteq \mathcal{M}_{d,fm}$ be the hyperbolic component consisting of marked quasi-Blaschke products. 
A marked quasi-Blaschke product $f$ is obtained by gluing together a pair of marked Blaschke products $\bp_+, \bp_-\in \BP_d$ using their markings on $\mathbb{S}^1$.

Many quasi post-critically finite degenerations are diverging for quasi-Blaschke products (see Figure \ref{fig:DB}). We will give criteria for convergence using the dual laminations, defined as follows.

A {\em lamination} $\mathcal{L}$ is a family of disjoint hyperbolic geodesics in $\D$ together with the two end points in $\mathbb{S}^1\cong \R/\Z$, whose union $|\mathcal{L}| :=\bigcup\mathcal{L}$ is closed.
We define two laminations are {\em parallel} if there exists
a simple closed loop in the union $|\mathcal{L}| \cup \overline{|\mathcal{L}'|}$ containing leaves in both laminations, where $\overline{|\mathcal{L}'|}$ is the image of $|\mathcal{L}'|$ under $z \mapsto \frac{1}{z}$.
We show the two quasi-invariant trees $\mathcal{T}_\pm$ give a pair of {\em dual laminations} $\mathcal{L}_\pm$, and prove:
\begin{theorem}\label{thm:dlm}
Let $[f_n] \in \QH_d \subseteq \mathcal{M}_{d,\fm}$ be a quasi post-critically finite degeneration.
Let $\mathcal{L}_{\pm}$ be the corresponding dual laminations.
If $\mathcal{L}_{\pm}$ are non-parallel, then $[f_n]$ has a convergent subsequence in $\mathcal{M}_{d, \fm}$.
\end{theorem}

We now sketch the idea of the proof.
Since the fixed points of $f_n$ are marked, the Julia set is marked $\eta_n: \R/\Z \longrightarrow \mathcal{J}(f_n)$.
By Theorem \ref{thm:crm}, $[f_n]$ converges to a rational map on $(\RT, \hat\C^\RV)$.
Suppose $[f_n]$ diverges.
We show the partition of the two limiting quasi-invariant graphs $\T_{\pm}$ into $\bigcup_{a \in \RV} \T_{\pm, a}$, where $\T_{\pm, a} = \T_ {\pm} \cap \hat\C_a$ gives two non-trivial partitions of angles 
$$
\R/\Z = \bigcup_{a \in \RV} I_{\pm, a}.
$$
We show that for all but finitely many rational angles $t\in \pm I_{\pm, a}$, the rescaled limit $A_{a,n}^{-1}(\eta_n(t))$ converges to a point in $\hat \C_a - \Xi_a$.
This implies the two partitions are conjugate, i.e., $I_{+, a} = -I_{-,a}$. We show this implies $\mathcal{L}_{\pm}$ are parallel (see Figure \ref{fig:Par}).

\subsection*{Applications on polynomial matings}
These convergence results can be coupled with the realization theorem proved in \cite{Luo21} to show the existence of polynomial matings or tunings.
In the Appendix, we give one such application.

Given a geometrically finite monic and centered polynomial $P \in \MP_d$ with connected Julia set, there exists a unique post-critically finite monic and centered polynomial $\hat P \in \MP_d$ with conjugate dynamics on the Julia sets with compatible markings \cite{Haissinsky00}.
The dynamics of $\hat P$ is described combinatorially by the dynamics on its Hubbard tree $\hat P: H \rightarrow H$.
We say $P$ has simplicial Hubbard tree if the corresponding dynamics $\hat P: H \rightarrow H$ is simplicial: there exists a finite simplicial structure on $H$ so that $\hat P$ sends an edge of $H$ to an edge of $H$. We show
{
\renewcommand{\thetheorem}{\ref{thm:mating}}
\begin{theorem}
Let $P_\pm \in \MP_d$ be two geometrically finite polynomials, where $P_-$ has simplicial Hubbard tree.
Then $P_\pm$ are mateable if and only if their laminations $\mathcal{L}_\pm$ are non-parallel.
\end{theorem}
\addtocounter{theorem}{-1}
}

The idea of the proof is to start with a geometrically finite polynomial $P_+$ and perform appropriate quasi post-critically finite degeneration on the unbounded Fatou component.
The realization theorem proved in \cite{Luo21} and convergence results in this paper give that the limit is the desired mating.

We do not insist to give the most general results.
The techniques can potentially be applied to the case where $P_+$ is not geometrically finite.
By inductively performing quasi post-critically finite degenerations on Fatou components, a similar technique can also be applied to the case where $P_-$ has zero core entropy on its Hubbard tree.

Since two geometrically finite polynomials are mateable if and only if the corresponding post-critically finite polynomials are mateable \cite{HT04}, the above result can also be proved by applying Thurston's characterization theorem of rational maps (cf. Theorem E in \cite{BD18} and arc intersecting obstruction theorem in \cite{PT98} and \cite{Park21}).
It is worth investigating if it is possible to prove similar results using Dylan Thurston's positive criterion \cite{Thurston20} (see \cite{DPWY20} for a degree $2$ application).

\subsection*{Notes and discussions}
Degenerations of Fatou components and their convergence/divergence properties have been studied extensively in terms of elementary moves like pinching and spinning \cite{Makienko00, TanLei02, HT04, PT04, CT18}.
When the Julia set is connected, such elementary operations provide examples of quasi post-critically finite degenerations.
On the other hand, many quasi post-critically finite degenerations cannot be constructed by these elementary operations.
Thus, the convergence/divergence results in the paper generalize those discussions for pinching and spinning.

The topology of various hyperbolic components is studied in \cite{BDK91, Milnor12, WY17}.
The hyperbolic components for quadratic rational maps are studied in details in \cite{Rees90, Rees92, Rees95}.
Results on unboundedness of hyperbolic components were obtained in \cite{Makienko00, TanLei02}. 
On the other hand, showing a hyperbolic component is bounded is generally hard, and is only known in some special cases \cite{Epstein00, NP19, NP20}.

Rescaling limits are first defined and formalized in \cite{Kiwi15}.
Predicting where a rescaling limit can occur is usually not easy.
The quasi post-critically finite degenerations provide dynamically meaningful sequences of normalizations, allowing us to construct all non-trivial rescaling limits (see also \cite{L20}).

Maps on tree of spheres appeared naturally in various setting of conformal dynamics \cite{Pilgrim03, Selinger12, Koch13, HK14, DeMF14, Arfeux17, BD20}.
The decomposition of tree of spheres is closely related to the decomposition along some multicurve giving Thurston's obstructions (see \S \ref{sec:qpcfd} and \S \ref{sec:cscrm}).

The finite rescaling trees $\RT$ are related to the construction using the barycentric extensions of rational maps as in \cite{L19p, L19b}, and are analogues of isometric group actions on $\R$-trees for Kleinian groups \cite{MorganShalen84, Bestvina88, Paulin88}.
Since the length of different edges in $\RT$ may grow to infinity at different rates, the standard construction by `rescaling the hyperbolic metric' as in \cite{L19p} may only reveal a quotient tree of $\RT$.

Matings of polynomials have been studied extensively in the literature, see the discussion in \cite{BEKMPRL12, PM12}.
Various different methods were used to show the existence of the matings \cite{TanLei92, Milnor93, PT98, Shishikura00, YZ00, HT04, BD18, DPWY20}. 
We hope the techniques developed in this paper can be combined with other methods to give more examples of matings.

\subsection*{Outline of the paper}
We will review the theory of Blaschke model spaces $\BP^\mathcal{S}$ for hyperbolic components developed by Milnor in \S \ref{sec:bms}, and study the quasi post-critically finite degenerations for $\BP^\mathcal{S}$ in \S \ref{sec:DBS}.
Quasi post-critically finite degenerations of rational maps are discussed in \S \ref{sec:qpcfd} where Theorem \ref{thm:crm} is proved.
We study the convergence of quasi-invariant forests in \S \ref{sec: cqif}.
The convergence results Theorem \ref{thm:schcb} and Theorem \ref{thm:dlm} are proved in \S \ref {sec:cscrm} and \S \ref{sec:ptq} respectively.
Finally, the application on polynomial matings (Theorem \ref{thm:mating}) is proved in the appendix.

\subsection*{Acknowledgement}
The author thanks Dima Dudko, Curt McMullen, Sabya Mukherjee, and Kevin Pilgrim for useful suggestions and discussions on this problem.
The author gratefully thanks the anonymous reviewers for valuable comments and suggestions.

\section{Blaschke model spaces $\BP^\mathcal{S}$}\label{sec:bms}
In this section, we review the theory of using Blaschke products to model hyperbolic components developed by Milnor in \cite{Milnor12}.

A {\em hyperbolic mapping scheme} $\mathcal{S} = (|\mathcal{S}|, \Phi, \delta)$ (or briefly a {\em mapping scheme}) consists of finite set $|\mathcal{S}|$ of `vertices', together with a map $\Phi = \Phi_{\mathcal{S}}: |\mathcal{S}| \longrightarrow |\mathcal{S}|$, and a degree function $\delta : |\mathcal{S}| \longrightarrow \Z_{\geq 1}$, satisfying two conditions:
\begin{itemize}
\item (Minimality.) Any vertex of degree $1$ is the iterated forward image of some vertex of degree $\geq 2$;
\item (Hyperbolicity.) Every periodic orbit under $\Phi$ contains at least one vertex of degree $\geq 2$.
\end{itemize}
We define the degree of the scheme as
$\deg(\mathcal{S}) = 1+ \sum_{s\in |\mathcal{S}|} (\delta(s)-1)$.

A proper holomorphic map $f: \D\longrightarrow \D$ of degree $d\geq 1$. 
It can be uniquely written as a {\em Blaschke product}
$$
f(z) = e^{i\theta}\prod_{i=0}^d \frac{z-a_i}{1-\overline{a_i}z},
$$
where $|a_i|<1$.
We define a Blaschke product $f$ is 
\begin{itemize}
\item {\bf \em 1-anchored} if $f(1) = 1$;
\item {\bf \em fixed point centered} if $f(0) = 0$;
\item {\bf \em zeros centered} if the sum
$$
a_1+... + a_d
$$
of the points of $f^{-1}(0)$ (counted with multiplicity) is equal to $0$.
\end{itemize}
Note that if $d=1$, then the only Blaschke product which is 1-anchored and centered in either sense is the identity map.
We define $\BP_{d, \fc}$ and $\BP_{d, \zc}$ as the space of all 1-anchored Blaschke products of degree $d$ which are respectively fixed point centered or zeros centered.
When $d=1$, $\BP_{1, \fc} = \BP_{1, \zc}$ consist of only the identity map.

There is another natural normalization of Blaschke products (see \cite{McM09}). Let
$$
\BP_d:=\{f(z) = z\prod_{i=1}^{d-1} \frac{z-a_i}{1-\overline{a_i}z}: |a_i| < 1\}.
$$
Let $p_d(z) = z^d$.
For each $f\in \BP_d$, there exists a unique quasisymmetric homeomorphism $\eta_f:\mathbb{S}^1 \longrightarrow \mathbb{S}^1$ with
\begin{enumerate}
\item $\eta_f \circ p_{d} = f \circ \eta_f$, and
\item $f\mapsto \eta_f$ is continuous on $\BP_d$ with $\eta_{p_d} = \mathrm{id}$.
\end{enumerate}
Therefore, by conjugating with a rotation fixing $0$ that sends $\eta_f(1)$ to $1$, we have a canonical identification of $\BP_{d, \fc}$ and $\BP_d$.

\begin{defn}
Let $\mathcal{S}=(|\mathcal{S}|, \Phi, \delta)$ be a mapping scheme. 
We associate the {\em Blaschke model space} $\BP^\mathcal{S}$ consisting of all proper holomorphic maps
$$
\bp: |\mathcal{S}| \times \D \longrightarrow |\mathcal{S}| \times \D
$$
such that $\mathcal{F}$ carries each $s\times \D$ onto $\Phi(s) \times \D$ by a 1-anchored Blaschke product 
$$
(s, z) \mapsto (\Phi(s), \bp_s(z))
$$
of degree $\delta(s)$ which is either fixed point centered or zeros centered according to $s$ is periodic or aperiodic under $\Phi$.
\end{defn}
We call $\mathcal{F} \in \BP^\mathcal{S}$ a {\em Blaschke mapping scheme}.
The degree of $\bp \in \BP^\mathcal{S}$ satisfies
$\deg (\bp) = \deg \mathcal{S}$, and there are exactly $\deg(\bp)-1$ critical points counted with multiplicities.

\subsection*{Markings of a Blaschke mapping scheme}
Since a Blaschke mapping scheme is $1$-anchored, we have a natural combinatorial marking on the dynamics on the circle defined as below.
It is natural to use $\R/\Z$ to represent these combinatorial angles, and let $m_d: \R/\Z \longrightarrow \R/\Z$ be the multiplication by $d$ map.

\begin{defn}
Let $\bp\in \BP^\mathcal{S}$. We say a family of homeomorphisms 
$$
\eta_s := \eta_{s, \bp}: \R/\Z \longrightarrow \mathbb{S}^1
$$
for $s\in |\mathcal{S}|$ is a {\em marking} for $\bp$ if
\begin{itemize}
\item $\eta_{\Phi(s)} \circ m_{\delta(s)} = \bp_s \circ \eta_s$.
\item $\eta_s(0) = 1$.
\end{itemize}
\end{defn}

For each $\bp\in \BP^\mathcal{S}$, there exists a unique marking.
Indeed, the markings can be constructed for pre-periodic points first. 
Since $\bp$ is expanding on $|\mathcal{S}| \times \mathbb{S}^1$, it extends to homeomorphisms.
The normalization $\eta_s(0) = 1$ guarantees the uniqueness.

\subsection*{Fixedpoint-marked rational maps}
By a {\em fixedpoint-marked rational map} $(f; z_0,..., z_d)$, we mean a rational map $f: \hat \C \longrightarrow \hat \C$ of degree $d\geq 2$, together with an ordered list of its $d+1$ (not necessarily distinct) fixed points $z_j$.

Let $\Rat_{d, \fm}$ be the space of all fixedpoint-marked rational maps of degree $d$.
The groups of M\"obius transformation $\PSL_2(\C)$ acts naturally on $\Rat_{d, \fm}$ by
$$
M\cdot (f; z_0,..., z_d) := (M \circ f \circ M^{-1}; M(z_0),..., M(z_d)),
$$
for $M \in \PSL_2(\C)$.
We define the {\em moduli space} $\mathcal{M}_{d, \fm} = \Rat_{d, \fm}/ \PSL_2(\C)$.
A {\em hyperbolic component} in $\mathcal{M}_{d, \fm}$ is a connected component in the open subset consisting of all conjugacy classes of hyperbolic fixedpoint-marked rational maps.
It is said to have connected Julia set if any representative map $f$ has connected Julia set.

To avoid cumbersome notations, we will simply use $f\in \Rat_{d,\fm}$ or $[f]\in \mathcal{M}_{d,\fm}$ to represent a fixedpoint-marked rational map or conjugacy class.

Let $f$ be a hyperbolic rational map with connected Julia set.
We can associate a scheme $\mathcal{S}_f$ in the following way.
Vertices $s_U \in |\mathcal{S}_f|$ are in correspondence with components $U$ of the Fatou set which contain a critical or post-critical point.
The associated map $F = F_f: |\mathcal{S}_f| \longrightarrow |\mathcal{S}_f|$ carries $s_U$ to $s_{f(U)}$ and the degree $\delta(s_U)$ is defined to be the degree of $f: U \longrightarrow f(U)$.
We have the following theorem
\begin{theorem}[\cite{Milnor12}, Theorem 9.3]\label{thm:BMS}
Let $\mathcal{H} \subseteq \mathcal{M}_{d, \fm}$ be a hyperbolic component with connected Julia set.
Then $\mathcal{H}$ is canonically homeomorphic to the Blaschke model space $\BP^\mathcal{S}$, where $\mathcal{S} = \mathcal{S}_f$ for a (and hence all) representative map $f$.
\end{theorem}

\section{Degenerations in $\BP^\mathcal{S}$}\label{sec:DBS}
In this section, we will study quasi post-critically finite degenerations on $\BP^\mathcal{S}$.
The key results in this section are the construction of the quasi-invariant forests $\mathcal{T}$ as in Theorem \ref{thm:qit}, and the construction of the rescaling limits.

We first define the `metric' $d_{\mathcal{S}}$ on $|\mathcal{S}| \times \D$ by
\begin{itemize}
\item $d_{\mathcal{S}}(x,y) = \infty$ if $x \in s\times \D$ and $y\in t\times\D$ with $s \neq t$;
\item $d_{\mathcal{S}}(x,y) = d_{s\times \D}(x, y)$ if $x, y\in s\times \D$,
\end{itemize}
where $d_{s\times \D}$ is the hyperbolic metric on the unit disk.

\begin{defn}\label{defn:qpcfbm}
Let $\mathcal{S}$ be a scheme of degree $d$ and $\BP^\mathcal{S}$ be the corresponding Blaschke model space.
Let $\bp_n \in \BP^\mathcal{S}$ be a sequence.
It is said to be ($K$-)quasi post-critically finite if we can label the critical points by $c_{1,n},..., c_{2d-2, n}$, and for any sequence $c_{i,n}$ of critical points, there exist $l_i$ and $q_i$ called {\em quasi pre-periods} and {\em quasi periods} respectively such that
$$
d_{\mathcal{S}}(\bp_n^{l_i}(c_{i,n}), \bp_n^{l_i+q_i}(c_{i,n})) \leq K.
$$
We say $\bp_n$ is {\em quasi post-critically finite} if it is $K$-quasi post-critically finite for some $K$.
The sequence $\bp_n$ is said to be {\em degenerating} if it leaves every compact set of $\BP^\mathcal{S}$.
\end{defn}
Note that if $c_{i, n} \in \D_s$, then $\Phi^{l_i}(s)$ is periodic and has period dividing $q_i$.

\subsection{Consctruction of quasi-invariant forests}
In \cite{Luo21}, we constructed quasi-invariant trees for quasi post-critically finite degenerations $f_n \in \BP_d$ (cf. Ribbon $\R$-tree in \cite{McM09}).
The construction of {\em quasi-invariant forests} for $\bp_n \in \BP^S$ is almost identical, with only notational changes.
For completeness, we give the construction here.

We remark that the construction of tree of Riemann spheres in \S \ref{sec:qpcfd} is a $3$-dimensional generalization of the following construction.

\subsection*{Convention on notations}
Since there are multiple subindices for the construction, we make the following convention throughout this section.
\begin{itemize}
\item $s$ is reserved to represent a point in the mapping scheme $|\mathcal{S}|$;
\item $n$ is reserved to index the objects associated to $\mathcal{F}_n$
\end{itemize}
If elements in a sequence have multiple subindices, we use $(-)_n$ to emphasize the index for the sequence is $n$. We shall drop the $n$ if there is no confusion.

We first define a finite forward invariant set of marked points 
$$
\widetilde{\mathcal{P}}_n = \bigcup_{s\in |\mathcal{S}|} \widetilde{\mathcal{P}}_{s,n} \subseteq |\mathcal{S}|\times \D
$$ 
inductively as follows:
\begin{itemize}
\item If $s\in |\mathcal{S}|$ is periodic, then $\widetilde{\mathcal{P}}_{s,n} = \{(s,0)\}$ consists of the attracting periodic point for $\bp_n$ in $s \times \D$;
\item If $s\in |\mathcal{S}|$ is strictly pre-periodic, then $\widetilde{\mathcal{P}}_{s,n} =  \bp_n^{-1}(\widetilde{\mathcal{P}}_{\Phi(s),n}) \cap (s\times \D)$ consists of the full preimage of $\widetilde{\mathcal{P}}_{\Phi(s),n}$ in $s \times \D$.
\end{itemize}
We label the points in $\widetilde{\mathcal{P}}_{s,n} = \{p_{s,1,n},..., p_{s,k_s,n}\}$, and let
$$
\widetilde{\mathcal{P}} = \bigcup_{s\in |\mathcal{S}|} \widetilde{\mathcal{P}}_s = \bigcup_{s\in |\mathcal{S}|} \{(p_{s,j,n})_n: j=1,..., k_s\}
$$
be the set of sequences of marked points.

\begin{figure}[ht]
  \centering
  \resizebox{0.6\linewidth}{!}{
    \def\svgwidth{\columnwidth}
    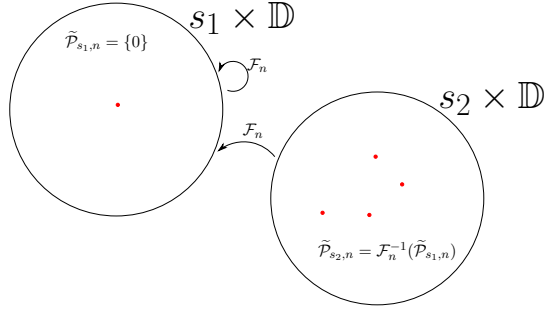

  }
  \caption{An illustration of $\widetilde{\mathcal{P}}$.}
  \label{fig:MarkedP}
\end{figure}

We also define
$$
\widetilde{\mathcal{Q}} = \widetilde{\mathcal{P}} \cup \{(\bp_n^j(c_{i,n}))_n: i = 1,..., 2d-2, j = 0,..., l_i+q_i-1\},
$$
as the union of the set of sequences of marked points and the set of sequences of critical and post-critical points up to $l_i+q_i-1$ iterates.

After passing to a subsequence, we assume that the degrees $\deg(v_n)$ are constant for $(v_n) \in \widetilde{\mathcal{Q}}$ and the limit
$$
\lim_{n\to\infty} d_\mathcal{S}(v_n, w_n)
$$
exists (which can possibly be $\infty$) for any pairs of sequences $(v_n), (w_n) \in \widetilde{\mathcal{Q}}$.
This defines an equivalence relation: 
$$
(v_n) \sim (w_n) \text{ if $\lim_{n\to\infty} d_\mathcal{S}(v_n, w_n) < \infty$.}
$$ 
An equivalence class $\mathcal{C}\in \widetilde{\mathcal{Q}} / \sim$ is called a {\em cluster set}, and its degree is
$$
\deg(\mathcal{C}) = 1+ \sum_{(w_n) \in \mathcal{C}} (\deg(w_n) - 1).
$$

Let $\mathcal{C}\in \widetilde{\mathcal{Q}} / \sim$ be an equivalence class, 
it is convenient to choose a representative $(v_n(\mathcal{C})) \in\mathcal{C}$ with the convention that
if $\mathcal{C} \in \widetilde{\mathcal{P}} / \sim$, then $(v_n(\mathcal{C})) \in \widetilde{\mathcal{P}}$.

\begin{figure}[ht]
  \centering
  \resizebox{0.4\linewidth}{!}{
    \def\svgwidth{\columnwidth}
\begingroup%
  \makeatletter%
  \providecommand\color[2][]{%
    \errmessage{(Inkscape) Color is used for the text in Inkscape, but the package 'color.sty' is not loaded}%
    \renewcommand\color[2][]{}%
  }%
  \providecommand\transparent[1]{%
    \errmessage{(Inkscape) Transparency is used (non-zero) for the text in Inkscape, but the package 'transparent.sty' is not loaded}%
    \renewcommand\transparent[1]{}%
  }%
  \providecommand\rotatebox[2]{#2}%
  \newcommand*\fsize{\dimexpr\f@size pt\relax}%
  \newcommand*\lineheight[1]{\fontsize{\fsize}{#1\fsize}\selectfont}%
  \ifx\svgwidth\undefined%
    \setlength{\unitlength}{792bp}%
    \ifx\svgscale\undefined%
      \relax%
    \else%
      \setlength{\unitlength}{\unitlength * \real{\svgscale}}%
    \fi%
  \else%
    \setlength{\unitlength}{\svgwidth}%
  \fi%
  \global\let\svgwidth\undefined%
  \global\let\svgscale\undefined%
  \makeatother%
  \begin{picture}(1,0.77272727)%
    \lineheight{1}%
    \setlength\tabcolsep{0pt}%
    \put(0,0){\includegraphics[width=\unitlength,page=1]{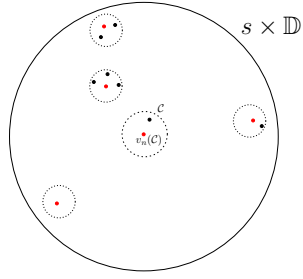}}%
    \put(0.73985778,0.66537185){\color[rgb]{0,0,0}\makebox(0,0)[lt]{\lineheight{1.25}\smash{\begin{tabular}[t]{l}{\Huge $s\times \D$}\end{tabular}}}}%
    \put(0.5135168,0.45155565){\color[rgb]{0,0,0}\makebox(0,0)[lt]{\lineheight{1.25}\smash{\begin{tabular}[t]{l}$\mathcal{C}$\end{tabular}}}}%
    \put(0.4535183,0.36724683){\color[rgb]{0,0,0}\makebox(0,0)[lt]{\lineheight{1.25}\smash{\begin{tabular}[t]{l}$v_n(\mathcal{C})$\end{tabular}}}}%
    \put(0,0){\includegraphics[width=\unitlength,page=2]{Cluster.pdf}}%
  \end{picture}%
\endgroup%

  }
  \caption{An illustration of clusters in $s\times \D$. We choose a representative (red points) for each cluster. As $n\to \infty$, the distance is bounded between points in the same cluster, while goes to infinity between points in different clusters.}
  \label{fig:Cluster}
\end{figure}

We denote 
$$
\mathcal{Q} := \{(v_n(\mathcal{C})): \mathcal{C} \in \widetilde{\mathcal{Q}} / \sim\}
$$ 
as the set of representative sequences. Let
$$
\mathcal{P} := \{(v_n) \in \mathcal{Q}: [(v_n)] \in \widetilde{\mathcal{P}} / \sim\} \subseteq \mathcal{Q}.
$$
We set $\mathcal{Q}_n  = \{v_n: (v_n) \in \mathcal{Q} \}$.
We denote $\mathcal{Q}_s:= \mathcal{Q} \cap (s \times \D)$ and $\mathcal{Q}_{s, n} := \mathcal{Q}_n \cap (s \times \D)$.
We denote $\mathcal{P}_n, \mathcal{P}_{s}$ and $\mathcal{P}_{s, n}$ accordingly.

Now fix $s\in |\mathcal{S}|$.
Label the elements of $\mathcal{Q}_s = \{(b_{0,n})_n, (b_{1,n})_n,..., (b_{m,n})_n\}$.
The sequence of quasi-invariant trees $\mathcal{T}_{s,n}$ (whose name will be justified in Theorem \ref{thm:qit}) is the `spine' of the degenerating hyperbolic polygon $\chull(\mathcal{Q}_{s, n})$ and is constructed inductively as follows.

\begin{figure}[ht]
  \centering
  \resizebox{0.9\linewidth}{!}{
    \def\svgwidth{\columnwidth}
    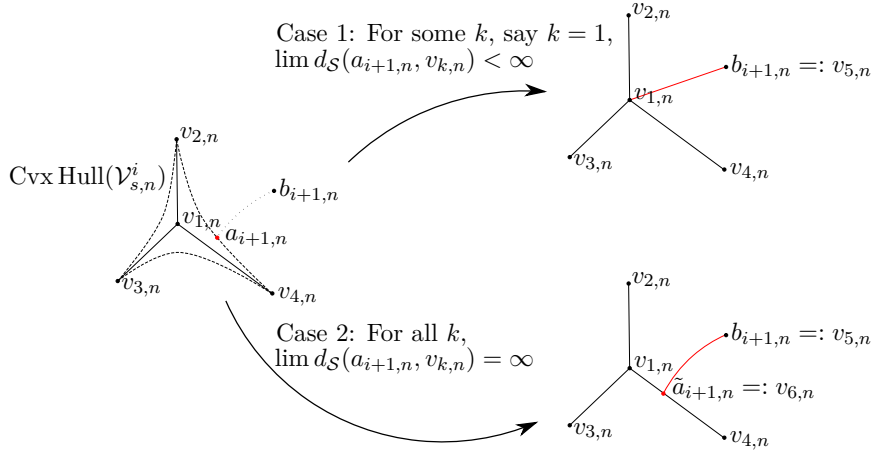

  }
  \caption{An illustration of the inductive construction of $\mathcal{T}_{s,n}$.}
  \label{fig:Construction}
\end{figure}

As the base case, we define $\mathcal{T}^0_{s,n} = \{b_{0,n}\}$ with vertex set $\mathcal{V}^0_{s,n}=\{b_{0,n}\}$.

Assume that $\mathcal{T}^i_{s,n}$ is constructed with vertex set $\mathcal{V}^i_{s,n}:= \{v_{1,n},..., v_{m_i,n}\}$ containing $\{b_{0,n},..., b_{i,n}\}$. 
Assume as induction hypotheses that 
\begin{align}\label{Eqn:H1}
\lim_{n\to\infty}d_{\mathcal{S}}(v_{k,n}, v_{k',n}) = \infty
\end{align}
for all $k\neq k'$ and 
\begin{align}\label{Eqn:H2}
\min_{k=i+1,..., m} d_{\mathcal{S}}(b_{k,n}, \chull(\mathcal{V}^{i}_{s, n})) \to \infty.
\end{align}
We also assume as an induction hypothesis that the edges of $\mathcal{T}^i_{s,n}$ are hyperbolic geodesics.

After passing to a subsequence and changing indices, we assume for all $n$,
$$
d_{\mathcal{S}}(b_{i+1,n}, \chull(\mathcal{V}^{i}_{s, n})) = \min_{k=i+1,..., m} d_{\mathcal{S}}(b_{k,n}, \chull(\mathcal{V}^{i}_{s, n})).
$$
Let $a_{i+1, n}$ be the projection of $b_{i+1,n}$ onto the convex hull $\chull(\mathcal{V}^i_{s, n})$.
After passing to a subsequence, we assume
$\lim_{n\to\infty} d_{\mathcal{S}}(a_{i+1, n}, v_{k,n})$
exist for all $k$ (which can possibly be $\infty$).
We also assume $a_{i+1, n}$ is on the hyperbolic geodesic $[v_{j_1,n}, v_{j_k,n}]$ and $v_{j_1,n}, v_{j_k,n}$ are connected by a sequence of geodesics $[v_{j_1,n}, v_{j_2,n}] \cup [v_{j_2,n}, v_{j_3,n}] \cup ... \cup [v_{j_{k-1}, n}, v_{j_k,n}]$ in $\mathcal{T}^i_{s, n}$ (see Figure \ref{fig:Construction}).

We have two cases, if $a_{i+1, n}$ stay a bounded distance from some vertex $v_{k,n} \in \mathcal{V}^i_{s, n}$, then we define
$$
\mathcal{T}^{i+1}_{s, n} := \mathcal{T}^{i}_{s, n} \cup [v_{k,n}, b_{i+1,n}],
$$
and set the vertex set $\mathcal{V}^{i+1}_{s, n}:=\mathcal{V}^i_{s, n} \cup \{b_{i+1,n}\}$ (see Case 1 in Figure \ref{fig:Construction}).

Otherwise, we can subdivide the polygon $v_{j_1,n},..., v_{j_k,n}$ into finitely many triangles. 
Since the hyperbolic triangles are thin, there exists $\tilde{a}_{i+1,n} \in [v_{j_1,n}, v_{j_2,n}] \cup ... \cup [v_{j_{k-1}, n}, v_{j_k,n}] \subseteq \mathcal{T}^i_{s, n}$ staying a bounded distance from $a_{i+1, n}$. Note that the bound depends only on the number of vertices.
We define 
$$
\mathcal{T}^{i+1}_{s, n} := \mathcal{T}^{i}_{s, n} \cup [\tilde{a}_{i+1,n}, b_{i+1,n}],
$$
and set the vertex set $\mathcal{V}^{i+1}_{s, n}:=\mathcal{V}^i_{s, n} \cup \{\tilde{a}_{i+1,n}, b_{i+1,n}\}$ (see Case 2 in Figure \ref{fig:Construction}).

By our construction, it is easy to verify the equations \ref{Eqn:H1} and \ref{Eqn:H2} are satisfied for $\mathcal{T}^{i+1}_{s, n}$ and the edges are hyperbolic geodesics.
Thus, by induction,  let $\mathcal{T}_{s, n} = \mathcal{T}^m_{s, n}$ be the finite tree when all $m$ points in $\mathcal{Q}_{s, n}$ are added.

By construction, $\mathcal{Q}_{s, n} \subseteq \mathcal{V}_{s, n}$, and any point in $\mathcal{V}_{s, n} - \mathcal{Q}_{s, n}$ is a branch point for $\mathcal{T}_{s, n}$. 
Moreover, $\lim_{n\to\infty} d_{\mathcal{S}}(v_{k,n}, v_{k',n}) = \infty$ for all $k\neq k'$.

We mark $\mathcal{T}_{s, n}$ by the set $\mathcal{P}_{s,n} \subseteq \mathcal{V}_{s,n}$, and $(\mathcal{T}_{s, n}, \mathcal{P}_{s,n})$ is regarded as a {\em marked tree} \footnote{Note that if $s$ is periodic, then there is only one marked point. We thus get a pointed tree as in \cite{Luo21}.}.

A {\em ribbon structure} on a finite tree is the choice of a planar embedding up to isotopy. The ribbon structure can be specified by a cyclic ordering of the edges incident to each vertex.
After passing to a further subsequence, we may assume that $\mathcal{T}_{s, n}$ are all canonically isomorphic with the same ribbon structure.
Thus, we denote $\mathcal{T}_s$ as the underlying ribbon finite tree, with the vertex set $\mathcal{V}_s$ and the isomorphisms 
$$
\phi_{s, n}: (\mathcal{T}_s, \mathcal{P}_s) \longrightarrow (\mathcal{T}_{s, n}, \mathcal{P}_{s,n}).
$$
The local degree $\delta(v)$ at a vertex $v\in \mathcal{Q}_{s} \subset \mathcal{V}_s$ is defined as the degree of the corresponding cluster, and $\delta(v) = 1$ if $v \in \mathcal{V}_s - \mathcal{Q}_s$.
We remark that this local degree can also be defined using the rescaling limits as in \S \ref{subsec:rl}.

We define a {\em forest} as a finite disjoint union of finite trees.
We define
$$
\phi_{n}: (\mathcal{T}, \mathcal{P}) := \bigcup_{s\in |\mathcal{S}|} (\mathcal{T}_s, \mathcal{P}_s) \longrightarrow (\mathcal{T}_n, \mathcal{P}_n) := \bigcup_{s\in |\mathcal{S}|} (\mathcal{T}_{s, n}, \mathcal{P}_{s,n})
$$
as the union of the maps $\phi_{s,n}$ and the vertex set
$$
\mathcal{V}:= \bigcup_{s\in |\mathcal{S}|} \mathcal{V}_s.
$$

A map $F$ between finite trees is said to be {\em simplicial} if $F$ sends an edge to an edge.
A map $F$ between forests is said to be {\em simplicial} if $F$ is simplicial restricting to each tree component.

Since $\bp_n$ is quasi post-critically finite, the vertices of $\mathcal{T}_n$ associated to the critical and post-critical points are uniformly quasi-invariant.
Since a proper holomorphic map converges to an isometry exponential fast away from the critical points (see Theorem 10.11 in \cite{McM09}), we have the following properties for the quasi-invariant forest.
The proof is similar as in Theorem 2.2 in \cite{Luo21}.

\begin{theorem}\label{thm:qit}
Let $\mathcal{S} = (|\mathcal{S}|, \Phi, \delta)$ be a mapping scheme and $\bp_n \in \BP^\mathcal{S}$ be quasi post-critically finite. 
After passing to a subsequence, there exists a constant $K > 0$ and a simplicial map 
$$
F = \cup_{s\in |\mathcal{S}|} F_s:  (\mathcal{T}, \mathcal{P})=\bigcup_{s\in |\mathcal{S}|} (\mathcal{T}_s, \mathcal{P}_s) \longrightarrow (\mathcal{T}, \mathcal{P})\footnote{We remark that for $F$ to be simplicial, we may need to add finitely many vertices on those strictly pre-periodic edges to our construction $\mathcal{T}$.}
$$ with
$$
F_s: (\mathcal{T}_s, \mathcal{P}_s) \longrightarrow (\mathcal{T}_{\Phi(s)}, \mathcal{P}_{\Phi(s)})
$$
on a marked finite ribbon forest with vertex set $\mathcal{V}=\bigcup_{s\in |\mathcal{S}|}\mathcal{V}_s$, and a sequence of isomorphisms
$$
\phi_{n}: (\mathcal{T}, \mathcal{P}) \longrightarrow (\mathcal{T}_n, \mathcal{P}_n)
$$
such that
\begin{itemize}
\item (Degenerating vertices.) If $v_1\neq v_2 \in \mathcal{V}$, then 
$$
d_\mathcal{S}(\phi_n(v_1), \phi_n(v_2)) \to \infty.
$$
\item (Geodesic edges.) If $E =[v_1,v_2] \subseteq \mathcal{T}$ is an edge, then the corresponding edge $\phi_n(E) \subseteq \mathcal{T}_n$ is a hyperbolic geodesic segment connecting $\phi_n(v_1)$ and $\phi_n(v_2)$.
\item (Critically approximating.) Any critical points of $\bp_n$ are within $K$ distance from the vertex set $\mathcal{V}_n:= \phi_n(\mathcal{V})$ of $\mathcal{T}_n$.
\item (Quasi-invariance on vertices.)
If $v\in \mathcal{V}$, then
$$
d_\mathcal{S}(\bp_n(\phi_n(v)), \phi_n(F(v))) \leq K \text{ for all } n.
$$

\item (Quasi-invariance on edges.) If $E\subseteq \mathcal{T}$ is an edge and $x_n \in \phi_n(E)$, then there exists $y_n \in \phi_n(F(E))$ so that
$$
d_\mathcal{S}(\bp_n(x_n), y_n) \leq K \text{ for all } n.
$$
If $E$ is a periodic edge of period $q$, then 
$$
d_\mathcal{S}(\bp_n^q(x_n), x_n) \leq K \text{ for all } n.
$$
\end{itemize}
\end{theorem}

We remark that Theorem \ref{thm:qit} states that the dynamics of $\bp_n$ on $\mathcal{T}_n$ is $K$-quasi-invariant, and is modeled by a simplicial map $F: (\mathcal{T}, \mathcal{P}) \longrightarrow (\mathcal{T}, \mathcal{P})$.
We also remark that this constant $K$ can be much bigger than the quasi post-critically finite constant for $\bp_n$.

For future reference, we define
\begin{defn}
Let $v\in \mathcal{V}$. It is called a {\em Fatou point} if $v$ is eventually mapped to a periodic cycle of local degree $\geq 2$ under $F$, and is called a {\em Julia point} otherwise.
\end{defn}

\subsection{Rescaling limits on vertices}\label{subsec:rl}
\subsection*{Degeneration of rational maps}
Let $\Rat_d$ be the space of rational maps of degree $d$.
By fixing a coordinate system, a rational map can be expressed as a ratio of two homogeneous polynomials $f(z:w) = (P(z,w):Q(z,w))$, where $P$ and $Q$ have degree $d$ with no common divisors.
Thus, using the coefficients of $P$ and $Q$ as parameters, we have
$$
\Rat_d = \Proj^{2d+1} \setminus V(\Res),
$$ 
where $\Res$ is the resultant of the two polynomials $P$ and $Q$, and $V(\mathrm{Res})$ is the hypersurface for which the resultant vanishes.
This embedding gives a natural compactification $\overline{\Rat_d} = \Proj^{2d+1}$, which will be called the {\em algebraic compactification}.
Maps $f\in \overline{\Rat_d}$ can be written as
$$
f = (P:Q) = (Hp: Hq),
$$
where $H = \gcd(P, Q)$.
We set $$\varphi_f := (p:q),$$ which is a rational map of degree at most $d$.
The zeroes of $H$ in $\hat\C$ are called the {\em holes} of $f$, and the set of holes of $f$ is denoted by $\mathcal{H}(f)$.

If $\{f_n\}\subseteq \Rat_d$ converges to $f \in \overline{\Rat_d}$, we say that $f$ is the {\em algebraic limit} of the sequence $\{f_n\}$. 
It is said to have degree $k$ if $\varphi_f$ has degree $k$.
Abusing notations, sometimes we shall refer to $\varphi_f$ as the algebraic limit of $\{f_n\}$, and write $f_n \to \varphi_f$. 
The following is useful in analyzing the limiting dynamics.

\begin{lem}[\cite{DeM05} Lemma 4.2]\label{lem:ac}
If $f_n$ converges to $f$ algebraically, then $f_n$ converges to $\varphi_f$ uniformly on compact subsets of $\widehat{\C}-\mathcal{H}(f)$.
\end{lem}

The following statement is also useful in many situations:
\begin{lem}[\cite{Luo21} Lemma 2.2]\label{lem:critc}
Let $f_n$ converges to $f$ algebraically with $\deg(f) \geq 1$.
If $x\in \mathcal{H}(f)$, then there exists a sequence of critical points $c_n$ for $f_n$ with $\lim_{n\to\infty}c_n = x$.
\end{lem}

We shall also use the following:
\begin{lem}[\cite{Luo21} Proposition 2.3]\label{lem:ac1}
Let $K \geq 0$, and let $f_n: \D\longrightarrow \D$ be a sequence of proper holomorphic map of degree $d$ with $d_{\Hyp^2}(0, f_n(0)) \leq K$.
Then after passing to a subsequence, $f_n$ converges compactly to a proper holomorphic map $f: \D \longrightarrow \D$ of potentially lower degree.

By Schwarz reflection, $f, f_n$ are rational maps, and $f_n$ converges algebraically to $f$ with holes contained in the unit circle $\mathbb{S}^1$.
\end{lem}

\subsection*{Rescaling limits}
Let $v \in \mathcal{V}_s \subseteq \mathcal{V}$. 
We define a {\em normalization at $v$} or a {\em coordinate at $v$} as a sequence $M_{v,n}\in \Isom(\Hyp^2)$ so that 
$$
\phi_n(v) = (s, M_{v,n}(0)).
$$
Note that different choices for the sequence $M_{v,n}$ are differed by pre-composing with rotations that fix $0$, which form a compact group.
We define the map 
\begin{align*}
\hat M_{v,n} : \D &\longrightarrow s\times \D\\
z &\mapsto (s, M_{v,n}(z)).
\end{align*}

Let us fix such a normalization $M_{v,n}$ for each vertex $v\in \mathcal{V}$ .
We can associate a limiting disk $\D_v$ and a limiting circle $\mathbb{S}^1_v$ to the vertex $v$.
More precisely, we say a sequence $z_n \in \overline{\D}$ converges to $z\in \overline{\D}_v$ in $v$-coordinate, denoted by $z_n \to_v z$ or $z = \lim_v z_n$ if
$$
\lim_{n\to\infty} M^{-1}_{v,n} (z_n) = z.
$$
More generally, a sequence $(s_n, z_n) \in |\mathcal{S}| \times \overline{\D}$ converges to $z\in \overline{\D}_v$ in $v$-coordinate, also denoted by $(s_n, z_n) \to_v z$ if $s_n = s$ for all sufficiently large $n$ and $z_n \to_v z$.

\begin{figure}[ht]
  \centering
  \resizebox{0.9\linewidth}{!}{
    \def\svgwidth{\columnwidth}
    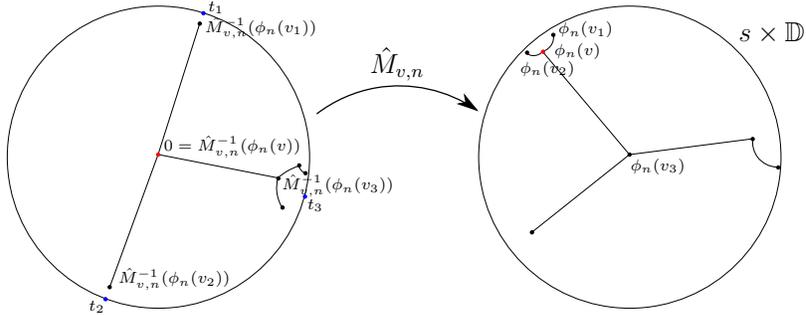

  }
  \caption{An illustration of the map $\hat M_{v,n}$, and the marked points $t_i \in \mathbb{S}^1_v$.}
  \label{fig:RL}
\end{figure}

By Theorem \ref{thm:qit}, there exists a constant $K$ such that
$$
d_{\mathcal{S}}(0,\hat M_{F(v),n}^{-1}\circ \bp_n \circ \hat M_{v,n}(0)) = d_{\mathcal{S}}(\hat M_{F(v),n}(0),\bp_n\circ \hat M_{v,n}(0)) \leq K.
$$
Thus, by Lemma \ref{lem:ac1}, after possibly passing to subsequences, the sequence of proper maps
$$
\hat M_{F(v),n}^{-1}\circ \bp_n \circ \hat M_{v,n}: \D \longrightarrow \D
$$
converges compactly to a proper holomorphic map
$$
\rl_v = \rl_{v\to F(v)}: \D_v \longrightarrow \D_{F(v)}.
$$

More generally,
$\hat M_{F^k(v),n}^{-1}\circ \bp^k_n \circ \hat M_{v,n}$
converges to 
$$
\rl^k_{v}:=\rl_{F^{k-1}(v)}\circ \rl_{F^{k-2}(v)}\circ...\circ \rl_{v}.
$$
In particular, if $v$ is a periodic point of period $q$, then
$$
\hat M_{v,n}^{-1}\circ \bp^q_n \circ \hat M_{v,n} \to \rl^q_v.
$$

\begin{defn}\label{defn:rlm}
We call $\rl_v$ the {\em rescaling limit} between $v$ and $F(v)$.
\end{defn}

Let $v_1,..., v_l$ be the list of adjacent vertices of $v$ in $\mathcal{T}$. 
Then $\{v_1,..., v_l\}$ can be identified with the tangent space $T_v\mathcal{T}$ of $\mathcal{T}$ at $v$.
By our construction, after passing to a subsequence, there exist $l$ distinct points $t_1,..., t_l \in \mathbb{S}^1_v$ with
$t_i = \lim_v \phi_n(v_i)$.
We denote this correspondence by the map
$$
\xi_v: T_v\mathcal{T} \longrightarrow \mathbb{S}^1_v.
$$

Since $F$ is simplicial, we have a well-defined tangent map $DF_v:T_v\mathcal{T} \longrightarrow T_{F(v)}\mathcal{T}$. 
On the other hand, the rescaling limit $\rl_{v}$ gives a maps from $\mathbb{S}^1_v$ to $\mathbb{S}^1_{F(v)}$.
We can verify that this marking is compatible with the dynamics of the tangent map (see the proof of Lemma 2.13 in \cite{Luo21}).
\begin{lem}\label{lem:ct}
Let $v\in \mathcal{V}$. Then
$$
\rl_{v} \circ \xi_v = \xi_{F(v)} \circ DF_v.
$$
\end{lem}

\subsection{Pullback and dual laminations}
\subsection*{Pullback of quasi-invariant forest}
Given the quasi-invariant forest $\mathcal{T}_n$ for $\bp_n$ modeled by $F:(\mathcal{T}, \mathcal{P}) \longrightarrow (\mathcal{T}, \mathcal{P})$, we can use the same method to construct the pullback quasi-invariant tree $\mathcal{T}^1_n$.
Indeed, start with $\mathcal{T}_n$, we add new vertices for new clusters of $\bp_n^{-1}(\phi_n(\mathcal{V}))$ and use the same method to construct $\mathcal{T}^1_n$.
Here, a cluster is new if the distance between any point in the cluster to the existing vertices $\phi_n(\mathcal{V})$ is unbounded.
A simplicial model of the dynamics can be constructed for the pullback and is denoted by 
$$
F: (\mathcal{T}^1, \mathcal{P})\longrightarrow (\mathcal{T}^0, \mathcal{P}) := (\mathcal{T}, \mathcal{P})\subseteq (\mathcal{T}^1, \mathcal{P}).
$$ 

The rescaling limits are defined similarly and the same proof shows the compatibility of the local dynamics with the tangent map.
We remark that all new vertices have local degree $1$, so the rescaling limits are defined by degree $1$ maps.
The pullback can be iterated, and we denote the $k$-th pullback simplicial model by 
$F: (\mathcal{T}^k, \mathcal{P}) \longrightarrow (\mathcal{T}^{k-1}, \mathcal{P}) \subseteq (\mathcal{T}^k, \mathcal{P})$.
We also denote $F: (\mathcal{T}^\infty, \mathcal{P}) \longrightarrow (\mathcal{T}^\infty, \mathcal{P})$ as the union of these maps.

\subsection*{Landing angles for edges}
Recall that any $\bp\in \BP^\mathcal{S}$ is equipped with a unique marking
$$
\eta_{s, \bp}: \R/\Z \longrightarrow \mathbb{S}^1.
$$
Denote $\eta_{s, n}: \R/\Z \longrightarrow \mathbb{S}^1$, $s\in |\mathcal{S}|$ as the marking for $\bp_n$.

Let $v\in \mathcal{T}_s$ be a vertex.
After passing to a subsequence, we assume $\lim_v \eta_{s,n}(t)$ exist for all rational angles.

Let $E$ be an edge of $\mathcal{T}_s$ incident at $v$ with corresponding tangent vector $x_E \in T_v\mathcal{T}_s$.
We consider the set
$$
I_E:=\{t \in \R/\Z: \lim_v \eta_{s,n}(t) = \xi_v(x_E) \in \mathbb{S}^1_v\}.
$$
Since $\eta_{s,n}$ is a homeomorphism, by the cyclic ordering, $I_E$ is an interval and we denote the two endpoints by $t^\pm_E$.

We say $t^\pm_E \in \R/\Z$ {\em land} at $E$, or $t^\pm_E$ are landing angles for $E$.
It can be checked that the landing angles are well defined: the construction does not depend on the choice of the end point of the edge $E$ (see Corollary 2.18 in \cite{Luo21}).
The ribbon structure of $\mathcal{T}_s$ gives each edge two {\em sides}, and each landing angle $t^\pm_E$ corresponds to one side of the edge $E$.
It is also easy to check the landing angles are compatible with the dynamics:
\begin{lem}\label{lem:lacs}
Let $E$ be an edge of $\mathcal{T}$.
There exist two angles $t^\pm_E \in \R/\Z$ landing at $E$ corresponding to two sides of $E$.
If $t_E^\pm$ are two angles landing at $E$, then $m_{\delta(s)}(t_E^\pm)$ are the two angles landing at the edge $F(E)$.
\end{lem}


The landing angles at an edge $E$ can also be constructed in the following way, providing a slightly different perspective.
We assume $E = [v, w]$ is periodic, the aperiodic case can be constructed by pull back.
After passing to an iterate, we also assume $E$ is fixed.
Let $x_n\in \phi_{s,n}(E)$ be such that
$d_\mathcal{S}(x_n, \phi_{s,n}(v)) \to \infty$ and $d_\mathcal{S}(x_n, \phi_{s,n}(w)) \to \infty$.
Let $M_{x, n} \in \Isom (\Hyp^2)$ with 
\begin{align*}
\hat M_{x, n}(0) &:= (s, M_{x, n}(0)) = x_n.
\end{align*}
After passing to a subsequence, the map $\hat M_{x,n}^{-1}\circ \bp_n \circ \hat M_{x,n}$ converges to a degree $1$ map $\mathcal{R}_x$ as $x_n$ is quasi-fixed and its distance to any critical point is going to infinity.
Under this normalization, $v\to_x t_1 \in \mathbb{S}^1_x$ and $w \to_x t_2 \in \mathbb{S}^1_x$ with $t_1 \neq t_2$ and $\mathcal{R}_x(t_i) = t_i$.
Define
$$
I_i:=\{ x \in \R/\Z: \eta_{s,n}(x) \to_v t_i\}, i=1, 2.
$$
Since the holes of $\mathcal{R}_x$ are contained in $\{t_1, t_2\}$, $I_1$ and $I_2$ are two intervals of $\R/\Z$ sharing two common boundary points, which are the landing angles $t^\pm_E$.

\subsection*{Dual laminations}
A {\em lamination} $\mathcal{L}$ is a family of disjoint hyperbolic geodesics in $\D$ together with the two endpoints in $\mathbb{S}^1\cong \R/\Z$, whose union $|\mathcal{L}| :=\bigcup\mathcal{L}$ is closed.
An element of the lamination is called a {\em leaf} of the lamination.

We define the leaf associated with the edge $E$ as the hyperbolic geodesics in $\D$ connecting $t^\pm_E$.
It is easy to check that leaves for different edges have disjoint interiors.

The {\em dual finite lamination} for $(\mathcal{T}_s,\mathcal{P}_s)$ is defined as the union of all leaves for edges of $\mathcal{T}_s$, and is denoted by $\mathcal{L}^F_s$.
The leaves can be constructed for edges of any pullbacks of $\mathcal{T}_s$.
We call the closure of the union of leaves for edges in $\mathcal{T}_s^\infty$ the {\em dual lamination} of $\mathcal{T}_s$, and is denoted by $\mathcal{L}_s$.

The lamination $\mathcal{L}_s$ gives an {\em equivalence relation} $\sim_s$ on $\mathbb{S}^1$: $a\sim_s b$ if there exists a finite chain of leaves connecting $a$ and $b$.
Note that different laminations may generate the same equivalence relation.

We denote the union by $\mathcal{L}^F:= \bigcup_{s\in |\mathcal{S}|} \mathcal{L}^F_s$ and $\mathcal{L}:= \bigcup_{s\in |\mathcal{S}|} \mathcal{L}_s$.
By Lemma \ref{lem:lacs}, the dual lamination $\mathcal{L}$ is generated by leaves in $\mathcal{L}^F$.

\section{Limits of quasi post-critically finite degenerations}\label{sec:qpcfd}
In this section, we shall study quasi post-critically finite degenerations of rational maps $[f_n]$.
The main goal of this section is to prove Theorem \ref{thm:crm}.

We first give the definition of rational maps on a tree of Riemann spheres in \S \ref{subsec:trs}.
We introduce dynamically meaningful rescalings for $[f_n]$ in \S \ref{subsec:rescalings}.
Using a 3-dimensional generalization of the construction of quasi-invariant trees in \S \ref{sec:DBS}, we construct the rescaling trees $\mathscr{T}_n \subseteq \Hyp^3$ in \S \ref{subsec:rt}.
We study the rescaling limits in \S \ref{subsec:rrm} and the induced map on the rescaling trees in \S \ref{subsec:rtm}.
Finally, the proof of Theorem \ref{thm:crm} is given in \S \ref{subsec:gfl}.

\subsection{Rational maps on trees of Riemann spheres}\label{subsec:trs}
Let $\RT$ be a finite tree with vertex set $\RV$.
A continuous map $F: (\RT, \RV) \longrightarrow (\RT, \RV)$ is called a {\em tree map} if $F$ is injective on each edge.
A tree map gives a tangent map, which we denote by $DF_a: T_a\RT \longrightarrow T_{F(a)}\RT$.

\begin{defn}\label{defn:trs}
A {\em tree of Riemann spheres} $(\RT, \hat\C^\RV)$ consists of a finite tree $\RT$ with vertex set $\RV$, a disjoint union of Riemann spheres $\hat \C^\RV := \bigcup_{a\in \RV}\hat \C_{a}$, together with markings $\xi_a: T_a\RT \xhookrightarrow{} \hat\C_a$ for $a\in \RV$.
The image $\Xi_a :=\xi_a(T_a\RT)$ is called the {\em singular set} at $a$, and $\Xi = \bigcup_{a\in \RV}\Xi_a$.

A {\em rational map} $(F, R)$ on a tree of Riemann spheres is a tree map $F: (\RT, \RV) \longrightarrow (\RT, \RV)$ and a union of maps $R:= \bigcup_{a\in \RV} R_a$ so that
\begin{itemize}
\item $R_a: \hat\C_a \longrightarrow \hat \C_{F(a)}$ is a rational map of degree $\geq 1$;
\item $R_a \circ \xi_a = \xi_{F(a)} \circ DF_a$.
\end{itemize}
It is said to have degree $d$ if $R$ has $2d-2$ critical points in $\hat\C^\RV-\Xi$.
The local degree of $a\in \RV$ is defined as the degree of $R_a$.

A sequence $f_n$ of degree $d$ rational maps is said to {\em converge} to $(F, R)$ on $(\RT, \hat\C^\RV)$ if there exist rescalings $A_{a,n} \in \PSL_2(\C)$ for $a \in \RV$ such that
\begin{itemize}
	\item $A_{F(a),n}^{-1} \circ f_n \circ A_{a,n}(z) \to R_a(z)$ compactly on $\hat \C_a- \Xi_a$;
	\item $A_{b,n}^{-1} \circ A_{a,n}(z)$ converges algebraically to the constant map $\xi_a(v)$, where $v\in T_a\RT$ is the tangent vector in the direction of $b$.
\end{itemize}
\end{defn}

The Fatou set $\Omega(R) \subset \hat\C^\RV$ is the largest open set for which the iterations $\{R^n|_{\Omega}: n \geq 1\}$ form a normal family. 
The Julia set $\mathcal{J}(R)$ is the complement of the Fatou set.
We say $(F, R)$ is {\em geometrically finite} if every critical point of $R$ in the Julia set has a finite orbit.



\subsection{Rescalings}\label{subsec:rescalings}
To set up the notations, let $\mathcal{H} \subseteq \mathcal{M}_{d, \fm}$ be a hyperbolic component with connected Julia set. 
Let $\BP^\mathcal{S}$ be the corresponding Blaschke model space.
By Theorem \ref{thm:BMS}, $[f] \in \mathcal{H}$ corresponds to a unique $\bp \in \BP^\mathcal{S}$.

For $s\in |\mathcal{S}|$, we denote the corresponding Fatou components of $f$ by $\U_s$ and $\U = \bigcup_{s\in |\mathcal{S}|} \U_s$.
The metric $d_\mathcal{S}$ on $|\mathcal{S}| \times \D$ induces a metric $d_\U$ on $\U$:
\begin{itemize}
\item $d_{\U}(x,y) = \infty$ if $x \in \U_s$ and $y\in \U_t$ with $s \neq t$;
\item $d_{\U}(x,y) = d_{\U_s}(x, y)$ if $x, y\in \U_s$,
\end{itemize}
where $d_{\U_s}$ is the hyperbolic metric on $\U_s$.

We say a sequence in $\mathcal{H}$ is {\em quasi post-critically finite} if the corresponding sequence in $\BP^\mathcal{S}$ is quasi post-critically finite; and is a {\em degeneration} if the sequence escapes every compact set of $\mathcal{H}$.

Throughout this section, we fix a quasi post-critically finite degeneration $([f_n])_n \in \mathcal{H}$ with the corresponding sequence $(\bp_n)_n \in \BP^\mathcal{S}$.
We fix $f_n \in \Rat_d$ as representatives for $[f_n]$.
The corresponding Fatou component for a vertex $s\in |\mathcal{S}|$ is denoted by $\U_{s,n}$ and $\U_n = \bigcup_{s\in |\mathcal{S}|} \U_{s,n}$.

Since $f_n$ on $\U_n$ is conjugate to the Blaschke mapping scheme $\bp_n$ on $|\mathcal{S}| \times \D$, we can define the quasi-invariant forests for $f_n$ using the conjugacy.
Abusing the notations, we shall denote the quasi-invariant forests for $f_n$ by $\T_n$ as well, and the marking by
$$
\psi_n: (\mathcal{T}, \mathcal{P}) \longrightarrow (\T_n, \mathcal{P}_{n}).
$$

Since a rational map does not have a natural normal form, the following `rescalings' at vertices are introduced to give a dynamical meaningful normalizations.
\begin{defn}
Let $v\in \mathcal{V}_s \subseteq \mathcal{V}$ and $v_n := \psi_n(v) \in \T_{s,n} \subseteq \U_{s,n}$.
A sequence $A_{v,n} \in \PSL_2 (\C)$ is defined to be a {\em rescaling} for $v$ if
\begin{itemize}
\item $A_{v,n}(0) = v_n$;
\item $A_{v,n}(B(0,1))\subseteq \U_{s,n}$;
\item $A_{v,n}(1), A_{v,n}(\infty)\in \hat\C - \U_{s,n}$.
\end{itemize}
\end{defn}

\begin{figure}[ht]
  \centering
  \resizebox{0.4\linewidth}{!}{
    \def\svgwidth{\columnwidth}
\begingroup%
  \makeatletter%
  \providecommand\color[2][]{%
    \errmessage{(Inkscape) Color is used for the text in Inkscape, but the package 'color.sty' is not loaded}%
    \renewcommand\color[2][]{}%
  }%
  \providecommand\transparent[1]{%
    \errmessage{(Inkscape) Transparency is used (non-zero) for the text in Inkscape, but the package 'transparent.sty' is not loaded}%
    \renewcommand\transparent[1]{}%
  }%
  \providecommand\rotatebox[2]{#2}%
  \newcommand*\fsize{\dimexpr\f@size pt\relax}%
  \newcommand*\lineheight[1]{\fontsize{\fsize}{#1\fsize}\selectfont}%
  \ifx\svgwidth\undefined%
    \setlength{\unitlength}{841.88976378bp}%
    \ifx\svgscale\undefined%
      \relax%
    \else%
      \setlength{\unitlength}{\unitlength * \real{\svgscale}}%
    \fi%
  \else%
    \setlength{\unitlength}{\svgwidth}%
  \fi%
  \global\let\svgwidth\undefined%
  \global\let\svgscale\undefined%
  \makeatother%
  \begin{picture}(1,0.70707071)%
    \lineheight{1}%
    \setlength\tabcolsep{0pt}%
    \put(0,0){\includegraphics[width=\unitlength,page=1]{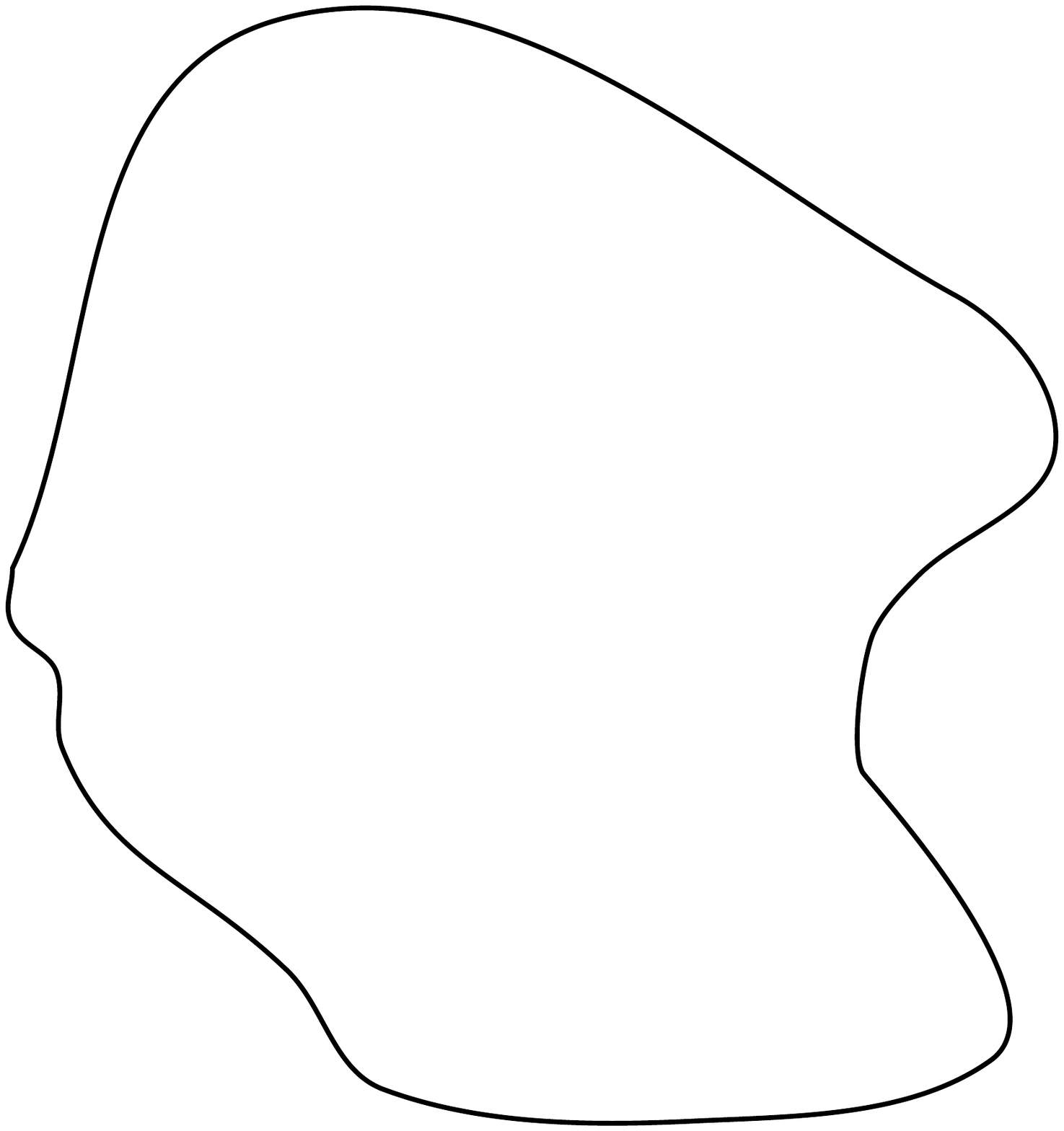}}%
    \put(0.42435396,0.42669801){\color[rgb]{0,0,0}\makebox(0,0)[lt]{\lineheight{1.25}\smash{\begin{tabular}[t]{l}{\Large$\mathcal{U}_{s,n}$}\end{tabular}}}}%
    \put(0,0){\includegraphics[width=\unitlength,page=2]{Norm.pdf}}%
    \put(0.27074791,0.59467904){\color[rgb]{0,0,0}\makebox(0,0)[lt]{\lineheight{1.25}\smash{\begin{tabular}[t]{l}$v_n=A_{v,n}(0)$\end{tabular}}}}%
    \put(0.27009743,0.66522373){\color[rgb]{0,0,0}\makebox(0,0)[lt]{\lineheight{1.25}\smash{\begin{tabular}[t]{l}$A_{v,n}(1)$\end{tabular}}}}%
    \put(0,0){\includegraphics[width=\unitlength,page=3]{Norm.pdf}}%
    \put(0.28269598,0.53113903){\color[rgb]{0,0,0}\makebox(0,0)[lt]{\lineheight{1.25}\smash{\begin{tabular}[t]{l}$A_{v,n}(B(0,1))$\end{tabular}}}}%
    \put(0.84470666,0.6085302){\color[rgb]{0,0,0}\makebox(0,0)[lt]{\lineheight{1.25}\smash{\begin{tabular}[t]{l}$A_{v,n}(\infty)$\end{tabular}}}}%
  \end{picture}%
\endgroup%

  }
  \caption{The rescaling at $v$.}
  \label{fig:Norm}
\end{figure}

The following lemma says differences between different choices are bounded.
\begin{lem}
Let $v \in \mathcal{V}$.
If $A_{v,n}, B_{v,n} \in \PSL_2(\C)$ are two rescalings for $v$. Then the sequence $B_{v,n}^{-1}\circ A_{v,n}$ is bounded.

Equivalently, if we identify the hyperbolic 3-space $\Hyp^3$ as the unit ball and $\PSL_2(\C) \cong \Isom(\Hyp^3)$, then 
$$
d_{\Hyp^3}(A_{v,n}({\bf 0}), B_{v,n}({\bf 0})) \text{ is bounded},
$$
where ${\bf 0} \in \Hyp^3$ is the center of the unit ball.
\end{lem}
\begin{proof}
Suppose this is not true. 
After passing to a subsequence, we may assume $C_n:=B_{v,n}^{-1}\circ A_{v,n}$ converges to a constant function $b$ compactly on $\hat \C-\{a\}$. 
The constant function $b\neq 0$ as $C_n(1), C_n(\infty) \notin B(0,1)$ for all $n$.
Since $C_n(0) = 0$ for all $n$, the hole $a = 0$. 
Therefore for any $x\neq b$, there is $N$ so that $x\in C_n(B(0,1))$ for $n\geq N$.
This is not possible as $1, \infty \notin C_n(B(0,1))$ for all $n$.
\end{proof}

Throughout the rest of this section, we shall fix a rescaling $A_{v,n}$ for each vertex $v\in \mathcal{V}$.

\subsection*{Carath\'eodory limits of pointed disks}
A {\em disk} is an open simply connected region in $\C$.
It is said to be {\em hyperbolic} if it is not $\C$.
For a sequence of pointed disks $(U_n, u_n)$, we say $(U_n, u_n)$ converges to $(U, u)$ in Carath\'eodory topology (see \S 5 in \cite{McM94}) if 
\begin{itemize}
\item $u_n \to u$;
\item For any compact set $K \subseteq U$, $K \subseteq U_n$ for all sufficiently large $n$;
\item For any open connected set $N$ containing $u$, if $N \subseteq U_n$ for all sufficiently large $n$, then $N \subseteq U$.
\end{itemize}

We say a sequence of proper holomorphic maps between pointed disks $f_n: (U_n, u_n) \longrightarrow (V_n, v_n)$ converges to $f:(U, u) \longrightarrow (V, v)$ if 
\begin{itemize}
\item $(U_n, u_n), (V_n,v_n)$ converge to $(U, u), (V, v)$ in Carath\'eodory topology;
\item For all sufficiently large $n$, $f_n$ converges to $f$ uniformly on compact subsets of $U$.
\end{itemize}

We have the following compactness result:
\begin{theorem}[\cite{McM94}, Theorem 5.2]\label{thm:cmpc}
The set of disks $(U_n, 0)$ containing $B(0,r)$ for some $r>0$ is compact in Carath\'eodory topology.
\end{theorem}

\subsection*{The rescalings are natural}
The above compactness theorem and the analysis of limits of proper maps (see Theorem 5.6 in \cite{McM94}) gives
\begin{lem}\label{lem:clr}
Let $v\in \mathcal{V}_s \subseteq \mathcal{V}$, with rescaling $A_{v,n}$. 
After passing to a subsequence, $A_{v,n}^{-1}((\U_{s,n}, v_n))$ converges in Carath\'eodory limit to a pointed hyperbolic disk $(\U_v, 0)$.
Moreover, the map 
$$
A_{F(v),n}^{-1} \circ f_n \circ A_{v,n}
$$ 
converges compactly to a proper holomorphic map
$$
R_v:\U_v \longrightarrow \U_{F(v)},
$$ 
conjugating to the corresponding rescaling limits for $\bp_n$
$$
\rl_v: \D_v \longrightarrow \D_{F(v)}.
$$
\end{lem}
\begin{proof}
The fact that $A_{v,n}^{-1}((\U_{s,n}, v_n))$ converges in Carath\'eodory limit follows immediately from Theorem \ref{thm:cmpc} as $B(0,1) \subseteq A_{v,n}^{-1}(\U_{s,n})$ and $A_{v, n}^{-1}(v_n) = 0$.
Since $d_{\U_n}(f_n(\psi_n(v)), \psi_n(F(v))) \leq K$, after passing to a subsequence, $A_{F(v),n}^{-1} \circ f_n \circ A_{v,n}(0)$ converges to a point in $\U_{F(v)}$.
Since $f_n: \U_{s,n} \longrightarrow \U_{\Phi(s),n}$ is proper,
$$
A_{F(v),n}^{-1} \circ f_n \circ A_{v,n} : A_{v,n}^{-1}(\U_{s,n}) \longrightarrow A_{F(v),n}^{-1}(\U_{\Phi(s),n})
$$ 
is also proper. 
Thus, the moreover part follows from case (3) of Theorem 5.6 in \cite{McM94}, and the fact $f_n$ is conjugate to $\bp_n$ on $s\times \D$.
\end{proof}

After passing to a subsequence, we may assume for different vertices $v, w\in \mathcal{V}$,
the rescalings $A_{v,n}^{-1} \circ A_{w,n}$ converge algebraically to either 
\begin{itemize}
\item A M\"obius transformation; or 
\item A constant map.
\end{itemize}
This defines an equivalence relation (cf. clusters in \S \ref{sec:DBS}):
\begin{defn}\label{defn:equivr}
Two vertices $v, w\in \mathcal{V}$ are said to be equivalent, denoted by $v\sim w$ if the rescalings $A_{v,n}^{-1} \circ A_{w,n}$ converge to a M\"obius transformation.
\end{defn}
\begin{lem}\label{lem:eq}
If $v \sim w$, then $F(v) \sim F(w)$.
\end{lem}
\begin{proof}
By Lemma \ref{lem:clr}, $A_{F(v),n}^{-1} \circ f_n\circ A_{v,n}$ and $A_{F(w),n}^{-1} \circ f_n\circ M_{w,n}$ both converge algebraically to a map of degree at least $1$.
Since $v\sim w$, then 
$$
A_{F(w),n}^{-1} \circ f_n\circ A_{v,n} = A_{F(w),n}^{-1} \circ f_n\circ A_{w,n} \circ (A_{w,n}^{-1} \circ A_{v,n})
$$ 
also converges algebraically to a map of degree at least $1$.
Thus $A_{F(v),n}^{-1} \circ A_{F(w),n}$ does not converge to a constant map.
Therefore, $F(v) \sim F(w)$.
\end{proof}

Let $\RP = \mathcal{V}/ \sim$ be the set of equivalence classes.
By Lemma \ref{lem:eq}, the dynamics $F$ on $\mathcal{V}$ induces a map on $\RP$. Abusing the notation, we use the same letter to denote the induced map
$$
F: \RP \longrightarrow \RP.
$$

\begin{defn}\label{defn:rld}
For each equivalence class $a\in \RP$, we choose a representative $v \in a$, and define the {\em rescaling} at $a$ by
$$
A_{a,n}:= A_{v,n} \in \PSL_2(\C).
$$
\end{defn}

\subsection{Rescaling trees}\label{subsec:rt}
We identify the hyperbolic 3-space $\Hyp^3$ as the unit ball in $\R^3$ and $\hat\C$ as the conformal boundary of $\Hyp^3$.
Denote ${\bf 0} \in \Hyp^3$ as the center of the unit ball.
Since $\PSL_2(\C) \cong \Isom(\Hyp^3)$, we denote $x_{a,n} = A_{a, n}({\bf 0}) \in \Hyp^3$.
We identify $\RP$ with the set of sequences of representatives
$$
\RP = \{(x_{a,n})_n \in \Hyp^3: a\in \RP\}.
$$
We also denote $\RP_n = \{x_{a, n} \in \Hyp^3: a\in \RP\}$.

\begin{figure}[ht]
  \centering
  \resizebox{0.6\linewidth}{!}{
    \def\svgwidth{\columnwidth}
\begingroup%
  \makeatletter%
  \providecommand\color[2][]{%
    \errmessage{(Inkscape) Color is used for the text in Inkscape, but the package 'color.sty' is not loaded}%
    \renewcommand\color[2][]{}%
  }%
  \providecommand\transparent[1]{%
    \errmessage{(Inkscape) Transparency is used (non-zero) for the text in Inkscape, but the package 'transparent.sty' is not loaded}%
    \renewcommand\transparent[1]{}%
  }%
  \providecommand\rotatebox[2]{#2}%
  \newcommand*\fsize{\dimexpr\f@size pt\relax}%
  \newcommand*\lineheight[1]{\fontsize{\fsize}{#1\fsize}\selectfont}%
  \ifx\svgwidth\undefined%
    \setlength{\unitlength}{680.31496063bp}%
    \ifx\svgscale\undefined%
      \relax%
    \else%
      \setlength{\unitlength}{\unitlength * \real{\svgscale}}%
    \fi%
  \else%
    \setlength{\unitlength}{\svgwidth}%
  \fi%
  \global\let\svgwidth\undefined%
  \global\let\svgscale\undefined%
  \makeatother%
  \begin{picture}(1,0.875)%
    \lineheight{1}%
    \setlength\tabcolsep{0pt}%
    \put(0,0){\includegraphics[width=\unitlength,page=1]{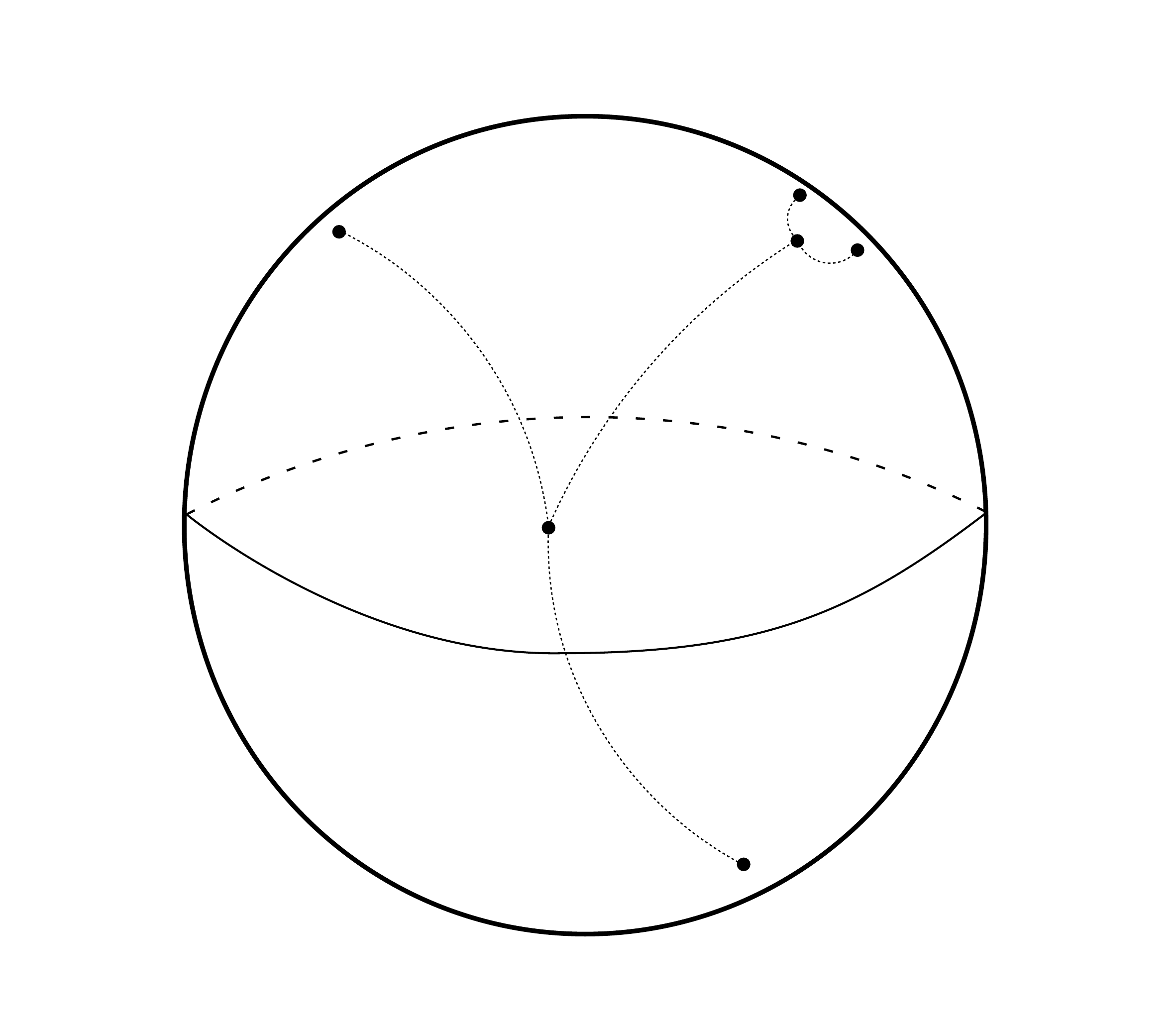}}%
    \put(0.48602708,0.4214963){\color[rgb]{0,0,0}\makebox(0,0)[lt]{\lineheight{1.25}\smash{\begin{tabular}[t]{l}{\tiny$x_{a,n}=A_{a,n}(\bf 0)=\Psi_n(a)$}\end{tabular}}}}%
    \put(0.31039569,0.6816699){\color[rgb]{0,0,0}\makebox(0,0)[lt]{\lineheight{1.25}\smash{\begin{tabular}[t]{l}{\tiny$x_{b,n}=A_{b,n}(\bf 0)=\Psi_n(b)$}\end{tabular}}}}%
    \put(0.63010057,0.15911776){\color[rgb]{0,0,0}\makebox(0,0)[lt]{\lineheight{1.25}\smash{\begin{tabular}[t]{l}{\tiny$x_{c,n}=A_{c,n}(\bf 0)=\Psi_n(c)$}\end{tabular}}}}%
    \put(0.69073404,0.67836261){\color[rgb]{0,0,0}\makebox(0,0)[lt]{\lineheight{1.25}\smash{\begin{tabular}[t]{l}{\tiny$x_{d,n}=A_{d,n}(\bf 0)=\Psi_n(d)$}\end{tabular}}}}%
    \put(0.69563363,0.75263021){\color[rgb]{0,0,0}\makebox(0,0)[lt]{\lineheight{1.25}\smash{\begin{tabular}[t]{l}{\LARGE $\Hyp^3$}\end{tabular}}}}%
    \put(0,0){\includegraphics[width=\unitlength,page=2]{RescalingTree.pdf}}%
  \end{picture}%
\endgroup%

  }
  \caption{An illustration of the rescaling tree in $\Hyp^3$.}
  \label{fig:RescalingTree}
\end{figure}

The same construction for quasi-invariant trees $\mathcal{T}_{s,n} \subseteq \Hyp^2$ in \S \ref{sec:DBS} generalizes immediately to $\Hyp^3$ to give a sequence of {\em rescaling trees} $\RT_n$ with vertex set $\RV_n$. 
After passing to a subsequence, we assume we have the markings
$$
\Psi_n: \RT \longrightarrow \RT_n,
$$
with vertex set $\RV$.
By construction, we have
\begin{align}\label{eqn:di}
d_{\Hyp^3}(\Psi_n(a), \Psi_n(b)) \to \infty
\end{align}
for $a \neq b \in \RV$.
By construction, $\RP\subseteq \RV$ and 
any point in $\RV - \RP$ is a branch point of $\RT$.

We can extend the definition of rescalings for $\RV - \RP$ as follows.
Let $a\in \RV-\RP$, a sequence $A_{a,n} \in \PSL_2(\C) \cong \Isom(\Hyp^3)$ is defined to be a {\em resacling} at $a$ if 
$$
A_{a,n}({\bf 0}) = \Psi_n(a).
$$
Note that different choices of rescalings at $a$ are differed by pre-composing with a rotation that fixes ${\bf 0}$.
For the rest of the section, we shall fix such a rescaling for each vertex $a\in \RV$.

Let 
$\hat\C^\RV := \bigcup_{a\in\RV} \hat \C_a$
be the disjoint union of Riemann spheres.
We remark that $\hat \C_a$ is interpreted as the Riemann sphere under the rescaling coordinate given by $A_{a,n}$.

Denote $\overline{\B} = \Hyp^3 \cup \hat \C$, and similarly $\overline{\B}_a = \Hyp^3_a \cup \hat \C_a$ for $a\in \RV$.
We say a sequence $z_n \in \overline{\B}$ converges to $z\in \overline{\B}_a$ {\em in $a$-coordinate} or {\em with respect to the rescaling at $a$}, denoted by $z_n \to_a z$ or $z = \lim_a z_n$ if
$$
\lim_{n\to\infty} A_{a,n}^{-1}(z_n) = z.
$$

\begin{figure}[ht]
  \centering
  \resizebox{0.6\linewidth}{!}{
    \def\svgwidth{\columnwidth}
\begingroup%
  \makeatletter%
  \providecommand\color[2][]{%
    \errmessage{(Inkscape) Color is used for the text in Inkscape, but the package 'color.sty' is not loaded}%
    \renewcommand\color[2][]{}%
  }%
  \providecommand\transparent[1]{%
    \errmessage{(Inkscape) Transparency is used (non-zero) for the text in Inkscape, but the package 'transparent.sty' is not loaded}%
    \renewcommand\transparent[1]{}%
  }%
  \providecommand\rotatebox[2]{#2}%
  \newcommand*\fsize{\dimexpr\f@size pt\relax}%
  \newcommand*\lineheight[1]{\fontsize{\fsize}{#1\fsize}\selectfont}%
  \ifx\svgwidth\undefined%
    \setlength{\unitlength}{1275.59055118bp}%
    \ifx\svgscale\undefined%
      \relax%
    \else%
      \setlength{\unitlength}{\unitlength * \real{\svgscale}}%
    \fi%
  \else%
    \setlength{\unitlength}{\svgwidth}%
  \fi%
  \global\let\svgwidth\undefined%
  \global\let\svgscale\undefined%
  \makeatother%
  \begin{picture}(1,0.46666667)%
    \lineheight{1}%
    \setlength\tabcolsep{0pt}%
    \put(0,0){\includegraphics[width=\unitlength,page=1]{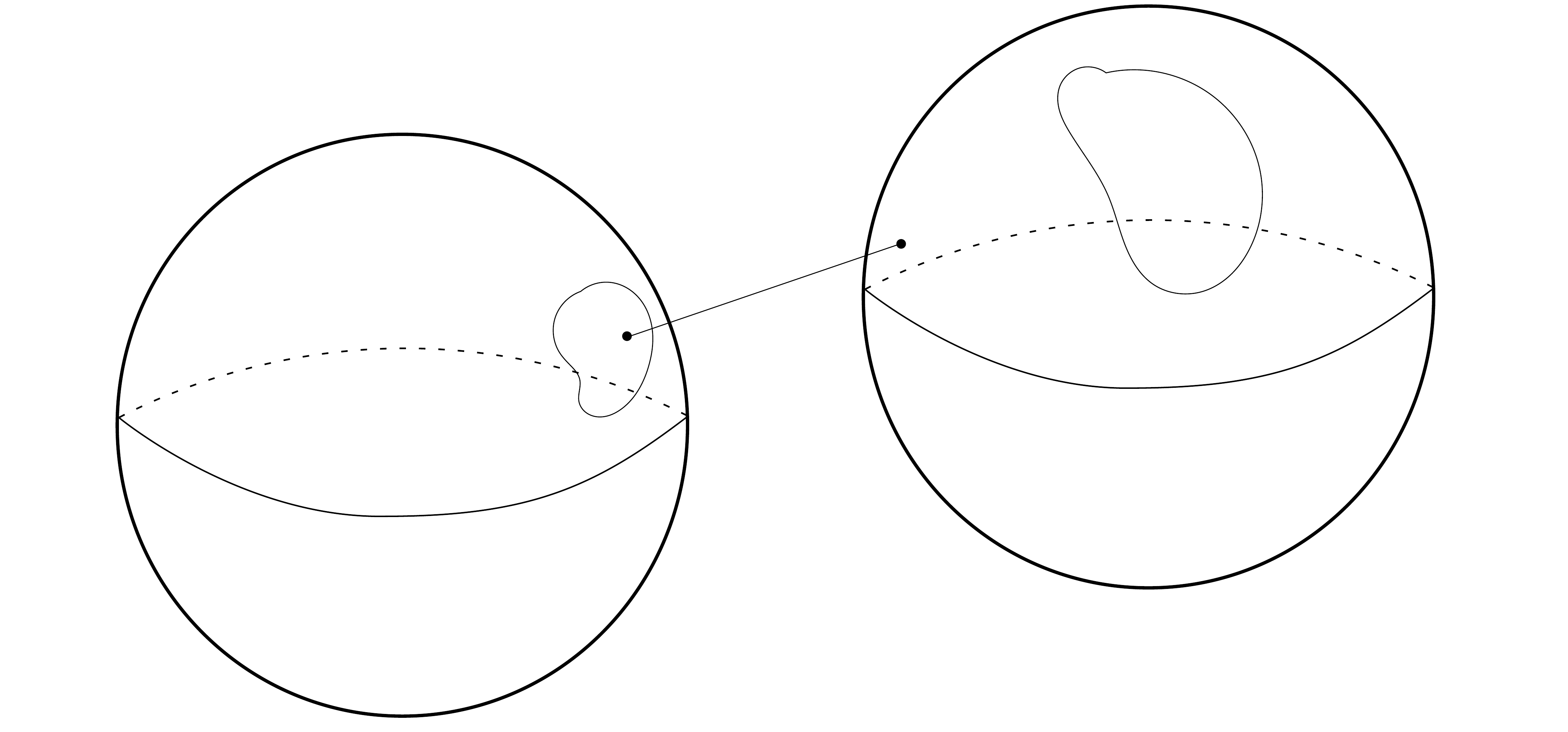}}%
    \put(0.30926853,0.38435185){\color[rgb]{0,0,0}\makebox(0,0)[lt]{\lineheight{1.25}\smash{\begin{tabular}[t]{l}{\LARGE $\hat\C_b$}\end{tabular}}}}%
    \put(0.83114577,0.07004893){\color[rgb]{0,0,0}\makebox(0,0)[lt]{\lineheight{1.25}\smash{\begin{tabular}[t]{l}{\LARGE $\hat\C_a$}\end{tabular}}}}%
    \put(0.40099074,0.2403009){\color[rgb]{0,0,0}\makebox(0,0)[lt]{\lineheight{1.25}\smash{\begin{tabular}[t]{l}{\tiny $x_a \in\Xi_b$}\end{tabular}}}}%
    \put(0.57132589,0.3218119){\color[rgb]{0,0,0}\makebox(0,0)[lt]{\lineheight{1.25}\smash{\begin{tabular}[t]{l}{\tiny $x_b \in\Xi_a$}\end{tabular}}}}%
    \put(0.73789354,0.37788426){\color[rgb]{0,0,0}\makebox(0,0)[lt]{\lineheight{1.25}\smash{\begin{tabular}[t]{l}$U$\end{tabular}}}}%
    \put(0.3710046,0.25206017){\color[rgb]{0,0,0}\makebox(0,0)[lt]{\lineheight{1.25}\smash{\begin{tabular}[t]{l}$V$\end{tabular}}}}%
  \end{picture}%
\endgroup%

  }
  \caption{An illustration of $U \subseteq_\RT V$.}
  \label{fig:Subset}
\end{figure}
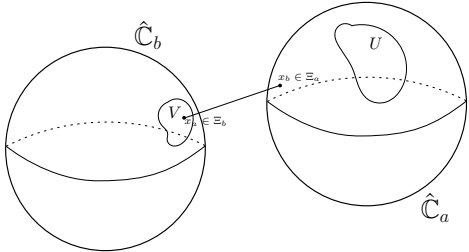

By Equation \ref{eqn:di}, $A_{a,n}^{-1} \circ A_{b,n}$ converges to a constant map $x_b$ for $a\neq b \in \RV$.
Thus, we can associate the point $x_b \in \hat \C_a$ to $b$.
It is interpreted that the Riemann sphere $\hat \C_b$ converges to $x_b$ in the rescaling coordinate $\hat \C_a$.
We denote 
$$
\Xi_a:= \bigcup_{b\neq a} x_b \subseteq \hat\C_a
$$
as the {\em singular set} at $a$ and $\Xi := \bigcup_{a\in \RV} \Xi_a \subseteq \hat \C^\RV$ as the {\em singular set} (see Figure \ref{fig:Subset}).

The following lemma follows immediately from the definition.
\begin{lem}\label{lem:h3l}
Let $a \neq b \in \RV$. Then $\Psi_n(b) \to_a x_b \in \hat\C_a$.
\end{lem}
\begin{proof}
Since $A_{a,n}^{-1} \circ A_{b,n}$ converges to a constant map $x_b$, $A_{a,n}^{-1} \circ A_{b,n} ({\bf 0}) \to x_b$.
Therefore, $A_{b,n}({\bf 0}) = \Psi_n(b) \to_a x_b$ by definition. 
\end{proof}

By the construction, the angle 
$$
\angle \Psi_n(a) \Psi_n(b) \Psi_n(c)
$$ 
is uniformly bounded below from $0$ for any distinct triple $a, b, c\in \RV$.
Thus, by Lemma \ref{lem:h3l}, the singular set $\Xi_a$ is in bijective correspondence with the tangent space $T_a \RT$ at $a$.
We denote this correspondence by
$$
\xi_a: T_a\RT \longrightarrow \Xi_a.
$$

\begin{defn}\label{defn:subset}
Let $a\neq b \in \RV$, and $x_b = \xi_a(v_b) \in \Xi_a \subseteq \hat\C_a$ where $v_b \in T_a \RT$ is the tangent vector at $a$ in the direction of $b$.
Similarly, let $x_a = \xi_b(v_a) \in \Xi_b \subseteq \hat\C_b$.
Let $U \subseteq \hat\C_a$, and $V \subseteq \hat\C_b$ be open sets. 
We define $U \subseteq_\RT V$ if $x_b \notin U$ and $x_a \in V$.
\end{defn}
This definition is natural as explained in the next lemma. The proof follows directly from our definition (see Figure \ref{fig:Subset}).
\begin{lem}
If $U \subseteq \hat\C_a$, and $V \subseteq \hat\C_b$ with $U \subseteq_\RT V$, then for sufficiently large $n$, we have
$$
A_{a,n}(U) \subseteq A_{b,n}(V).
$$
\end{lem}

\subsection{Rescaling rational maps}\label{subsec:rrm}
We now show the rescalings give abundant of rescaling limits of rational maps.
\begin{lem}\label{lem: extmap}
Let $a\in \RV$. 
After passing to a subsequence, there exists a unique $b\in \RV$ so that
$$
A_{b,n}^{-1} \circ f_n \circ A_{a,n}
$$
converges algebraically to a rational map $R_a = R_{a\to b}$ of degree at least $1$.

Moreover, if $a\in \RV-\RP$, then $b$ is a branch point.
\end{lem}
\begin{proof}
The uniqueness follows from Equation \ref{eqn:di}. 

To prove existence, if $a\in \RP \subseteq \RV$, then we let $b = F(a) \in \RP$.
It follows from Lemma \ref{lem:clr} that $A_{F(a),n}^{-1} \circ f_n \circ A_{a,n}$ converges to algebraically to a rational map of degree $\geq 1$.

If $a\in \RV - \RP$, then $a$ is a branch point.
There exists $B_n \in \PSL_2(\C)$ (see Lemma 2.1 in \cite{DeMF14}) so that
$$
B_n^{-1} \circ f_n \circ A_{a,n}
$$
converges algebraically to a rational map $R$ of degree $\geq 1$.
After passing to a subsequence, we assume $B_n^{-1}\circ A_{b,n}$ converges algebraically for all $b\in \RV$.

We claim that for each $t\in \Xi_a$, there exists $b\in \RP$ so that $B_n^{-1}\circ A_{b, n}$ converging to the constant map $x\mapsto R(t)$, or equivalently, $\lim_{n\to\infty} B_n^{-1}(\Psi_n(b)) = R(t)$ in $\overline{\B} = \Hyp^3 \cup \hat \C$.

\begin{figure}[ht]
  \centering
  \resizebox{0.6\linewidth}{!}{
    \def\svgwidth{\columnwidth}
    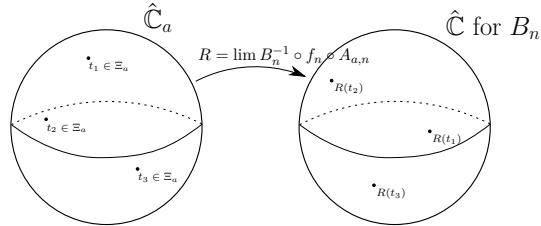

  }
  \caption{For each $t\in \Xi_a$, there exists some vertex $b\in \RP$ so that $B_n^{-1}(\Psi_n(b))$ converges to $R(t)$.}
  \label{fig:BranchClaim}
\end{figure}

\begin{proof}[Proof of the claim]
We consider two cases.

Case (1): If there exists a sequence of critical values $z_n$ of $f_n$ such that
$B_n^{-1}(z_n) \to R(t)$.
Let $U_n$ be the Fatou component of $f_n$ containing $z_n$.
After passing to a subsequence, the sequence $(z_n)$ corresponds to some vertex $b \in \RP$.
Thus, $A_{b,n}^{-1}((U_n, z_n))$
converges in Carath\'eodory limit to some pointed hyperbolic disk $(U, u)$ in $\hat\C_{b}$.

Suppose $B_n^{-1}\circ A_{b, n}$ does not converge to $R(t)$. Then $B_n^{-1}\circ A_{b, n}(U)$ contains a neighborhood of $R(t)$ for sufficiently large $n$. Thus, $B_n^{-1}(U_n)$ contains an open set $V$ containing $R(t)$ for sufficiently large $n$.
Since $B_n^{-1} \circ f_n \circ A_{a,n}$ converges algebraically to $R$ with degree $\geq 1$, there exists an open set $t\in V'$ so that $A_{a,n} (V')$ is contained in a Fatou component of $f_n$ for all sufficiently large $n$, contradicting $t\in \Xi_a$.

Case(2): Otherwise, there exists a small disk $U$ of $R(t)$ so that $B_n(U)$ contains no critical values of $f_n$ for all sufficiently large $n$.
So $f_n^{-1}(B_n(U))$ consists of $d$ disks for all sufficiently large $n$.
Thus no critical points of $B_n^{-1} \circ f_n \circ A_{a,n}$ converges to $t$.
By Lemma \ref{lem:critc}, $t$ is not a hole of $R$.
Let $b' \in \RP$ with $A_{a,n}^{-1}\circ A_{b', n}$ converging to $t$.
Let $b = F(b')$.
Since $t$ is not a hole for $R$, $B_n^{-1}\circ A_{b, n}$ converges to the constant map $R(t)$.
\end{proof}

Since $a\in \RV-\RP$, if $c_n$ is a critical point of $f_n$, then $\lim_a c_n \in \Xi_a$.
Therefore, $\Xi_a$ contains all critical points of $R$.
Since $a$ is a branch point, $|\Xi_a| \geq 3$.
Therefore, $|R(\Xi_a)| \geq 3$ as $\Xi_a$ contains all critical points of $R$.
By the claim, there exist at least three vertices $b_1, b_2, b_3 \in \RP$ with $d_{\Hyp^3}(B_n({\bf 0}), \Psi_n(b_i)) \to \infty$ so that the angles
$$
\angle \Psi_n(b_i) B_n({\bf 0}) \Psi_n(b_j)
$$
are bounded below from $0$ for all $1\leq i < j \leq 3$.
Thus $B_n({\bf 0})$ is within a bounded distance of a branch point of $\RT_n$, and the lemma follows.
\end{proof}

\subsection{Rescaling tree map}\label{subsec:rtm}
We define $F: \RV \longrightarrow \RV$ by setting $F(a)$ as the unique vertex $b$ in Lemma \ref{lem: extmap} which extends the map $F: \RP\subseteq \RV \longrightarrow \RP$.

We define the {\em rescaling tree map} $F: \RT \longrightarrow \RT$ for $([f_n])_n$ by extending continuously on any edge $[a,b]$ to the arc $[F(a), F(b)]$.

By Lemma \ref{lem: extmap},
$A_{F(a),n}^{-1} \circ f_n \circ A_{a,n}$
converges algebraically to a rational map $R_a=R_{a\to F(a)}$ of degree at least $1$.
We define the rescaling limit as follows (cf. Definition \ref{defn:rlm}):
\begin{defn}
We define $R_a = \lim A_{F(a),n}^{-1} \circ f_n \circ A_{a,n}$ as {\em rescaling limit} at $a$.
We define the union of the rescaling limits as
$$
R: \hat\C^\RV\longrightarrow \hat\C^\RV.
$$

The local degree at a vertex $a\in \RV$ is denoted by
$$
\delta(a) := \deg (R_a).
$$
\end{defn}

For future reference, we have the following lemma controlling the holes of the rescaling limit $R_a$.
\begin{lem}\label{lem:hcs}
Let $a\in \RV$. Then the holes of $R_a$ are contained in $\Xi_a$.
\end{lem}
\begin{proof}
Let $x \in \hat\C_a$ be a hole of $R_a$. 
Since $\deg(R_a) \geq 1$, by Lemma \ref{lem:critc}, there exists critical points $c_n$ of $f_n$ with $c_n \to_a x$.
After passing to a subsequence, the sequence $c_n$ corresponds to a point $b \in \RP\subset \RV$.
If $b = a$, $x$ is not a hole, so $b\neq a$.
Then the corresponding point of $b$ on $\hat \C_a$ is $x$, so $x\in \Xi_a$.
\end{proof}

\subsection*{Dynamics on edges}
Let $E = [a,b]$ be an edge of $\RT$.
In the following, we shall associate it with a sequence of annuli $\mathcal{A}_{E,n}$ and define the local degree of $E$.

Let $x_b \in \Xi_a \subseteq \hat\C_a$ and $x_a\in \Xi_b \subseteq \hat\C_b$ be the points associated to $b$ and $a$ respectively.
Choose a small closed curve $C_a \subseteq \hat\C_a$ around $x_b$.
We assume $C_a$ bounds no holes nor critical points of $R_a$ other than possibly $x_b$ and $R_a: C_a \longrightarrow R_a(C_a)$ is a covering map of degree $\deg_{x_b}(R_a)$.
Similarly, we define $C_b \subseteq \hat\C_b$ around $x_a$.

\begin{figure}[ht]
  \centering
  \resizebox{0.6\linewidth}{!}{
    \def\svgwidth{\columnwidth}
    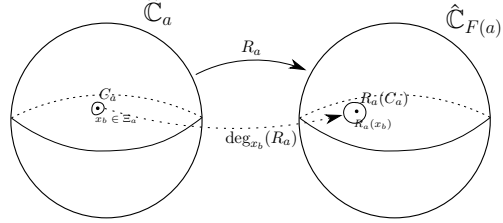

  }
  \caption{An illustration of the construction of $C_a$.}
  \label{fig:EdgeD}
\end{figure}

Since $A_{F(a), n}^{-1} \circ f_n \circ A_{a, n}$ converges uniformly to $R_a$ in a neighborhood of $C_a$, we can find a sequence of closed curves $C_{a,n}$ such that
\begin{itemize}
\item $f_n: C_{a,n} \longrightarrow f_n(C_{a,n})$ is a covering of degree $\deg_{x_b}(R_a)$;
\item $A_{a,n}^{-1}(C_{a,n})$ converges in Hausdorff topology to $C_a$.
\end{itemize}
Similarly, let $C_{b,n}$ be the corresponding closed curves for $C_b$.
Note that $C_{a, n}$ and $C_{b,n}$ are nested for sufficiently large $n$ and bounds an annulus $\mathcal{A}_{E,n}\subseteq \hat\C$.
We call $\mathcal{A}_{E,n}$ a sequence of annuli associated to $E$.

We claim that $\mathcal{A}_{E,n}$ contains no critical points of $f_n$ for sufficiently large $n$.
Indeed, otherwise, after passing to a subsequence, let $c_n \in \mathcal{A}_{E,n}$, and let $c\in \RP$ correspond to the sequence $(c_n)$.
Consider the projection 
$$
\proj_{[\Psi_n(a), \Psi_n(b)]}(\Psi_n(c))
$$ 
of $\Psi_n(c) \in \Hyp^3$ onto the geodesic $[\Psi_n(a), \Psi_n(b)]$.
Since we assume $C_a$ bounds no holes nor critical points other than possibly $\xi_a(v_b)$ and similarly for $C_b$, 
$$
d_{\Hyp^3}(\proj_{[\Psi_n(a), \Psi_n(b)]}(\Psi_n(c)), \partial [\Psi_n(a), \Psi_n(b)]) \to \infty.
$$
This contradicts that $[a, b]$ is an edge of $\RT$.

Therefore, $f_n: \mathcal{A}_{E,n} \longrightarrow f_n(\mathcal{A}_{E,n})$ is a covering map.
We define the local degree at $E$
$$
\delta(E) := \deg_{x_b}(R_a)= \deg_{x_a}(R_b),
$$
as the degree of this covering map.

We have the following modulus estimate.
\begin{prop}\label{prop: moduluse}
Let $\mathcal{A}_{E,n}$ associated to an edge $E = [a,b]$ of $\RT$.
There exists a constant $K$ such that
the modulus 
$$
|m(\mathcal{A}_{E,n})-\frac{d_{\Hyp^3}(\Psi_n(a), \Psi_n(b))}{2\pi}| \leq K,
$$
and 
$$
|m(f_n(\mathcal{A}_{E,n}))-\frac{d_{\Hyp^3}(\Psi_n(F(a)), \Psi_n(F(b)))}{2\pi}| \leq K.
$$
\end{prop}
\begin{proof}
Pre-composing $A_{a,n}$ and $A_{b,n}$ with rotations if necessary, we may assume
$$
A_{a,n}^{-1} \circ A_{b,n}(z) = \frac{z}{\rho_n},
$$
where $\rho_n = e^{d_{\Hyp^3}(\Psi_n(a), \Psi_n(b))} \to \infty$.
Then $x_b = 0\in \hat\C_a$.
Choose small radii $0< r_1 < r_2$ so that the curve $C_a$ is contained in the annulus $B(0, r_2) - \overline{B(0, r_1)}$.
Similarly, $A_{b,n}^{-1} \circ A_{a,n}(z) = \rho_n z$, so $x_a = \infty \in \hat\C_b$.
Choose large radii $0< r_1' < r_2'$ so that the curve $C_b$ is contained in the annulus $B(0, r_2') - \overline{B(0, r_1)}$.

Therefore, for sufficiently large $n$, $A_{a,n}^{-1} (C_{a,n})$ is contained in $B(0, r_2) - \overline{B(0, r_1)}$ and $A_{b,n}^{-1} (C_{b,n})$ is contained in the annulus $B(0, r_2') - \overline{B(0, r_1)}$.
So the annulus $\mathcal{A}_{E,n}$ contains the annulus $\mathcal{A}^1_n$ bounded by $A_{a,n}(\partial B(0, r_1))$ and $A_{b,n}(\partial B(0, r_2'))$ and is contained in the annulus $\mathcal{A}^2_n$ bounded by $A_{a,n}(\partial B(0, r_2))$ and $A_{b,n}(\partial B(0, r_1'))$.
Moreover, the inclusions $\mathcal{A}^1_n\xhookrightarrow{} \mathcal{A}_{E,n} \xhookrightarrow{} \mathcal{A}^2_n$ induce isomorphisms on the fundamental groups, so $m(\mathcal{A}^1_n) \leq m(\mathcal{A}_{E,n})\leq m(\mathcal{A}^2_n)$.

The modulus 
$$
m(\mathcal{A}^1_n) = \frac{\log(\rho_nr_1/r_2')}{2\pi} = \frac{d_{\Hyp^3}(\Psi_n(a), \Psi_n(b))}{2\pi} + \frac{\log r_1 - \log r_2'}{2\pi}
$$
and
$$
m(\mathcal{A}^2_n) = \frac{\log(\rho_nr_2/r_1')}{2\pi}= \frac{d_{\Hyp^3}(\Psi_n(a), \Psi_n(b))}{2\pi} + \frac{\log r_2 - \log r_1'}{2\pi}
$$
Therefore, $m(\mathcal{A}_{E,n}) = \frac{d_{\Hyp^3}(\Psi_n(a), \Psi_n(b))}{2\pi} + O(1)$.
The same argument gives the estimate of $m(f_n(\mathcal{A}_{E,n}))$.
\end{proof}

The above modulus estimate gives the following corollaries.
\begin{cor}\label{cor:inje}
Let $E = [a,b]$ be an edge of $\RT$. Then
$$
d_{\Hyp^3}(\Psi_n(F(a)), \Psi_n(F(b))) = \delta(E) d_{\Hyp^3}(\Psi_n(a), \Psi_n(b)) + O(1).
$$
\end{cor}
\begin{proof}
This follows from $m(f_n(\mathcal{A}_{E,n})) = \delta(E) m(\mathcal{A}_{E,n})$ and the estimates in Proposition \ref{prop: moduluse}.
\end{proof}

\begin{cor}\label{cor:inje}
If $[a,b]$ is an edge of $\RT$, then $F(a) \neq F(b)$.
\end{cor}

Corollary \ref{cor:inje} gives that the map $F$ is injective on edges, and allows to define the tangent map $DF_a: T_a \RT \longrightarrow T_{F(a)}\RT$.
The following compatibility result also follows from the modulus estimate.
\begin{cor}\label{cor:comtm}
Let $a \in \RT$. Then
$$
R_a \circ \xi_a = \xi_{F(a)} \circ DF_a.
$$
\end{cor}

As an immediate corollary of the above, we have
\begin{cor}\label{cor:sinv}
The singular set $\Xi \subseteq \hat\C^\RV$ is invariant under $R$.
\end{cor}

It now follows from our construction that
\begin{cor}\label{cor:rm}
	The map $(F,R)$ is a degree $d$ rational map on the tree of Riemann spheres $(\RT, \hat\C^\RV)$.
\end{cor}

\subsection*{Matrix encoding}
Index the set of edges of $\RT$ by
$\{E_1,..., E_k\}$.
We define the following two matrices to encode the dynamics $F: \RT\longrightarrow \RT$.
\begin{itemize}
\item (Markov matrix): $M_{i,j} = \begin{cases} 1 &\mbox{if } E_i \subseteq F(E_j) \\ 
0 & \mbox{otherwise } \end{cases}$
\item (Degree matrix): $D_{i,j} = \begin{cases} \delta(E_i) &\mbox{if } i = j \\ 
0 & \mbox{otherwise } \end{cases}$
\end{itemize}
We remark that the matrix $D^{-1} M$ is the analogue of the Thurston's matrix for branched coverings (see \S \ref{sec:cscrm} for more details).
The following proposition says that if the tree $\RT$ is non-trivial, then there is always a `canonical Thurston's obstruction'.

\begin{prop}\label{prop:mdto}
Let $M$ and $D$ be the Markov matrix and the degree matrix respectively.
If $\RT$ is not trivial, then there exists a non-negative vector $\vec{v} \neq \vec{0}$ so that
$$
M\vec v = D \vec v.
$$
\end{prop}
\begin{proof}
Let 
$\rho_n = \max_{i=1,..., k} l_{\Hyp^3} (\Psi_n(E_i)) \to \infty$,
where $l_{\Hyp^3} (\Psi_n(E_i))$ is the hyperbolic length of the geodesic segment $\Psi_n(E_i) \subseteq \Hyp^3$.
After passing to a subsequence, we assume the limit $\lim_{n\to\infty} l_{\Hyp^3} (\Psi_n(E_i))/\rho_n$ exists for all $i=1,..., k$.
Define $\vec{v}$ whose $i$-th entry is $\lim_{n\to\infty} l_{\Hyp^3} (\Psi_n(E_i))/\rho_n$. 
Then $\vec{v}$ is non-negative and $\vec{v} \neq \vec{0}$.

If suffices to check $M\vec v = D \vec v$.
If $a, b \in \RV$ are connected by a sequence of edges $E_{i_1} \cup E_{i_2} \cup... E_{i_m}$, 
since the angles between different incident edges at a vertex of $\RT_n$ are uniformly bounded below from $0$, we have
$$
d_{\Hyp^3}(\Psi_n(a), \Psi_n(b)) = \sum_{j=1}^m l_{\Hyp^3}(\Psi_n(E_{i_j})) + O(1).
$$ 
Thus, by Corollary \ref{cor:inje}, if $F(E_i) = E_{i_1} \cup E_{i_2} \cup... E_{i_m}$, then 
$$
\delta(E_i) l_{\Hyp^3} (\Psi_n(E_i)) = \sum_{j=1}^m l_{\Hyp^3}(\Psi_n(E_{i_j})) + O(1).
$$
Dividing both sides by $\rho_n$ and taking limits, we conclude the result.
\end{proof}

It is worth noting that from the proof, the vector $\vec{v}$ has a concrete geometric interpretation as `rescaled edge lengths'.

\subsection{Geometrically finite limits}\label{subsec:gfl}
Let $R: \hat\C^\RV \longrightarrow \hat\C^\RV$ be the rescaling rational map.
Recall that the Fatou set $\Omega(R) \subset \hat\C^\RV$ is the largest open set for which the iterates $\{R^n|_{\Omega}: n \geq 1\}$ form a normal family, and the Julia set $\mathcal{J}(R)$ is the complement of the Fatou set.

Let $a\in \RV$ be a periodic point of period $q$. 
If $\deg(R_a^q) \geq 2$, then $\Omega_a \subseteq \hat\C_a$ and $\mathcal{J}_a$ are the Fatou set and the Julia set of $R_a^q$.
If $\deg(R_a^q) =1$, then the Julia set $\mathcal{J}_a \subseteq \hat\C_a$ is the set of the repelling fixed point or parabolic fixed point if $R_a^q$ is hyperbolic or parabolic, and is empty if $R_a^q$ is elliptic.

If $a\in \RV$ is an aperiodic point with pre-period $l$, then the Fatou set and Julia set are the pull back of the periodic one: $\Omega_a = (R_a^l)^{-1}(\Omega_{F^l(a)})$ and $\mathcal{J}_a = (R_a^l)^{-1}(\mathcal{J}_{F^l(a)})$.

A rational map $R: \hat\C^\RV \longrightarrow \hat\C^\RV$ is {\em geometrically finite} if every critical point in the Julia set is pre-periodic.
We will now show the rescaling rational map is geometrically finite. 

To prepare for the proof, we need the following two lemmas.
The first estimate follows from the Schwarz lemma and Koebe's theorem.
\begin{lem}\label{lem:hme}
Let $U$ be a hyperbolic disk in $\C$ with hyperbolic metric $\rho_U |dz|$. Then for $z\in U$,
$$
\frac{1}{2d_{\R^2}(z, \partial U)} \leq \rho_U(z) \leq \frac{2}{d_{\R^2}(z, \partial U)}.
$$
\end{lem}

We also need the following lemma, which is useful in its own right.
\begin{lem}\label{lem:key}[cf. Lemma 6.5 in \cite{Luo21}]
Let $f_n: \hat \C \longrightarrow \hat \C$ be a sequence of rational maps converging algebraically to $f: \hat\C \longrightarrow \hat \C$ with $\deg(f) \geq 1$.
Let $U_n$ be a sequence of invariant hyperbolic disks for $f_n$, and $x_n \in U_n$.
After passing to a subsequence, let $x = \lim_{n\to\infty} x_n$.
If there exists $K$ so that $d_{U_n}(x_n, f_n(x_n)) \leq K$ for all $n$, then
\begin{itemize}
\item $x$ is a hole for $f$; or
\item $x$ is fixed by $f$; or
\item $(U_n, x_n)$ converges in Carath\'eodory topology to $(U, x)$. Thus $f_n: (U_n, x_n) \longrightarrow (U_n, f_n(x_n))$ converges to $f:(U, x)\longrightarrow (U, f(x))$.
\end{itemize}
\end{lem}
\begin{proof}
We may assume $U_n \subseteq \C$ for all sufficiently large $n$ and $x\in \C$.
Suppose $x$ is not a hole nor fixed by $f$, we claim there exists a Euclidean ball of radius $r$ so that $B(x_n, r) \subseteq U_n$ for all sufficiently large $n$.
Suppose not, then the hyperbolic metric $\rho_{U_n}(x_n)|dz|$ at $x_n$ is going to infinity by Lemma \ref{lem:hme}.
Since $x$ is not a hole, $f_n$ converges to $f$ uniformly near $x$. Thus, for sufficiently large $n$, the Euclidean distance between $x_n$ and $f_n(x_n)$ is bounded below from $0$.
Thus, $d_{U_n}(x_n, f_n(x_n))$ is unbounded which is a contradiction.
Therefore, after passing to a subsequence, the pointed disks $(U_n, x_n)$ converge in Carath\'eodory topology to $(U, x)$ by Theorem \ref{thm:cmpc}.
\end{proof}

\begin{prop}\label{prop:gfl}
The rescaling rational map $R: \hat\C^\RV \longrightarrow \hat\C^\RV$ is geometrically finite.
\end{prop}
\begin{proof}
Let $a\in \RV$ and $c \in \hat\C_a$ be a critical point. If its iterates contain a singular point, then it is pre-periodic by Corollary \ref{cor:sinv}.

Otherwise, let $c_n$ be a sequence of critical points for $f_n$ with $c_n \to_a c$.
Since $([f_n])$ is quasi post-critically finite, after passing to a subsequence, we can find an iterate $x_n = f_n^l(c_n)$ which is quasi periodic.
More precisely, let $U_n$ be a Fatou component of $f_n$ containing $x_n$. 
Then there exists a quasi periodic $q$ and a constant $K$ so that
$$
d_{U_n}(x_n, f_n^q(x_n))\leq K.
$$
Let $b = F^l(a)$ and $x = R^l(c) \in \hat \C_b$.
Since iterates of $c$ do not intersect the singular set, by Lemma \ref{lem:hcs}, the iterates of $c$ avoid the holes of $R$.
Therefore, $x_n \to_b x$ and $x$ is not a hole of $R^q$.

If $x$ is fixed by $R^q$, then $c$ is pre-periodic.
Otherwise, applying Lemma \ref{lem:key} to $A_{b,n}^{-1} \circ f_n^q \circ A_{b,n}$, we conclude that 
$$
A_{b,n}^{-1} \circ f_n^q \circ A_{b,n}: (A_{b,n}^{-1}(U_n), A_{b,n}^{-1}(x_n)) \longrightarrow (A_{b,n}^{-1}(U_n), A_{b,n}^{-1}(f_n^q(x_n)))
$$
converges in Carath\'eodory topology to
$$
R^q : (U, x) \longrightarrow (U, R^q(x)).
$$
Therefore $x\in U$ is in the Fatou set of $R$. So $c$ is contained in the Fatou set as well.

Run the above argument for all critical points, we conclude $R$ is geometrically finite.
\end{proof}

\begin{proof}[Proof of Theorem \ref{thm:crm}]
Let $F: (\RT, \RV) \longrightarrow (\RT, \RV)$ be the rescaling tree map and $R: \hat\C^\RV \longrightarrow \hat\C^\RV$ be the rescaling rational map. Then Corollary \ref {cor:rm} gives that $(F, R)$ is a degree $d$ rational map on $(\RT, \hat\C^\RV)$. Proposition \ref{prop:gfl} gives that it is geometrically finite.
Since the holes are contained in the singular set by Lemma \ref{lem:hcs}, $f_n$ converges to $(F, R)$.

To prove the moreover part, if there exists a fixed vertex of degree $d$, then there exists $a\in \RV$ so that $A_{a, n}^{-1} \circ f_n \circ A_{a,n}$ converging algebraically to a degree $d$ rational map. 
Thus $[f_n]$ converges in $\mathcal{M}_{d, \fm}$.

Conversely, if $[f_n]$ converges in $\mathcal{M}_{d, \fm}$, then let $A_n \in \PSL_2 (\C)$ so that $A_n^{-1} \circ f_n \circ A_n$ converges algebraically to a degree $d$ rational map.
Let $v\in \mathcal{V}_s \subseteq \mathcal{V}$ be a periodic point of degree at least $2$. 
Then the pointed Fatou component $A_n^{-1}(U_{s, n}, \psi_n(v))$ converges in Carath\'eodory limit. So $A_n^{-1}\circ A_{v,n}$ converges to a M\"obius transformation. Thus, there exists a fixed vertex of degree $d$.
\end{proof}

\section{Convergence of quasi-invariant forests}\label{sec: cqif}
In this section, we shall prove the convergence of quasi-invariant forests.
Let $\mathcal{H} \subseteq \mathcal{M}_{d, \fm}$ be a hyperbolic component with connected Julia set.
Let $([f_n]) \in \mathcal{H}$ be quasi post-critically finite. 
We have the corresponding quasi-invariant forests
$$
\psi_n: (\mathcal{T}, \mathcal{P}) = \bigcup_{s\in |\mathcal{S}|} (\mathcal{T}_s, \mathcal{P}_s) \longrightarrow (\T_n, \mathcal{P}_{n}),
$$
with vertex set $\mathcal{V} = \bigcup_{s\in |\mathcal{S}|} \mathcal{V}_s$.

Let $F:(\RT, \RV) \longrightarrow (\RT, \RV)$ be the rescaling tree map for $([f_n])$.
For $a\in \RV$, we have the corresponding rescaling $A_{a,n} \in \PSL_2(\C)$ at $a$.
After passing to a subsequence, we assume that the quasi-invariant tree in $a$-coordinate $A_{a, n}^{-1} (\T_{s,n})$ converges in Hausdorff topology for each $a \in \RV$ and $s\in |\mathcal{S}|$.
We denote these limits by
\begin{align}\label{eqn:HL}
\T_{s,a} := \lim_{n\to\infty} A_{a, n}^{-1} (\T_{s,n}) \subseteq \hat \C_a.
\end{align}

We denote the set of Fatou points by $\RV^\mathcal{F} \subseteq \RV$ and Julia points by $\RV^\mathcal{J} \subseteq \RV$.
Note that if $v\in \mathcal{V}$ is a Fatou point for $\mathcal{T}$, then $a=[v]\in \RP$ is a Fatou point for $\RT$.

Let $G\subseteq \hat\C$ be a closed set.
We say it is a {\em graph} with a finite vertex set $V\subset G$ if each component of $G-V$ is an open arc, i.e., homeomorphic to an open interval of $\R$.
A component of $G-V$ is called an {\em edge} of $G$, and an edge $I$ of $G$ is a loop if the two endpoints of $I$ are the same.

The main result in this section is the following.
\begin{theorem}\label{thm:hl}
Let $\mathcal{H} \subseteq \mathcal{M}_{d, \fm}$ be a hyperbolic component with connected Julia set.
Let $([f_n])_n \in \mathcal{H}$ be quasi post-critically finite converging to $(F,R)$ on $(\RT, \hat\C^\RV)$.
For any $s\in |\mathcal{S}|$ and any Fatou point $a\in \RV$,
$\T_{s,a}$ is a graph with finitely many vertices and finitely many non-loop edges.
\end{theorem}

We say $T$ is a {\em $\nu$ star-shaped tree} if $T$ is an union of the $\nu$ arcs $[c, x_i]$ $i=1,..., \nu$, glued at $c$.
We say $T$ is an {\em open $\nu$ star-shaped tree} if $T$ is an union of the $\nu$ arcs $[c, x_i)$ $i=1,..., \nu$, glued at $c$.
With $\U_v$ as defined in Lemma \ref{lem:clr}, we first show

\begin{lem}\label{lem:cfm}
Let $v\in \mathcal{V}_s$ with $a=[v] \in \RV^\mathcal{F}$.
The intersection $\T_{s,a} \cap \U_v$ is an open $\nu$ star-shaped tree where $\nu$ is the valence of $v$ in $\mathcal{T}_s$.

The limit points of $\T_{s,a} \cap \U_v$ on $\partial \U_v$ are pre-periodic points of $R$ and each tangent direction $x \in T_v\mathcal{T}_s$ gives one limit point $\beta_x \in \partial \U_v$.

Moreover, if $x$ is periodic, then $\T_{s,a} \cap \U_v$ gives an attracting direction for $\beta_x$ if $x$ is in the direction of $\p_s \in \mathcal{T}_s$, and a repelling direction otherwise.
\end{lem}

\begin{figure}[ht]
  \centering
  \resizebox{0.4\linewidth}{!}{
    \def\svgwidth{\columnwidth}
\begingroup%
  \makeatletter%
  \providecommand\color[2][]{%
    \errmessage{(Inkscape) Color is used for the text in Inkscape, but the package 'color.sty' is not loaded}%
    \renewcommand\color[2][]{}%
  }%
  \providecommand\transparent[1]{%
    \errmessage{(Inkscape) Transparency is used (non-zero) for the text in Inkscape, but the package 'transparent.sty' is not loaded}%
    \renewcommand\transparent[1]{}%
  }%
  \providecommand\rotatebox[2]{#2}%
  \newcommand*\fsize{\dimexpr\f@size pt\relax}%
  \newcommand*\lineheight[1]{\fontsize{\fsize}{#1\fsize}\selectfont}%
  \ifx\svgwidth\undefined%
    \setlength{\unitlength}{566.92913386bp}%
    \ifx\svgscale\undefined%
      \relax%
    \else%
      \setlength{\unitlength}{\unitlength * \real{\svgscale}}%
    \fi%
  \else%
    \setlength{\unitlength}{\svgwidth}%
  \fi%
  \global\let\svgwidth\undefined%
  \global\let\svgscale\undefined%
  \makeatother%
  \begin{picture}(1,0.9)%
    \lineheight{1}%
    \setlength\tabcolsep{0pt}%
    \put(0,0){\includegraphics[width=\unitlength,page=1]{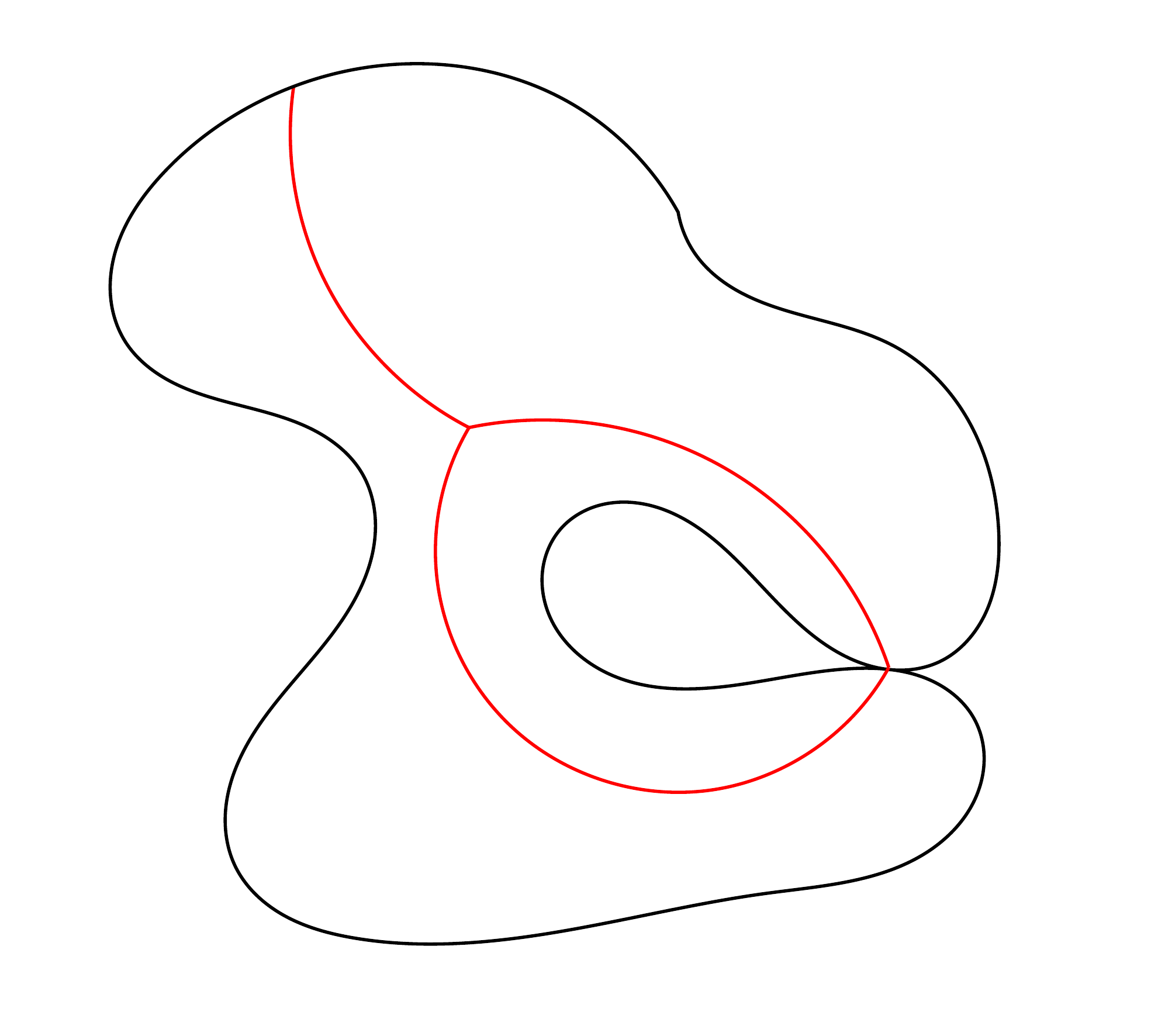}}%
    \put(0.29340603,0.1909864){\color[rgb]{0,0,0}\makebox(0,0)[lt]{\lineheight{1.25}\smash{\begin{tabular}[t]{l}{\Large$\mathcal{U}_{v}$}\end{tabular}}}}%
    \put(0,0){\includegraphics[width=\unitlength,page=2]{CarT.pdf}}%
    \put(0.41141507,0.54734604){\color[rgb]{0,0,0}\makebox(0,0)[lt]{\lineheight{1.25}\smash{\begin{tabular}[t]{l}$v$\end{tabular}}}}%
    \put(0,0){\includegraphics[width=\unitlength,page=3]{CarT.pdf}}%
    \put(0.27929665,0.67282091){\color[rgb]{0,0,0}\makebox(0,0)[lt]{\lineheight{1.25}\smash{\begin{tabular}[t]{l}$\T_{s,a} \cap \U_v$\end{tabular}}}}%
  \end{picture}%
\endgroup%

  }
  \caption{A possible configuration of $\T_{s,a} \cap \U_v$.}
  \label{fig:CarT}
\end{figure}
\begin{proof}
We first note that the quasi-invariant tree $\mathcal{T}_{s,n}$ for $\bp_n \in \BP^\mathcal{S}$ converges to a $\nu$ star-shaped tree in $\D_v$.
For any compact subset $K \subseteq \U_v$, $\T_{s,a} \cap K$ is a star-shaped tree by Lemma \ref{lem:clr}.
Let $x$ be a limit point of $\T_{s,a} \cap \U_v$ in $\partial \U_v$.
Then there exists a sequence $x_n \to x$ with $x_n\in A_{a,n}^{-1}(\T_{s,n} \cap \U_{s,n})$.

We may assume $\U_{s,n} \subseteq \C_a \subseteq \hat \C_a$, and $x\in \C_a$.
Let $m\in\N$ be such that each edge of $\mathcal{T}$ is either fixed or pre-fixed by $F^m$.
Then $x_n$ is either quasi-fixed or quasi-prefixed by $f_n^m$.

If $x_n$ is quasi-fixed, then 
$$
d_{A_{a,n}^{-1}(\U_{s,n})}(x_n, A_{a,n}^{-1}\circ f_n^m \circ A_{a,n}(x_n)) \leq K.
$$ 
Since $x_n \to \partial \U_v$, the hyperbolic metric 
$$
\rho_{\U_{s,n}}(x_n)|dz| \to \infty
$$ 
by Lemma \ref{lem:hme}.
So the Euclidean distance 
$$
d_{\R^2}(x_n, A_{a,n}^{-1}\circ f_n^m \circ A_{a,n}(x_n)) \to 0.
$$ 
Let $x\in D$ be a small disk.
Since the holes of $R^m$ are discrete, we may assume $D$ contains no holes of $R^m$ except possibly at $x$.
Thus, $A_{a,n}^{-1}\circ f_n^m \circ A_{a,n}$ converges compactly to $R^m$ on $D - \{x\}$.
Thus $x$ is fixed by $R^m$.

The same argument shows $x$ is pre-fixed by $R^m$ if $x_n$ is quasi pre-fixed.

Since the holes and these fixed or pre-fixed points are discrete, each tangent direction gives a unique limit point.
The moreover part follows from comparing with the dynamics of $\rl_v^m$ on $\D_v$.
\end{proof}

We remark that different tangent directions of $\mathcal{T}_s$ at $v$ may give the same limit point (see Figure \ref{fig:CarT}). 
Thus the closure $\overline{\T_{s,a} \cap \U_v}$ is not necessarily a tree.
More generally, we have
\begin{lem}\label{lem:htcg}
Let $s\in |\mathcal{S}|$ and $a\in \RV^\mathcal{F} \subseteq \RV$.
Then the Hausdorff limit $\T_{s,a} \subseteq \hat\C_a$ is a connected graph with finitely many vertices.
\end{lem}
\begin{proof}
It suffices to prove the case when $a$ is periodic.
The aperiodic case can be proved by pull back.
Let $m\in\N$ be such that each edge of $\mathcal{T}$ is either fixed or pre-fixed by $F^m$.
Let $R_a^m: \hat\C_a \longrightarrow \hat\C_a$ be the limiting rational map.

Let 
$$
V_1:=\{z\in \hat\C_a: R_a^{2m}(z) = R_a^m(z)\},
$$
and
$$
V_2 := \{\lim_{a} \psi_n(v): v\in \mathcal{V}\}
$$
be the limits of vertices of $\mathcal{T}$ under rescaling at $a$.
Let $V := V_1 \cup V_2$.
Then $\Xi_a \subseteq V$.
Since $a \in \RV^\mathcal{F}$, $R_a^{2m}$ is not identically equal to $R_a^ m$\footnote{This is the only place we use the assumption that $a$ is a Fatou point.}.
Thus, $V$ is finite.
We shall prove that $\T_{s,a}$ is a union of arcs connecting points in $V$.

Let $x\in \T_{s,a} - V$.
Let $x_n \in A_{a,n}^{-1}(\T_{s,n})$ with $\lim_{n\to\infty} x_n = x$.
We assume $x_n$ is quasi-fixed by $f_n^m$, the quasi pre-fixed case can be proved similarly.
Since $x \notin V$, $x$ is neither a hole nor a fixed point of $R_a$.
Thus, Lemma \ref{lem:key} gives that $(A_{a,n}^{-1}(\U_{s,n}), x_n)$ converges in Carath\'eodory topology to $(\U_x, x)$.
Thus, the same proof in Lemma \ref{lem:cfm} gives that the $\T_{s,a} \cap \U_x$ is a union of arcs connecting points in $V$.
\end{proof}

\begin{figure}[ht]
  \centering
  \resizebox{1\linewidth}{!}{
    \def\svgwidth{\columnwidth}
    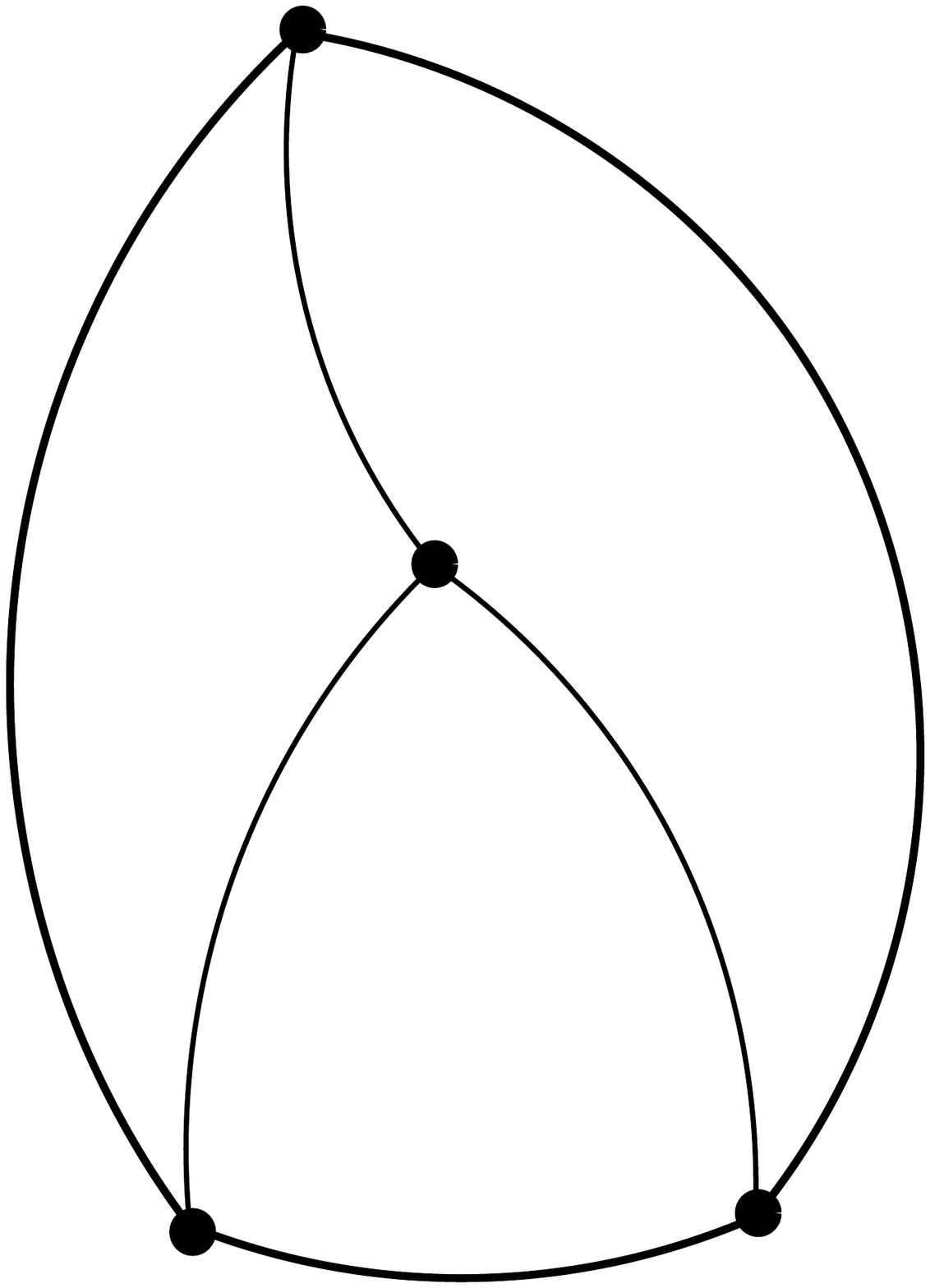

    \def\svgwidth{\columnwidth}
\begingroup%
  \makeatletter%
  \providecommand\color[2][]{%
    \errmessage{(Inkscape) Color is used for the text in Inkscape, but the package 'color.sty' is not loaded}%
    \renewcommand\color[2][]{}%
  }%
  \providecommand\transparent[1]{%
    \errmessage{(Inkscape) Transparency is used (non-zero) for the text in Inkscape, but the package 'transparent.sty' is not loaded}%
    \renewcommand\transparent[1]{}%
  }%
  \providecommand\rotatebox[2]{#2}%
  \newcommand*\fsize{\dimexpr\f@size pt\relax}%
  \newcommand*\lineheight[1]{\fontsize{\fsize}{#1\fsize}\selectfont}%
  \ifx\svgwidth\undefined%
    \setlength{\unitlength}{841.88976378bp}%
    \ifx\svgscale\undefined%
      \relax%
    \else%
      \setlength{\unitlength}{\unitlength * \real{\svgscale}}%
    \fi%
  \else%
    \setlength{\unitlength}{\svgwidth}%
  \fi%
  \global\let\svgwidth\undefined%
  \global\let\svgscale\undefined%
  \makeatother%
  \begin{picture}(1,0.70707071)%
    \lineheight{1}%
    \setlength\tabcolsep{0pt}%
    \put(0,0){\includegraphics[width=\unitlength,page=1]{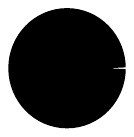}}%
    \put(0.44280744,0.11961272){\color[rgb]{0,0,0}\makebox(0,0)[lt]{\lineheight{1.25}\smash{\begin{tabular}[t]{l}\LARGE{$v$}\end{tabular}}}}%
    \put(0,0){\includegraphics[width=\unitlength,page=2]{Sp.pdf}}%
  \end{picture}%
\endgroup%

  }
  \caption{On the left, we have an illustration of a Carath\'eodory disk. It contains three frontier vertices $v_i$ and one interior vertex $v$. On the right, we have a possible configuration of infinitely many loops.}
  \label{fig:FV}
\end{figure}

The same proof of Lemma \ref{lem:htcg} also gives (see Figure \ref{fig:FV}):
\begin{lem}\label{lem:cedge}
Let $E$ be an edge of $\T_{s,a}$, and $x\in E$.
Then there exists a hyperbolic disk $U_E$, which is the Carath\'eodory limit of pointed disks $A_{a,n}^{-1}((\U_{s,n}, x_n))$ with $x_n \in \T_{s,n}$ and $x_n \to_a x$, such that $E$ is a hyperbolic geodesic or a hyperbolic geodesic ray in $U_E$.
\end{lem}

\begin{defn}\label{defn:cd}
Let $E$ be an edge of $\T_{s,a}$. We call $U_E$ in Lemma \ref{lem:cedge} the corresponding {\em Carath\'eodory disk} for $E$.

A vertex $v$ of $\T_{s,a}$ is called {\em interior} if $v \in U_E$ for some $E$; and is called {\em frontier} otherwise.

We say $\T_{s,a}$ is {\em trivial} if it is a single point, and {\em non-trivial} otherwise.
\end{defn}

Note that the Carath\'eodory disks for different edges are either the same or disjoint.
For future references, we record the following lemma.
\begin{lem}\label{lem:ncfh}
If $\T_{s,a}$ is non-trivial, then there are no periodic critical frontier vertices for $\T_{s,a}$.
\end{lem}
\begin{proof}
Suppose for contradiction that $x$ is a periodic critical frontier vertex.
After passing to an iterate, we may assume $x$ is fixed.
Since $\T_{s,a}$ is non-trivial, there exists an edge $E$ of $\T_{s,a}$ with $x$ as a boundary point.

Let $D$ be a small disk around $x$.
By shrinking $D$ if necessary, we assume $A_{a,n}^{-1} \circ f_n \circ A_{a,n}$ converges compactly on $D-\{x\}$.
We choose $y \in E$ close to $x$, and let $y_n \in \T_{s,n}$ with $y_n \to_a y$.

Since $x$ is a critical fixed point of $R_a$, the local dynamics near $x$ is $z\to z^e$ for some $e\geq 2$.
By Lemma \ref{lem:hme}, for sufficiently large $n$, the hyperbolic metric $d_{A_{a,n}(\U_{s,n})}$ at $z$ is bounded below by $\frac{1}{2|z-x|} |dz|$.
Since $\lim_{u\to 0^+} \int_{u^e}^{u} \frac{1}{t}dt = \infty$, and $A_{a,n}^{-1} \circ f_n \circ A_{a,n}$ converges compactly on $D-\{x\}$, by choosing $y$ sufficiently close to $x$, we can make $d_{\U_{s,n}}(y_n, f_n(y_n))$ arbitrarily large for sufficiently large $n$. This is a contradiction to Theorem \ref{thm:qit}.
\end{proof}

Let $I$ be an edge in $\T_{s,a}$.
Recall $I$ is a {\em loop} if the two endpoints of $I$ are the same.

\begin{lem}\label{lem:fg}
Let $s\in |\mathcal{S}|$ and $a\in \RV^\mathcal{F} \subseteq \RV$.
Then $\T_{s, a}\subseteq \hat\C_a$ has only finitely many non-loop edges\footnote{We do not know if $\T_{s, a}$ can have infinite loops or not (see Figure \ref{fig:FV}).}.
\end{lem}
\begin{proof}
Let $V \subseteq \hat\C_a$ be the set of vertices for $\T_{s,a}$.
Suppose for contradiction that there are infinitely many non-loop edges. 
Then there exist $v, w \in V$ with infinitely many edges connecting them.

Let $C$ be a small circle centered at $v$ so that any edge connecting $v,w$ must intersect $C$.
For each edge, we choose an intersection and denote it by $z_k$.
Passing to a finite iterate and change $C$ by its iterate if necessary, we may assume there are infinitely many of these edges that are quasi-fixed. 
We may also assume there are no fixed points for $R_a$ on $C$ and $C\subseteq \C_a \subseteq \hat \C_a$.
Thus there exists $\epsilon >0$ with $d_{\R^2}(x, R_a(x)) > \epsilon$ for all $x\in C$.

On the other hand, since there are infinitely many edges intersecting $C$, the hyperbolic metrics $d_{\U_{s,n}}$ at $z_k$ can be arbitrarily large by Lemma \ref{lem:hme}.
This is a contradiction to the fact that $f_n$ is quasi post-critically finite.
\end{proof}
\begin{proof}[Proof of Theorem \ref{thm:hl}]
Combining Lemma \ref{lem:htcg} and Lemma \ref{lem:fg}, we have the result.
\end{proof}

\subsection*{Pullback of limiting quasi invariant graph}
Let $\mathcal{T}^{k}_{s}$ be the $k$-th pull back of $\mathcal{T}_{s}$, and let 
$$
\psi_{s,n}: \mathcal{T}^k_{s} \longrightarrow \T^k_{s,n}
$$ 
be the corresponding pull back quasi invariant tree of $\T_{s,n}$ for $f_n$.
The same proof shows that $\T_{s,n}^k$ converges in Hausdorff topology to a graph $\T^k_{s, a}$ with finitely many vertices and finitely many non-loop edges.
Since $\T^k_{s,n} \subseteq \T^{k+1}_{s,n}$, $\T^k_{s,a} \subseteq \T^{k+1}_{s,a}$.
Thus we denote
$$
\T^{\infty}_{s,a} = \bigcup_{k=0}^\infty \T^{k}_{s,a}.
$$

\begin{defn}\label{defn:lqig}
Let $s\in |\mathcal{S}|$ and $a\in \RV^\mathcal{F} \subseteq \RV$.
We call $\T_{s,a}$ the {\em limiting quasi-invariant graph} of $\mathcal{T}_s$ at $a$.
We denote 
\begin{align*}
\T_{a} &:= \bigcup_{s\in |\mathcal{S}|} \T_{s,a} \subseteq\hat\C_a.
\end{align*}
Similarly, we define $\T^k_{a}$ and $\T^\infty_{a}$ respectively.
\end{defn}

\section{Boundedness of Sierpinski carpet rational maps}\label{sec:cscrm}
A closed set $K \subseteq \hat\C$ is called a {\em Sierpinski carpet} if it is obtained from the sphere by removing a countable dense set of open disks, bounded by disjoint Jordan curves whose diameters tend to zero.
We will call a hyperbolic rational map with Sierpinski carpet Julia set a {\em Sierpinski carpet rational map}.
In this section, we will prove quasi post-critically finite degenerations for Sierpinski carpet rational maps are bounded in $\mathcal{M}_{d, \fm}$.

The proof is by contradiction.
We first show the limiting quasi-invariant graphs are disjoint (Proposition \ref{prop:dlqig}).
This is the most technical part of the argument and allows us to construct a curve system $\Sigma_n$ in the complement of critical and post-critical Fatou components for sufficiently large $n$, where curves in $\Sigma_n$ are in correspondence with edges of an invariant subtree $\RT^r\subseteq \RT$ (Lemma \ref{lem:CSC}).

Suppose $[f_n]$ diverges. Then $\RT^r$ is non-trivial (Lemma \ref{lem:ntrt}). 
We show Proposition \ref{prop:mdto} gives the spectral radius $\lambda(\Sigma_n)$ of the Thurston's matrix associated to $\Sigma_n$ is greater or equal to $1$ for all sufficiently large $n$, giving a contradiction (Proposition \ref{prop:mtc}).

To set up notations in this section, let $\mathcal{H} \subseteq \mathcal{M}_{d,\fm}$ be a hyperbolic component with Sierpinski carpet Julia set.
Let $\mathcal{S} = (|\mathcal{S}|, \Phi, \delta)$ be the corresponding mapping scheme.

Let $[f]\in \mathcal{H}$. Then it corresponds to a Blaschke mapping scheme $\bp \in \BP^\mathcal{S}$ by Theorem \ref{thm:BMS}.
Let 
$\U_f = \bigcup_{s\in |\mathcal{S}|} \U_{s, f}$ 
be the union of critical and post-critical Fatou components.
The dynamics $f$ on $\U_f$ is conjugate to $\bp$ on $|\mathcal{S}|\times \D$ that is compatible with the marking.
Let $P_f\subseteq \U_f$ be the post-critical set.

We define a finite forward invariant set of marked points $Q_f \subseteq \U_f$ inductively as follows (see the definition of $\widetilde{\mathcal{P}}$ in \S \ref{sec:DBS}):
\begin{itemize}
\item If $s\in |\mathcal{S}|$ is periodic, $Q_f \cap \U_{s, f}$ is the set of the attracting periodic point in $\U_{s, f}$;
\item Otherwise, $Q_f \cap \U_{s, f} := f^{-1}(Q_f \cap \U_{\Phi(s), f}) \cap \U_{s, f}$.
\end{itemize}
We change the notation from $\widetilde{\mathcal{P}}$ to $Q_f$ to avoid potential confusions with the post-critical set $P_f$.
We denote $V_f$ as the set of critical values for $f$.

\subsection{Lifting arcs and simple closed curves}
Let $X$ be a topological space. We define an {\em arc} $I\subseteq X$ as a continuous map 
$$
h: [0,1] \subseteq \R \longrightarrow I = h([a,b]) \subseteq X.
$$
We shall not differentiate the continuous map with its image $I\subseteq X$.
Two arcs $h_0, h_1$ are said to be {\em homotopic relative to its boundary} if
there exists a continuous map 
$$
H: [0,1]\times [0,1] \longrightarrow X
$$ 
with $H(t, 0) = h_0(t)$, $H(t,1) = h_1(t)$, $H(0, s) = h_0(0) = h_1(0)$ and $H(1,s) = h_0(1) = h_1(1)$.
We shall denote it by $I_0 \sim_\partial I_1$.
An arc $I$ is said to connect $x, y\in X$ if $h(0) = x$ and $h(1) = y$.
It is said to be simple if $h$ is an embedding.
If $f: Y\longrightarrow X$ is a covering, we say $I'$ is a {\em lift} of $I$ if $I'$ is a component of $f^{-1}(I)$.

Similarly, a closed curve in $X$ is defined as a continuous map 
$$
h: \mathbb{S}^1 \longrightarrow \gamma = h(\mathbb{S}^1) \subseteq X.
$$
We do not differentiate the continuous map with its image $\gamma\subseteq X$, and two closed curves are said to be {\em homotopic} if the two continuous maps are homotopic.
We denote this by $\gamma_0 \sim \gamma_1$.
It is said to be {\em simple} if $h$ is an embedding.
If $f: Y\longrightarrow X$ is a covering, we say $\gamma'$ is a {\em lift} of $\gamma$ if $\gamma'$ is a component of $f^{-1}(\gamma)$ and the {\em degree} of $\gamma'$ is defined as the degree of the covering $f: \gamma' \longrightarrow \gamma$.

\begin{lem}\label{lem:isc}
Let $f$ be a hyperbolic Sierpinski carpet rational map.
Let $I \subseteq \hat\C - (V_f\cup Q_f)$ be an arc connecting $x, y \in \hat\C - \U_f$, such that $I\cap \U_{s, f}$ does not separate $V_f\cup Q_f$ in $\U_{s, f}$ for any $s\in |\mathcal{S}|$.
Let $I'$ be a lift of $I$.
Then there exists an arc $J\subseteq \hat\C - \U_f$ so that
\begin{itemize}
\item $J \sim_\partial I$ in $\hat\C - (V_f\cup Q_f)$, and
\item $J$ has a lift $J'\sim_\partial I'$ in $\hat\C - Q_f$.
\end{itemize}

Similarly, let $\alpha \subseteq \hat\C - (V_f\cup Q_f)$ be a simple closed curve such that $\alpha \cap \U_{s, f}$ does not separate $V_f\cup Q_f$ in $\U_{s, f}$ for any $s\in |\mathcal{S}|$.
Let $\alpha'$ be a lift of degree $e$. Then there exists a simple closed curve $\beta\subseteq \hat\C - \U_f$ so that
\begin{itemize}
\item $\beta \sim \gamma$ in $\hat\C - (V_f\cup Q_f)$, and
\item $\beta$ has a lift $\beta'$ of degree $e$ homotopic to $\alpha'$ in $\hat\C - Q_f$.
\end{itemize}
\end{lem}
\begin{figure}[ht]
  \centering
  \resizebox{0.6\linewidth}{!}{
    \def\svgwidth{\columnwidth}
\begingroup%
  \makeatletter%
  \providecommand\color[2][]{%
    \errmessage{(Inkscape) Color is used for the text in Inkscape, but the package 'color.sty' is not loaded}%
    \renewcommand\color[2][]{}%
  }%
  \providecommand\transparent[1]{%
    \errmessage{(Inkscape) Transparency is used (non-zero) for the text in Inkscape, but the package 'transparent.sty' is not loaded}%
    \renewcommand\transparent[1]{}%
  }%
  \providecommand\rotatebox[2]{#2}%
  \newcommand*\fsize{\dimexpr\f@size pt\relax}%
  \newcommand*\lineheight[1]{\fontsize{\fsize}{#1\fsize}\selectfont}%
  \ifx\svgwidth\undefined%
    \setlength{\unitlength}{1020.47244094bp}%
    \ifx\svgscale\undefined%
      \relax%
    \else%
      \setlength{\unitlength}{\unitlength * \real{\svgscale}}%
    \fi%
  \else%
    \setlength{\unitlength}{\svgwidth}%
  \fi%
  \global\let\svgwidth\undefined%
  \global\let\svgscale\undefined%
  \makeatother%
  \begin{picture}(1,0.5)%
    \lineheight{1}%
    \setlength\tabcolsep{0pt}%
    \put(0,0){\includegraphics[width=\unitlength,page=1]{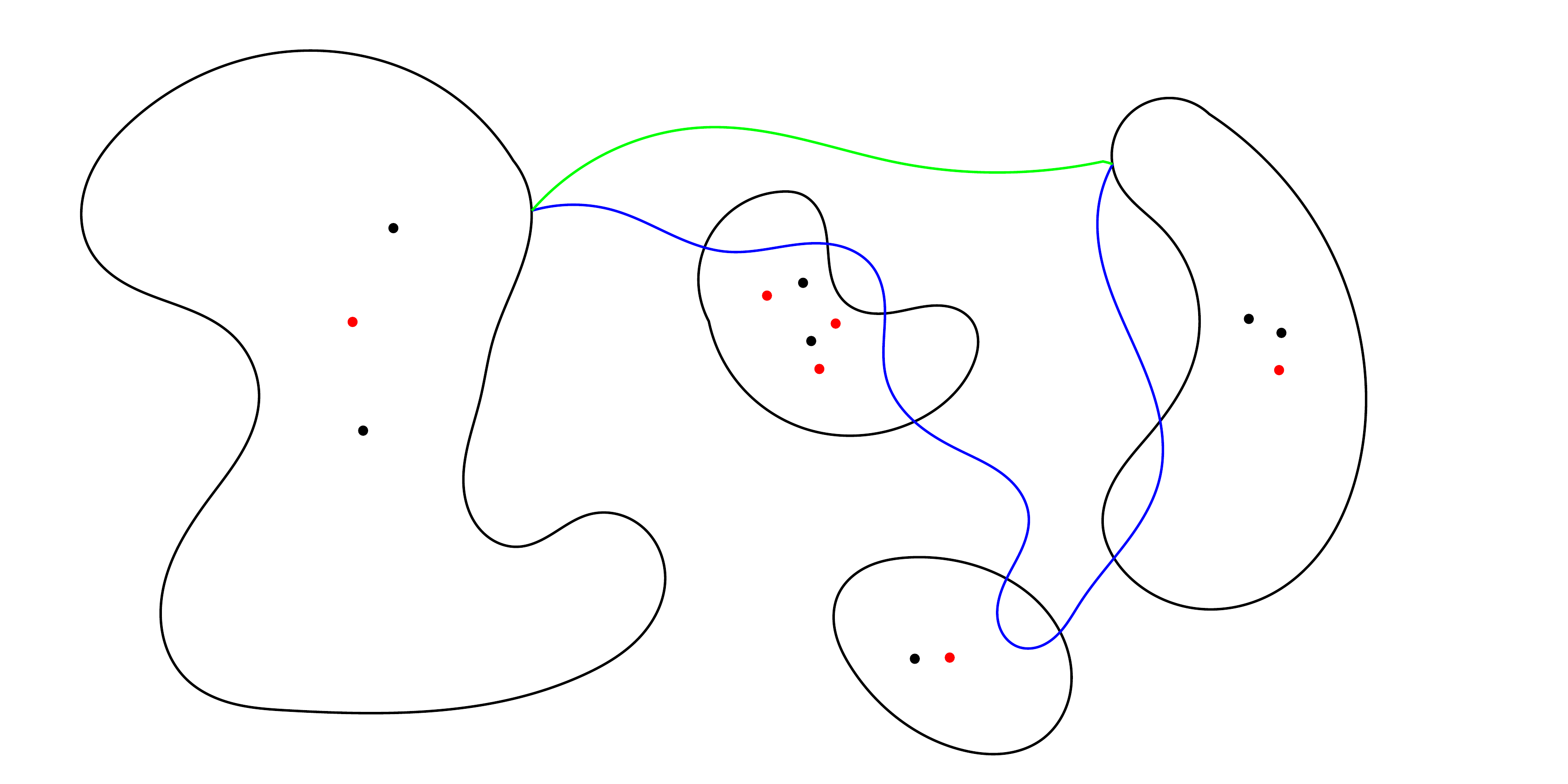}}%
    \put(0.49296375,0.42427364){\color[rgb]{0,0,0}\makebox(0,0)[lt]{\lineheight{1.25}\smash{\begin{tabular}[t]{l}$J$\end{tabular}}}}%
    \put(0.61843194,0.22159426){\color[rgb]{0,0,0}\makebox(0,0)[lt]{\lineheight{1.25}\smash{\begin{tabular}[t]{l}$I$\end{tabular}}}}%
    \put(0,0){\includegraphics[width=\unitlength,page=2]{HomLem.pdf}}%
  \end{picture}%
\endgroup%

  }
  \caption{An illustration of Lemma \ref{lem:isc}. The black and red points represent elements in $V_f$ and $Q_f$ respectively. The blue and green arcs represent $I$ and $J$ in Lemma \ref{lem:isc}.}
  \label{fig:HomLem}
\end{figure}
\begin{proof}
Let $L$ be a component of $I \cap \U_{s, f}$.
Since $L$ does not separate $V_f\cup Q_f$ in $\U_{s, f}$, 
we can homotope $L$ to an arc in $\partial \U_{s, f}$ in $\hat\C - (V_f\cup Q_f)$ relative to its boundary.
Let $J$ be the arc after performing all such homotopy modifications, then $J \sim_\partial I$ in $\hat\C - (V_f\cup Q_f)$, and $J \subseteq \hat\C - \U_f$.
Since $J\sim_\partial I$ in $\hat\C - V_f$, we can lift the homotopy, and get a lift $J'$ of $J$ connecting $\partial I'$. 
Since $Q_f$ is forward invariant, $Q_f \subseteq f^{-1}(Q_f)$. Thus, $J' \sim_\partial I'$ in $\hat\C - Q_f$.

The same proof also works for simple closed curves $\alpha$.
\end{proof}

\begin{lem}\label{lem:LT}
Let $f$ be a hyperbolic Sierpinski carpet rational map.
Let $x\neq y$ be repelling fixed points.
Let $I \subseteq \hat\C - (V_f\cup Q_f)$ be an arc connecting $x, y$ such that $I\cap \U_{s, f}$ does not separate $V_f\cup Q_f$ in $\U_{s, f}$ for any $s\in |\mathcal{S}|$.
Let $I'$ be a lift of $I$. Then $I'\not\sim_\partial I$ in $\hat\C - Q_f$.

Similarly, if $\alpha \subseteq \hat\C - (V_f\cup Q_f)$ is a homotopically non-trivial simple closed curve such that $\alpha \cap \U_{s, f}$ does not separate $V_f\cup Q_f$ in $\U_{s, f}$ for any $s\in |\mathcal{S}|$.
Suppose $\alpha'$ is a lift of $\alpha$ by $f$ of degree $1$. Then $\alpha' \not\sim \alpha$ in $\hat\C - Q_f$.
\end{lem}
\begin{proof}
Suppose for contradiction that $I' \sim_\partial I$.
Applying Lemma \ref{lem:isc}, we may assume that $I \subseteq \hat\C - \U_f$.
Then $I' \sim_\partial I$ in $\hat\C - Q_f$ implies $I' \sim_\partial I$ in $\hat\C - \U_f$.
Since the post-critical set $P_f \subseteq \U_f$, $I' \sim_\partial I$ in $\hat\C - P_f$.
Therefore, we can inductively take preimage $I^{(k)}$ so that 
$$
f^{k}(I^{(k)}) = I \text{ and }I^{(k)} \sim_\partial I \text{ in } \hat\C - P_f.
$$
Since $f$ is hyperbolic, $f$ is uniformly expanding on the compact set $\hat \C - \U_f$ with respect to the hyperbolic metric on $\hat \C - P_f$. 
Therefore, the spherical diameters of $I^{(k)}$ shrink to $0$, contradicting to $x \neq y$.

The case of simple closed curves can be proved similarly.
\end{proof}

\subsection{Disjoint limiting quasi-invariant graphs}
Let $([f_n])_n \in \mathcal{H}$ be quasi post-critically finite, converging to $(F, R)$ on $(\RT, \hat\C^\RV)$.
By replacing $\mathcal{T}$ with some pull back $\mathcal{T}^m$ if necessary, we assume that the following technical assumption holds throughout this section:
\begin{assumption}\label{assumption:2}
There exists a Fatou point $v\in \mathcal{V}_s$ for each $s\in |\mathcal{S}|$.
\end{assumption}

Recall that $\T_{s, a}$ for $a\in \RV$ and $s \in |\mathcal{S}|$ is defined in Equation \ref{eqn:HL}.
The main goal of this subsection is to show:
\begin{prop}\label{prop:dlqig}
The limiting quasi-invariant graphs are disjoint: For each $s\in |\mathcal{S}|$, there exists a unique $a = a(s) \in \RV^\mathcal{F}$ so that $\T_{s, a} \subseteq \hat \C_a$ is non-trivial.
Moreover, for all $t \neq s\in |\mathcal{S}|$, $\T_{s, a} \cap \T_{t, a} = \emptyset$.
\end{prop}

\begin{figure}[ht]
  \centering
  \resizebox{0.6\linewidth}{!}{
    \def\svgwidth{\columnwidth}
    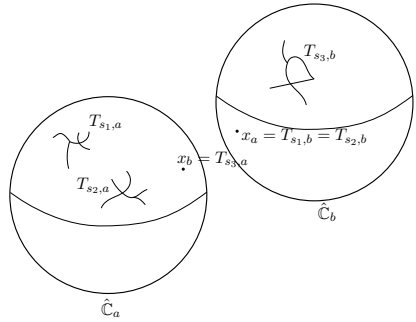

  }
  \caption{An illustration of limiting quasi-invariant graphs. The point $x_b \in \Xi_a$ (or $x_a \in \Xi_b$) is the point in the singular set corresponds to the direction of $b$ (or $a$).}
  \label{fig:DisjointG}
\end{figure}

We start with the following lemma.
\begin{lem}\label{lem:mK}
There exists $m\in \N$ so that each edge $E$ of $\mathcal{T}$ is either fixed or quasi-fixed by $F^m$, i.e., $F^m(E)$ is fixed by $F^m$.

Let $V_n$ be the set of critical values of $f_n^m$, and $Q_n := Q_{f_n}$.
Then there exists $K > 0$ so that for all $n$ and for all $x_n\in V_n \cup Q_n$,
$$
d_{\U_n}(x_n, \mathcal{V}_n) \leq K,
$$
where $\mathcal{V}_n$ is the vertex set of $\T_n$.
\end{lem}
\begin{proof}
The first part follows immediately from the fact that $F : \mathcal{T} \longrightarrow \mathcal{T}$ is simplicial and that $\mathcal{T}$ is a finite forest.

By the construction, $Q_n$ is within a uniform bounded hyperbolic distance from the marked set $\mathcal{P}_n$ of $\T_n$.
If $x_n$ is a critical value of $f_n^m$, then $x_n = f_n^k(c_n)$ for some critical point of $f_n$ and $k\leq m$.
So the second part follows from the fact $[f_n]$ is quasi post-critically finite.
\end{proof}

\subsection*{Periodic frontier vertices}
We now analyze a periodic frontier vertex $v$\footnote{See Definition \ref{defn:cd}.} of a non-trivial limiting quasi-invariant graph $\T_{s,a}$, with $s\in |\mathcal{S}|$ and $a\in \RV^\mathcal{F}$.
The strategy is to consider a small circle $C \subseteq \hat\C_a$ around $v$.
The incident edges of the limiting quasi-invariant graph $\T_a = \bigcup_{s\in |\mathcal{S}|}\T_{s,a}$ cut the circle into finitely many arcs.
We study the lifts of these arcs.

Throughout this discussion, we fix the constants $m$ and $K$ as in Lemma \ref{lem:mK}.
Note the period of $v$ divides $m$.
We assume the quasi-invariant constant for $f_n^m$ on $\T_n$ as in Theorem \ref{thm:qit} is also bounded by $K$.
By Lemma \ref{lem:ncfh}, $v$ is not a critical point of $R_a^m$.
To avoid complicated indices, we denote 
$$
G: = R_a^m: \hat \C_a \longrightarrow \hat\C_a.
$$

Let $U$ be a small disk containing $v$ so that the only possible hole of $G$ in $U$ is $v$. 
We also assume $U$ contains no other vertices of $\T_a$ other than $v$.
Let $C, C' \subseteq U$ be a small simple closed curve around $v$ so that $G (C') = C$.
Denote $D$ and $D'$ as the disk bounded by $C$ and $C'$ that contains $v$.
Since $v$ has local degree $1$ by Lemma \ref{lem:ncfh}, $G: D' \longrightarrow D$ is a degree $1$ map.

\begin{figure}[ht]
  \centering
  \resizebox{0.8\linewidth}{!}{
    \def\svgwidth{\columnwidth}
    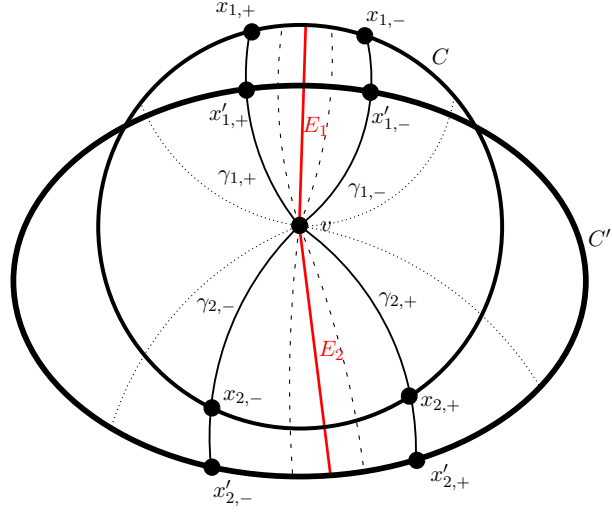

  }
  \caption{A schematic picture of the neighborhood of $v$. The dashed lines correspond to $\partial X_{E_i}$ and the doted lines correspond to $\partial U_{E_i}$.}
  \label{fig:NHV}
\end{figure}

Shrinking $U$ if necessary, we assume that $C$ only intersects edges of $\T_a$ incident at $v$.
Let $E$ be an edge of $\T_a$ intersected by $C$. 
Let $\hat E$ be the corresponding edge of $\mathcal{T}$.
\begin{lem}\label{lem:IEP}
$\hat E$ is fixed by $F^m$.
\end{lem}
\begin{proof}
Let $U_E$ be the Carath\'eodory disk for $E$ then $G(U_{E})$ is invariant under $G = R_a^m$.
If $\hat E$ is not fixed by $F^m$, then $G^2(U_E) = G(U_E) \neq U_E$.
Since $v$ is fixed by $G$ and is on the boundary of both $U_{E}$ and $G(U_{E})$, contradicting to the degree at $v$ is $1$.
\end{proof}

Let $E$ be an edge of $\T_{a}$, and let $U_E$ be the corresponding Carath\'eodory disk.
Define 
$$
X_E:=\{z\in U_E: d_{U_E}(z, \T_{a} \cap U_E) \leq 2K\}.
$$
Note that $\partial X_E$ are piecewise smooth curves.
Similar proof as in Lemma \ref{lem:fg} shows that there are finitely many edges $E$ of $\T_{a}$ with $C \cap X_E \neq \emptyset$.
Modify $C$ if necessary, we assume $C$ intersects $X_E$ if and only if $C$ intersects $E$.
To avoid cumbersome notations, if $C$ intersects a loop, we introduce a vertex in the middle of the edge, and treat it as two different non-loop edges.

If $E$ intersects $C$, we choose two invariant curves $\gamma_{E, \pm} \subseteq U_E\cap \overline{(D \cup D')}$ so that (see Figure \ref{fig:NHV})
\begin{itemize}
\item $\gamma_{E, \pm}$ connects $v$ to $\partial (D \cup D')$, with $\Int(\gamma_{E, \pm}) \subseteq D \cup D'$;
\item $\gamma_{E, \pm}$ bounds $X_E$ in $D \cup D'$.
\item If $x\in \gamma_{E, \pm} \cap D'$, then $G(x)\in \gamma_{E, \pm}$;
\item If $y\in \gamma_{E, \pm} \cap D$, then there exists $x \in \gamma_{E, \pm} \cap D'$ with $y = G(x)$.
\end{itemize}
We assume $\gamma_{E, -}, X_E, \gamma_{E, +}$ are ordered counterclockwise at $v$.

Let $E_1,..., E_k$ be the edges intersecting $C$, where the indices are labeled counterclockwise at $v$, and assume that $E_i$ is an edge of $\T_{s_i, a}$, with $s_i \in |\mathcal{S}|$.
Let $\gamma_{i, \pm}$ be the corresponding invariant curve for $E_i$.
Let $x_{i, \pm}$ and $x_{i, \pm}'$ be the intersection points of $\gamma_{i, \pm}$ with $C$ and $C'$ respectively.
Thus, we can decompose
\begin{align*}
C &= H_1 \cup I_1 \cup.... \cup H_k \cup I_k,
\end{align*}
where $H_i$ connects $x_{i, -}$ with $x_{i, +}$ and $I_i$ connects $x_{i, +}$ and $x_{i+1, -}$ counterclockwise (see Figure \ref{fig:NHV}).
Similarly, we have 
$$
C'= H_1' \cup I_1' \cup.... \cup H_k' \cup I_k'.
$$
Modify $C$ if necessary, we assume 
\begin{itemize}
\item $H_i \subseteq U_{E_i}$;
\item $I_i \cap X_E = \emptyset$ for all edges $E$ of $\T_a$.
\end{itemize}

Let $t_{i, \pm} \in \R/\Z$ be the angle landing on the corresponding side of $\hat E_i$ determined by $x_{i, \pm}$.
By Lemma \ref{lem:lacs}, $\eta_{s_i,n}(t_{i, \pm}) \in s_i\times \mathbb{S}^1$ is a fixed point of $\bp_n^m$, which corresponds to a fixed points $q_{i, \pm, n} \in \partial \U_{s_i,n}$ of $f_n^m$.

Let $x_{i, \pm, n} = A_{a,n}(x_{i, \pm})$. 
Then $x_{i, \pm, n} \in \U_{s_i,n}$ for all sufficiently large $n$.
Let $I_{i, n} = A_{a,n}(I_i)$ connecting $x_{i, +, n} $ and $x_{i+1, -, n} $.
Since $A_{a,n}^{-1} \circ f_n^m \circ A_{a,n}$ converges compactly to $G$ on $U - \{v\}$ and $I_i' \subseteq U - \{v\}$, $I_{i,n}$ has a lift $I_{i,n}'$ by $f_n$ connecting $x_{i, +, n}'$ and $x_{i+1, -, n}'$, with $x_{i, \pm, n}'\to_a x_{i,\pm}'$.

\begin{lem}\label{lem:carc}
For sufficiently large $n$, there exists an arc $L_{i,\pm, n}$ connecting $x_{i,\pm, n}$ and $q_{i, \pm, n}$ so that
\begin{itemize}
\item $\Int(L_{i, \pm, n}) \subseteq \U_{s_i,n}$;
\item $L_{i,\pm, n}$ has a lift $L_{i, \pm, n}'$ connecting $x_{i,\pm,n}'$ and $q_{i,\pm,n}$;
\item $d_{\U_{s_i,n}}(L_{i,\pm,n}, \T_{s_i,n}) \geq 2K$.
\end{itemize} 
\end{lem}
\begin{proof}
Since $f_n^m$ on $\U_{s_i,n}$ is conformally conjugate to $\bp_n^m$ on $s_i\times \D$, it suffices to construct such an arc in $s_i\times \D$.
Abusing the notations, we will denote the quasi-invariant trees in $s\times \D$ by $\mathcal{T}_s$ as well.
We use the notational convention that the elements in $s\times \D$ are labeled by a hat.

Let $\rl_i: \D \longrightarrow \D$ be the corresponding rescaling limit for $G: U_{E_i} \longrightarrow U_{E_i}$.
More precisely, we have a rescaling $M_{E_i, n} \in \Isom(\Hyp^2)$ with $\hat M_{E_i, n}(z) = (s_i, M_{E_i, n}(z)) \in s_i \times \D$ so that
$$
\hat M_{E_i, n}^{-1} \circ \bp_n^m \circ \hat M_{E_i, n}
$$
converges compactly to $\rl_i$.

Let $\hat v_i \in \mathbb{S}^1$ and $\hat x_{i, \pm} \in \D$ be the corresponding point for $v$ and $x_{i, \pm}$, and $\hat \gamma_{i, \pm}$ be the corresponding invariant curve.
Let $\hat u_{i, \pm} \in \mathbb{S}^1$ close to $\hat v_i$ from above and below respectively so that the hyperbolic geodesic $\hat{\sigma}_{i, \pm} \subseteq \D$ connecting $\hat x_{i, \pm}$ and $\hat u_{i, \pm}$ is at least $2K$ distance to $ \mathcal{T}_{s_i,E_i}$,
where $\mathcal{T}_{s_i,E_i}$ is the limit $\mathcal{T}_{s_i,n}$ under normalization $\hat M_{E_i, n}$.

Let $[\hat u_{i, -}, \hat v_i], [\hat v_i, \hat u_{i, +}] \subseteq \mathbb{S}^1$.
Let $[\hat u_{i, -}', \hat v_i], [\hat v_i, \hat u_{i,+}'] \subseteq \mathbb{S}^1$ be the corresponding lift by $\rl_i$ respectively.
By choosing $\hat u_{i, \pm}$ closer to $\hat v_i$, we assume that both $[\hat u_{i, -}', \hat v_i)$ and $(\hat v_i, \hat u_{i,+}']$ contain no holes of $\rl_i$.

Let $\hat \beta_{i, +} = [\hat v_i, \hat u_{i, +}] \cup \hat{\sigma}_{i, +}$ and $\hat \beta_{i, -} = [\hat u_{i,-}, \hat v_i] \cup \hat{\sigma}_{i, -}$.
Note $\hat \beta_{i, \pm} \sim_\partial \hat \gamma_{i, \pm}$ in $\hat \C - V(\rl_i)$ where $V(\rl_i)$ is the set of critical values of $\rl_i$ viewed as a rational map.
So $\hat \beta_{i, \pm}$ has a lift $\hat \beta_{i, \pm}'$ connecting $\hat v_i$ and $\hat x_{i, \pm}'$.
Thus, there exists a lift $\hat{\sigma}_{i, \pm}'$ of $\hat{\sigma}_{i, \pm}$ connecting $\hat x_i'$ and $\hat u_i'$ (see Figure \ref{fig:LID}).

\begin{figure}[ht]
  \centering
  \resizebox{0.8\linewidth}{!}{
    \def\svgwidth{\columnwidth}
    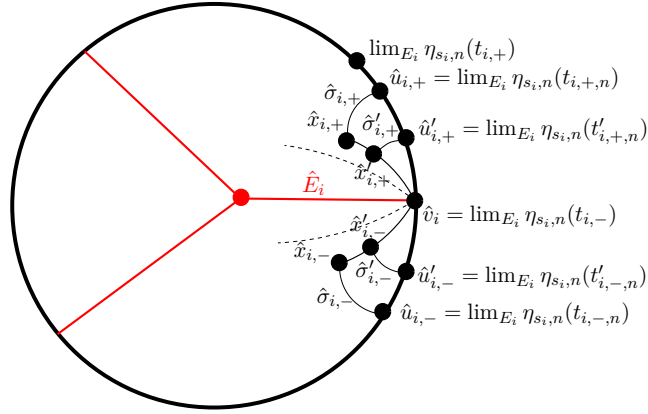

  }
  \caption{A schematic picture of the curves and points defined for $\rl_i: \D\rightarrow\D$. Dashed lines correspond to $\partial \hat{X}_{E_i}$.}
  \label{fig:LID}
\end{figure}

Let $t_{i,\pm,n}, t_{i,\pm,n}' \in \R/\Z$ with $t_{i, \pm, n} = m_{e_i}(t_{i,\pm,n}')$, $\eta_{s_i,n}(t_{i,\pm,n}) \to_{E_i} \hat u_{i, \pm}$ and $\eta_{s_i,n}(t_{i,n}') \to_{E_i} \hat u_{i, \pm}'$, where $e_i = \deg \mathcal{F}_n^m$ on $s_i\times \D$.
Since the angle landing at the corresponding edge $E_i$ is $t_{i, \pm}$,
by making $\hat u_{i,\pm}$ closer to $\hat v_i$ if necessary, we can assume both $t_{i,\pm, n}$ and $t_{i,\pm, n}'$ are within $\frac{1}{e_i}$ distance (in the Euclidean metric on $\R/\Z$) to $t_{i, \pm}$ for sufficiently large $n$.
Thus the shorter arc $\tau_{i, \pm,n}'$ between $t_{i,\pm}$ and $t_{i, \pm, n}'$ is a lift of the shorter arc $\tau_{i, \pm,n}$ of $t_{i,\pm}$ and $t_{i, \pm, n}$ by $m_{e_i}$ for sufficiently large $n$.\footnote{We remark that the shorter arc can be either $[t_{i,\pm}, t_{i, \pm, n}]$ or $[t_{i,\pm,n}, t_{i, \pm}]$ depending on whether $\hat v_i$ is an attracting, a singly parabolic, a doubly parabolic or a repelling fixed point of $\rl_i$ (see Figure \ref{fig:LID}).}
Thus, $\eta_{s_i,n}(\tau_{i,\pm,n})$ has a lift connecting $\eta_{s_i,n}(t_{i,\pm,n}')$ and $\eta_{s_i,n}(t_{i,\pm})$, where $\eta_{s_i,n}: \R/\Z \longrightarrow \mathbb{S}^1$ is the marking for $s_i \times \D$.

Let $\hat{\sigma}_{i, \pm, n}$ be the geodesic connecting $\hat x_{i,\pm,n}$ to $\eta_{s_i,n}(t_{i,\pm,n})$.
Then $\hat M_{E_i, n}^{-1}(\hat{\sigma}_{i,\pm,n})$ converges to $\hat{\sigma}_{i,\pm}$.
Since $\hat{\sigma}_{i,\pm}'$ contains no holes of $\rl_i$, $\hat M_{E_i, n}^{-1} \circ \bp_n^m \circ \hat M_{E_i, n}$ converges uniformly near $\hat{\sigma}_{i,\pm}'$.
Therefore, $\hat{\sigma}_{i, \pm, n}$ has a lift $\hat{\sigma}_{i, \pm, n}'$ connecting $x_{i,\pm,n}'$ to $\eta_{s_i,n}(t_{i,\pm,n}')$ for all large $n$.

Let $\hat L_{i,\pm,n}$ be constructed from $\hat{\sigma}_{i,\pm,n}\cup \eta_{s_i,n}(\tau_{i,\pm,n})$ by pushing the arc on $\mathbb{S}^1$ slightly towards $\D$.
Then for sufficiently large $n$, $\hat L_{i,\pm,n}$ satisfies the three conditions.
\end{proof}

Consider the arc $J_{i, n} = L_{i,+,n} \cup I_{i,n} \cup L_{i+1,-,n}$ connecting $q_{i,+,n}$ and $q_{i+1,-,n}$.
Then by Lemma \ref{lem:carc}, for sufficiently large $n$, $J_{i, n}$ has a lift $J_{i, n}'$ again connecting $q_{i,+,n}$ and $q_{i+1,-,n}$.
\begin{lem}\label{lem:JNS}
For sufficiently large $n$, $J_{i, n}$ does not separate $V_n \cup Q_n$ in $\U_{s,n}$ for all $s\in |\mathcal{S}|$.
\end{lem}
\begin{proof}
For sufficiently large $n$, the intersection of $J_{i, n}$ with $\U_{s, n}$ is at least $K$ distance away from $\T_{s,n}$ for all $s\in |\mathcal{S}|$.
On the other hand, any point in $V_n\cup Q_n$ is within $K$ distance away from $\T_{s,n}$ by Lemma \ref{lem:mK}.
Thus, $J_{i, n}$ does not separate $V_n \cup Q_n$ in $\U_{s,n}$ for all $s\in |\mathcal{S}|$.
\end{proof}

\begin{lem}\label{lem:JTE}
For sufficiently large $n$, $J_{i, n}' \sim_\partial J_{i, n}$ in $\hat \C - Q_n$.
\end{lem}
\begin{proof}
We shall construct the homotopy piecewisely.

Let $\chi_{i, \pm} \subseteq \gamma_{i, \pm}$ be the arc connecting $x_{i, \pm}$ and $x_{i, \pm}'$ (see Figure \ref{fig:NHV}).
By the local dynamics of $G$ near $v$, we can perform a homotopy 
$$
H_i: [0,1] \times [0,1] \longrightarrow U-\{v\}
$$
between $I_i$ and $I_i'$ in $U- \{v\}$ 
with $H(0, s)$ and $H(1,s)$ parameterizing $\chi_{i,+}$ and $\chi_{i+1,-}$ respectively.
Since $U$ contains no vertices nor holes other than $v$, for sufficiently large $n$, $I_{i,n}$ is homotopic to $I_{i,n}'$ in $\hat \C - Q_n$ where the boundary points stay on $A_{a,n}(\chi_{i,+})$ and $A_{a,n}(\chi_{i+1,-})$ throughout the homotopy.

Note that $Q_n \cap \U_{s_i,n} = \{p_{i,n}\}$ where $p_{i,n}$ is the attracting periodic point in $\U_{s_i,n}$.
For sufficiently large $n$, we have a homotopy $H_{i,\pm,n}$ between $L_{i, \pm, n}$ to $L_{i, \pm, n}'$ in $\overline{\U_{s_i, n}} - \{p_{i,n}\}$,
where $H_{i,\pm,n}(0, s) =q_{i, \pm, n}$ and $H_{i,\pm,n}(1,s)$ parameterizes $A_{a,n}(\chi_{i, \pm})$.

By gluing the homotopies together, we conclude the result.
\end{proof}

\begin{cor}\label{cor:qeq}
For sufficiently large $n$, $q_{i,+,n} = q_{i+1,-,n}$ and $J_{i, n}$ is homotopically trivial in $\hat \C - (V_n\cup Q_n)$.
\end{cor}
\begin{proof}
If $q_{i,+,n} \neq q_{i+1,-,n}$, by Lemma \ref{lem:JNS} and Lemma \ref{lem:LT}, $J_{i,n} \not\sim_\partial J_{i,n}'$ in $\hat \C - Q_n$ contradicting to Lemma \ref{lem:JTE}.
Therefore $q_{i,+,n} = q_{i+1,-,n}$.
Since $J_{i,n}'$ is a degree $1$ lift of $J_{i,n}$, if $J_{i, n}$ were not homotopically trivial, Lemma \ref{lem:LT} would give that $J_{i,n} \not\sim J_{i,n}'$ in $\hat \C - Q_n$ contradicting to Lemma \ref{lem:JTE}.
\end{proof}

\begin{cor}\label{cor:HCU}
For sufficiently large $n$, $C_n = A_{a,n}(C)$ is homotopic to a simple closed curve $\tilde{C}_n \subseteq \U_{s,n}$ in $\hat\C - (V_n \cup Q_n)$.
\end{cor}
\begin{proof}
Since $J_{i, n}$ is homotopically trivial in $\hat\C - (V_n \cup Q_n)$, we can homotope $I_{i,n}$ into an arc in $\U_{s,n}$ in $\hat\C - (V_n \cup Q_n)$ while the boundary points stay in $\U_{s,n}$.
Therefore $C_n$ is homotopic to a simple closed curve $\tilde{C}_n \subseteq \U_{s,n}$ in $\hat\C - (V_n \cup Q_n)$.
\end{proof}

\begin{cor}\label{cor:NSC}
For sufficiently large $n$, $C_n = A_{a,n}(C)$ does not separate critical points of $f_n^m$.
\end{cor}
\begin{proof}
Note $C_n$ has a degree $1$ lift $C_n'$ with $C_n' \sim C_n$ in $\hat\C - (V_n \cup Q_n)$.
Let $\tilde{C}_n \subseteq \U_{s,n}$ be as in Corollary \ref{cor:HCU}.
Then it has a degree $1$ lift $\tilde{C}_n' \subseteq \U_{s,n}$ with $\tilde C_n' \sim C_n'$ in $\hat\C - f_n^{-m}(V_n \cup Q_n)$.
Let $\tilde D_n$ and $\tilde D_n'$ in $\U_{s,n}$ be the disks bounded by $\tilde{C}_n$ and $\tilde{C}_n'$.
Since $f_n^m: \U_{s,n} \longrightarrow \U_{s,n}$ is proper, $\tilde D_n'$ is a component of $f^{-m}_n(\tilde D_n)$.
Since $\tilde{C}_n'$ is a degree $1$ lift of $\tilde C_n$, $f_n^m: \tilde D_n' \longrightarrow \tilde D_n$ has degree $1$.
Thus $\tilde{D}_n'$ contains no critical points of $f_n^m$.
Hence, $\tilde C_n'$ and $C_n'$ does not separate critical points.
\end{proof}

We are now ready to prove Proposition \ref{prop:dlqig}:
\begin{proof}[Proof of Proposition \ref{prop:dlqig}]
Let $s\in |\mathcal{S}|$. 
Choose $a\in \RV^\mathcal{F}$ so that $\T_{s,a}$ is non-trivial.
Suppose for contradiction that $t \neq s \in |\mathcal{S}|$ with $\T_{s,a} \cap \T_{t,a} \neq \emptyset$.
Let $m\in \N$ be the constant as in Lemma \ref{lem:mK}.

Assume first they intersect at a periodic frontier vertex $v$.

If $\T_{t,a}$ is also non-trivial, then we can find a small circle $C \subseteq \hat\C_a$ around $v$ intersecting both edges in $\T_{s,a}$ and $\T_{t,a}$.
Thus, we assume $C$ intersects $E_s, E_t$ of $\T_{s,a},\T_{t,a}$ consecutively in counterclockwise orientation at $v$.
Let $q_{s,+,n}, q_{t,-,n}$ be the corresponding fixed points of $f_n^m$.
Since $f_n$ is a Sierpinski rational map, $\partial \U_{s,n} \cap \partial \U_{t,n} = \emptyset$.
Thus $q_{s,+,n}\neq q_{t,-,n}$, contradicting to Corollary \ref{cor:qeq}.

If $\T_{t,a} = \{v\}$ is trivial, then there exists a critical point $v_n$ of $f_n^m$ with $v_n \to_a v$ by Lemma \ref{lem:critc}.
Since $a$ is a Fatou point, $R_a^m$ has degree $\geq 2$.
Let $c \neq v$ be a critical point of $R_a^m$.
Then there exists critical points $c_n$ of $f_n^m$ with $c_n \to_a c$.
Since $C$ separates $v$ and $c$, $C_n$ separates critical points of $f_n^m$ for sufficiently large $n$, which is a contradiction to Corollary \ref{cor:NSC}.

Consider now that $v$ is aperiodic. Then $w:=R_a^m(v) \in \hat\C_b$ is fixed by $R_b^m$ where $b = F^m(a) \in \RV^\mathcal{F}$.
Since $\T_{s,a}$ is non-trivial, $\T_{r,b}$ is non-trivial where $r = \Phi^m(s)\in |\mathcal{S}|$.
Let $C_a \subseteq \hat\C_a$ and $C_b \subseteq \hat\C_b$ be small simple closed curves around $v$ and $w$ respectively.
We assume that $R_a^m(C_a) = C_b$.

If $\T_{t,a}$ is non-trivial, then $C_a$ intersects $E_s, E_t$ of $\T_{s,a},\T_{t,a}$ consecutively in counterclockwise orientation at $v$.
We can connect the corresponding pre-fixed points $q_{s,+,n}, q_{t,-,n}$ of $f_n^m$ with $J_n$.
By Corollary \ref{cor:qeq}, $R_a^m(J_n)$ is a closed curve homotopically trivial in $\hat\C-(V_n \cup Q_n)$.
Thus, $q_{s,+,n} = q_{t,-,n}$, contradicting to $\partial \U_{s,n} \cap \partial \U_{t,n} = \emptyset$.

If $\T_{t,a} = \{v\}$ is trivial, then there exists $v_n \in \U_{t,n} \cap Q_n$ with $v_n \to_a v$.
Since $a$ is a Fatou point, $R_a^m$ has degree $\geq 2$, so there exists a sequence 
$$
q_n \in f_n^{-m}(Q_n) - \U_{s,n} \cap f_n^{-m}(Q_n)
$$ 
with $q_n\to_a q \neq v$.
By Corollary \ref{cor:HCU}, $C_{b,n}$ is homotopic to a simple closed curve of $\U_{r,n}$ in $\hat\C-(V_n \cup Q_n)$, so $C_{a,n}$ is homotopic to a simple closed curve of $\U_{s,n}$ in $\hat\C-f_n^{-m}(Q_n)$.
Thus, $C_{a,n}$ does not separate the set $f_n^{-m}(Q_n) - \U_{s,n} \cap f_n^{-m}(Q_n)$ for all sufficiently large $n$.
On the other hand, $C_a$ separates $v$ and $q$, so $C_{a,n}$ separates the set $f_n^{-m}(Q_n) - \U_{s,n} \cap f_n^{-m}(Q_n)$ for all sufficiently large $n$, which is a contradiction.

Finally, we show that the choice of $a$ is unique.
Suppose there exists another Fatou point $b$ with $\T_{s,b}$ non-trivial.
Let $x_b \in \Xi_a \subseteq \hat \C_a$ associated to $b$, and $x_a \in \Xi_b \subseteq \hat \C_b$ associated to $a$.
Since $\T_{s,n}$ is connected, $x_b \in \T_{s,a}$ and $x_a \in \T_{s,b}$.

Let $t \neq s \in |\mathcal{S}|$.
Since $\T_{s,b} \cap \T_{t, b} = \emptyset$, $x_a \notin \T_{t, b}$.
Therefore, $\T_{t,a}$ is trivial and $\T_{t, a} = \{x_b\}$. This is a contradiction to $\T_{s,a} \cap \T_{t, a} = \emptyset$.
\end{proof}

\subsection{Reduced rescaling tree maps}
Let $\RT^{r} \subseteq \RT$ be the convex hull of Fatou points $\RP^\mathcal{F}$.
The reason why we make such modification is that the tree of Riemann spheres $(\RT, \hat\C^\RV)$ may `blow up' some degree one periodic points, giving Julia endpoints in $\RT$ and creating non-essential degeneracies.

Let $\RV^{r} := \RV \cap \RT^{r}$ be the set of vertices for $\RT^r$ and $\Pi^r = \Pi \cap \RT^{r}$.
We first show $\RT^{r}$ is invariant under $F$.

\begin{lem}\label{lem:lf}
We have the following
\begin{itemize}
\item If $a\in \RP^\mathcal{J}$, then $a\notin \RT^{r}$. Thus $\Pi^r = \Pi^\mathcal{F}$.
\item If $a\in \RV^r - \RP^r$, then $T_a\RT^r = T_a \RT$. In particular, $a$ is a branch point of $\RT^r$.
\end{itemize}
\end{lem}
\begin{proof}
To prove the first statement, we assume $a=[v]$ wih $v\in \mathcal{V}_s$ for some $s\in |\mathcal{S}|$.
By Assumption \ref{assumption:2}, let $w\in \mathcal{V}_s$ be a Fatou point and $b = [w] \in \RP^\mathcal{F}$.
Let $x_a \in \hat\C_b$ be the associated point of $a$.
Then $\T_{s,b}$ is non-trivial and $x_a \in \T_{s,b}$.

Suppose that $a\in \RT^r$. Then there exists $c \in \RP^\mathcal{F}$ corresponding to the same tangent direction in $T_b \RT^r$ as $a$.
Let $t\in |\mathcal{S}|$ be such that $\T_{t,c}$ is non-trivial.
Then $\T_{t,b} = \{x_a\}$, contradicting to Proposition \ref{prop:dlqig}.

To prove the second statement, suppose that there exists $x\in T_a \RT$ but $x\notin T_a\RT^r$.
Let $b'\in \Pi$ be in the component of $\RT - \{a\}$ associated to the direction of $x$.
Then $b'$ is a Julia point as $x\notin T_a\RT^r$.
Assume $b'= [v]$ for some $v\in \mathcal{V}_s$. 
By Assumption \ref{assumption:2}, let $w\in \mathcal{V}_s$ be a Fatou point and $b = [w] \in \RP^\mathcal{F}$.
Then $b$ is not in the component of $\RT - \{a\}$ in the direction of $x$.
Let $x_a \in \hat\C_b$ be the associated point of $a$.
Then $\T_{s,b}$ is non-trivial and $x_a \in \T_{s,b}$.
Then the same argument as above gives a contradiction.
\end{proof}

\begin{lem}
If $a\in \RV^r$, then $F(a) \in \RV^r$.
\end{lem}
\begin{proof}
If $a\in \RP^r$, then $F(a) \in \RP^r \subset \RV^r$ by Lemma \ref{lem:lf}.
If $a\in \RV^r - \RP$, then it is a branch point.
The same proof of Lemma \ref{lem: extmap} shows that $F(a)$ is a branch point of $\RT^r$. In particular, $F(a) \in \RV^r$.
\end{proof}

\begin{cor}
$\RT^r$ is invariant under $F$.
\end{cor}
We call $\RT^r$ the {\em reduced rescaling tree}, and $F:\RT^r \longrightarrow \RT^r$ the {\em reduced rescaling tree map}.

\begin{lem}\label{lem:ntrt}
If $[f_n]$ diverges, then $\RT^r$ is non-trivial, i.e., it contains more than $1$ point.\footnote{We remark that this is not true for quasi post-critically finite degenerations in other hyperbolic components. There are examples of diverging quasi post-critically finite degenerations with only one Fatou point in the rescaling tree.}
\end{lem}
\begin{proof}
If $\RT^r = \{a\}$, then $F(a) = a$.
By Proposition \ref{prop:dlqig}, $T_{s,a}$ are disjoint and non-trivial for all $s \in |\mathcal{S}|$.
Thus $R_a$ has no holes and $[f_n]$ converges.
\end{proof}

Label the edges of $\RT^r$ by $E_1,..., E_k$, we denote the {\em reduced Markov matrix} by 
$$
M_{i,j}^r = \begin{cases} 1 &\mbox{if } E_i \subseteq F(E_j) \\ 
0 & \mbox{otherwise } \end{cases},
$$
and the {\em reduced degree matrix} by
$$
D_{i,j}^r = \begin{cases} \delta(E_i) &\mbox{if } i = j \\ 
0 & \mbox{otherwise } \end{cases}
$$
where $\delta(E_i)$ is the degree at the edge $E_i$.
The same proof as Proposition \ref{prop:mdto} gives 
\begin{lem}\label{lem:mdto2}
If $\RT^r$ is not-trivial, then there exists a non-negative vector $\vec{v} \neq \vec{0}$ so that
\begin{align*}
M^r \vec v = D^r \vec v.
\end{align*}
\end{lem}

\subsection{Thurston's criterion and the proof of Theorem \ref{thm:schcb}}
Let $f: \hat\C \longrightarrow \hat\C$ be a rational map with post-critical set $P_f$.
Let $P_f \subseteq U$ be a forward invariant set, i.e., $f(U) \subseteq U$. 
A simple closed curve $\sigma$ on $\hat \C - U$ is {\em essential} if it does not bound a disk in $\hat C - U$, and a curve is {\em peripheral} if it encloses a single point of $U$.
Two simple curves are {\em parallel} if they are homotopic in $\C- U$.

A {\em curve system} $\Sigma = \{\sigma_i\}$ in $\hat\C - U$ is a finite nonempty collection of disjoint simple closed curves, each essential and non-peripheral, and no two parallel.
A curve system determines a {\em transition matrix} $A(\Sigma) : \R^\Sigma \longrightarrow \R^{\Sigma}$ by the formula
$$
A_{\sigma\tau} = \sum_\alpha \frac{1}{\deg(f: \alpha \to \tau)}
$$
where the sum is taken over components $\alpha$ of $f^{-1}(\tau)$ isotopic to $\sigma$.

Let $\lambda(\Sigma)\geq 0$ denote the spectral radius of $M(\Sigma)$.
Since $A(\Sigma) \geq 0$, the Perron-Frobenius theorem guarantees that $\lambda(\Sigma)$ is an eigenvalue for $A(\Sigma)$ with a non-negative eigenvector.

We need the following modification of the necessary direction of Thurston's criterion, which follows from the same proof of Theorem B.4 in \cite{McM94}.
\begin{prop}\label{prop:mtc}
Let $f$ be a hyperbolic Sierpinski carpet rational map.
Let $\U_f$ be the union of critical and post-critical Fatou components.
Let $\Sigma$ be a curve system in $\C - \U_f$, then $\lambda(\Sigma) < 1$. 
\end{prop}

We will now construct a curve system $\Sigma$ in $\hat\C - \U_n$ for $f_n$, where each curve $\sigma$ corresponds to an edges of $\RT^r$.

Let $E=[a,b]$ be an edge of $\RT^r$.
Let $x_b \in \hat\C_a$ and $x_a \in \hat\C_b$ be the point associated to $b$ and $a$ respectively.
By Proposition \ref{prop:dlqig}, let $C_a \subseteq \hat\C_a$ and $C_b \subseteq \hat \C_b$ be a small circle around $x_b$ and $x_a$ so that $C_a$ and $C_b$ do not cut $\T_a$ and $\T_b$ respective.
Modify $C_a$ and $C_b$ if necessary, we assume $R_a(C_a)$ and $R_b(C_b)$ are simple closed curves.
Denote $D_a$ and $D_b$ be the disk bounded by $C_a$ and $C_b$ that contain $x_b$ and $x_a$ respectively.
Shrinking $C_a, C_b$ if necessary, we assume $D_a, D_b$ contain no holes of $R_a, R_b$ except possibly at $x_b, x_a$.

Since $a, b$ bound an edge, either $\T_{s, a} \in \hat\C_a - \overline{D_a}$ or $\T_{s, b} \in \hat\C_a - \overline{D_b}$ for any $s\in |\mathcal{S}|$.
Shrinking $C_a, C_b$ if necessary, we assume for sufficiently large $n$, there is an annulus $\mathcal{A}_{E,n}$ bounded by $C_{a,n}$ and $C_{b,n}$ such that
\begin{enumerate}
\item $A_{a,n}^{-1}(C_{a,n})$ and $A_{b,n}^{-1}(C_{b,n})$ converges to $C_a\subseteq\hat\C_a$ and $C_b \subseteq\hat\C_b$ in Hausdorff topology; 
\item $f_n(C_{a,n})$ and $f_n(C_{b,n})$ are simple closed curves;
\item $\mathcal{A}_{E,n}$ contains no critical points of $f_n$;
\item Any essential simple closed curve $\alpha \subseteq \mathcal{A}_{E,n}$ does not separate $V_n \cup Q_n$ in $\U_{s, n}$ for any $s\in |\mathcal{S}|$.
\end{enumerate}
By condition (2) and (3), $f_n(\mathcal{A}_{E,n})$ is again an annulus.
By condition (4), any essential simple closed curve of $\mathcal{A}_{E,n}$ is homotopic to a simple closed curve $\sigma_{E, n} \subseteq \hat\C - \U_n$ in $\hat\C - (V_n \cup Q_n)$.
Since each component $\U_{s,n}$ of $\U_n$ contains a point in $Q_n$, the homotopy class of $\sigma_{E, n}$ in $\hat C - \U_n$ is well defined.

We call $\sigma_{E, n}$ the corresponding simple closed curve for $E$, although technically speaking, only the homotopy class is well defined.
Let 
$$
\Sigma_n =\{\sigma_{E,n}: E \text{ is an edge of } \RT^r\}.
$$

\begin{lem}\label{lem:CSC}
For sufficiently large $n$, $\Sigma_n$ is a curve system in $\hat\C - \U_n$.
\end{lem}
\begin{proof}
We first show that $\sigma_{E,n}$ is essential and non-peripheral.
Note that there exist $a, b \in \RP^r = \RP^\mathcal{F}$ that are in two different components of $\RT^r - \Int(E)$.
There exist $s, t \in |\mathcal{S}|$ so that $\T_{s,a}\subseteq \hat\C_a$, $\T_{t,b}\subseteq \hat\C_b$ are non-trivial.
Thus, for sufficiently large $n$, $\sigma_{E,n}$ separates $\T_{s,n}$ and $\T_{t,n}$.
Thus, $\sigma_{E_n}$ is essential and non-peripheral.

We now show different simple closed curves are not parallel.
Let $E_1, E_2$ be two edges, then there exists $a\in \RV^r$ separating $\Int(E_1)$ and $\Int(E_2)$.

If $a \in \RP^r$, then there exists $s \in |\mathcal{S}|$ so that $\T_{s,a}\subseteq \hat\C_a$ is non-trivial. Then for sufficiently large $n$, the annulus bounded by $\sigma_{E_1,n}$ and $\sigma_{E_2,n}$ contains $\T_{s,n}$, so $\sigma_{E_1,n}$ is not parallel to $\sigma_{E_2,n}$.

If $a \in \RV^r - \RP^r$, then it is a branch point by Lemma \ref{lem:lf}. Thus there exists $b \in \RP^r$ corresponding to a different direction in $T_a \RT^r$ other than those corresponding to $\Int(E_1)$ and $\Int(E_2)$.
Then a similar argument as above shows that $\sigma_{E_1,n}$ is not parallel to $\sigma_{E_2,n}$.
\end{proof}

\begin{lem}\label{lem:lto}
If $E \subseteq F(E')$, then for sufficiently large $n$, $\sigma_{E, n}$ has a lift of degree $\delta(E')$ homotopic to $\sigma_{E',n}$ in $\hat\C-\U_n$.
\end{lem}
\begin{proof}
Let $\mathcal{A}_{E,n}$ and $\mathcal{A}_{E',n}$ be annuli associated to $E$ and $E'$ respectively.
Modify the boundaries of $\mathcal{A}_{E,n}$ if necessary, we may assume $\mathcal{A}_{E,n} \subseteq f_n(\mathcal{A}_{E',n})$, where the inclusion induces an isomorphism on the fundamental group.
So for sufficiently large $n$, there exists an essential simple closed curve $\gamma_n \subseteq \mathcal{A}_{E,n}$ that has a degree $\delta(E)$ lift $\gamma_n' \subseteq \mathcal{A}_{E',n}$.
Since $\gamma_n$ does not separate $V_n \cup Q_n$ in $\U_{s, n}$ for any $s\in |\mathcal{S}|$, by Lemma \ref{lem:isc}, $\sigma_{E,n}$ has a degree $\delta(E)$ lift $\sigma_n' \subseteq \hat\C - \U_n$ homotopic to $\sigma_{E',n}$ in $\hat\C - \U_n$.
\end{proof}

\begin{lem}\label{lem:srg1}
Let $\RT^r$ be non-trivial, and let $\Sigma_n$ be the curve system associated with the edges of $\RT^r$ in $\hat\C - \U_n$.
Then for sufficiently large $n$, the spectral radius $\lambda(\Sigma_n) \geq 1$.
\end{lem}
\begin{proof}
Let $M^r$ and $D^r$ be the reduced Markov matrix and reduced degree matrix.
By Lemma \ref{lem:lto}, the transition matrix $A_n=A(\Sigma_n)$ satisifies
$$
A_n \geq (D^r)^{-1}M^r,
$$
where the inequality means each entry of $A$ is greater or equal to the corresponding entry of $(D^r)^{-1}M^r$.
Since both matrices are non-negative, the spectral radius 
$$
\lambda(\Sigma_n) = \lambda(A_n) \geq \lambda((D^r)^{-1}M^r).
$$
By Lemma \ref{lem:mdto2}, $\lambda((D^r)^{-1}M^r) \geq 1$, so $\lambda(\Sigma_n) \geq 1$.
\end{proof}

We now prove Theorem \ref{thm:schcb}:
\begin{proof}[Proof of Theorem \ref{thm:schcb}]
Suppose for contradiction that $[f_n]$ diverges.
Then by Lemma \ref{lem:ntrt}, $\RT^r$ is non-trivial.
Thus by Lemma \ref{lem:srg1}, for sufficiently large $n$, there exists a curve system $\Sigma_n$ in $\hat\C- \U_n$ with spectral radius $\lambda(\Sigma_n) \geq 1$, contradicting to Proposition \ref{prop:mtc}.
\end{proof}

\section{A double limit theorem for $\QH_d$}\label{sec:ptq}
In this section, we study the convergence of simultaneous degeneration on both Fatou components in $\QH_d$ and prove Theorem \ref{thm:dlm}.

A hyperbolic rational map $f$ is said to be a {\em quasi-Blaschke product} if $J(f)$ is a Jordan curve and $f$ fixes the two Fatou components.
Let $\QH_d \subseteq \mathcal{M}_{d,fm}$ denote the space of conjugacy classes of degree $d$ quasi-Blaschke products. 

The corresponding mapping scheme $\mathcal{S} = (|\mathcal{S}|, \Phi, \delta)$ has rather simple description: it consists of two fixed points of degree $d$.
To simplify the notations,
we identify the Blaschke model space $\BP^\mathcal{S}$ with $\BP_d\times \BP_d$, and $\bp\in \BP^\mathcal{S}$ with $(\bp_+, \bp_-)$.
Thus, a marked quasi-Blaschke product $f$ is obtained by gluing together a pair of maps $\bp_+, \bp_-\in \BP_d$ using their markings on $\mathbb{S}^1$.

Let $[f_n] \in \QH_d$ be quasi post-critically finite.
Since $[f_n] \in \mathcal{M}_{d,\fm}$, the dynamics $\bp_{\pm,n}$ induces two {\em markings}
$$
\eta_{\pm,n} : \R/\Z \longrightarrow \mathcal{J}(f_n)
$$
conjugating $m_d(t) = dt$ with $f_n$.
The two markings satisfy
$$
\eta_{+,n}(t) = \eta_{-,n}(-t).
$$
We denote $\eta_n = \eta_{+,n}$ and call it the {\em marking} on $\mathcal{J}(f_n)$.
The pre-periodic points on $\mathcal{J}(f_n)$ are thus labeled by a rational angle in $\Q/\Z$.

Recall that a {\em lamination} $\mathcal{L}$ is a family of disjoint hyperbolic geodesics in $\D$ together with the two end points in $\mathbb{S}^1\cong \R/\Z$, whose union $|\mathcal{L}| :=\bigcup\mathcal{L}$ is closed.
The hyperbolic geodesics are called the {\em leaves} of the lamination.

We say two laminations are {\em parallel} if there exists
a simple closed loop in the union $|\mathcal{L}| \cup \overline{|\mathcal{L}'|}$ containing leaves in both laminations, where $\overline{|\mathcal{L}'|}$ is the image of $|\mathcal{L}'|$ under $z \mapsto \frac{1}{z}$.

A lamination $\mathcal{L}$ generates an equivalence relation $\sim_\mathcal{L}$ on $\mathbb{S}^1$ by $a\sim_{\mathcal{L}} b$ if $a,b \in \mathbb{S}^1$ are connected by a finite chain of leaves in $\mathcal{L}$.
Using this language, $\mathcal{L}$ and $\mathcal{L}'$ are parallel if there exist $a_0,b_0,..., a_{k-1}, b_{k-1} \in \mathbb{S}^1$ with $a_i\sim_{\mathcal{L}} b_i$ and $-b_i \sim_{\mathcal{L}'} -a_{i+1}$ for all $i$, where the indices are modulo $k$ (cf. Definition 1.4 in \cite{BD18}).


\subsection*{Sketch of the proof}
We first sketch the idea of the proof.
By Theorem \ref{thm:crm}, $[f_n]$ converges to a geometrically finite rational map $(F, R)$ on a tree of Riemann spheres $(\RT, \hat\C^\RV)$.
If $[f_n]$ diverges, then the partition of the limiting quasi-invariant graphs $\T_{\pm}$ into $\bigcup_{a} \T_{\pm, a}$ gives two non-trivial partitions of angles 
$$
\R/\Z = \bigcup_{a} I_{\pm, a},
$$
where $I_{\pm, a}$ is a finite union of closed intervals, and the boundaries of $I_{\pm, a}$ are `formed' by leaves of $\mathcal{L}_\pm$.
We show that for all but finitely many rational angles $t\in \pm I_{\pm, a}$, the rescaled limit $\eta_n(t) \to_a \hat \C_a - \Xi_a$.
This implies the two partitions are conjugate, i.e., $I_{+, a} = -I_{-,a}$, showing $\mathcal{L}_{\pm}$ are parallel.

To set up the notations,
we use $+$ and $-$ to denote the two points in the mapping scheme $|\mathcal{S}|$.
We use $\pm$ in a statement to mean the statement holds if we replace {\em all} $\pm$ by $+$ or {\em all} $\pm$ by $-$.

\subsection{Partition of the pre-periodic points}
Let $\mathcal{T}_\pm$ be the quasi-invariant trees.
Let $m$ be the largest pre-period of the edges of $\mathcal{T}_\pm$.
By replacing $\mathcal{T}_\pm$ by a subtree in the $m$-th pull back $\mathcal{T}_\pm^m$ \footnote{We also need to construct the rescaling tree $\RT$ using the corresponding replacement.}, we assume the following technical `balanced tree' assumption holds throughout this section:
\begin{assumption}\label{assumption:1}
All edges of pre-period $\leq m$ in $\mathcal{T}_\pm^\infty$ are contained in $\mathcal{T}_\pm$.\footnote{Since $\mathcal{T}_\pm$ contains all critical points, it is easy to check that the union of all edges of pre-period $\leq m$ in $\mathcal{T}^\infty$ form a tree, i.e, it is connected.}
\end{assumption}

Let $v\in \mathcal{V}_\pm$. The angles landing at $v$ are defined using the Blaschke product model
$$
I_v := \{t\in \R/\Z: \lim_v \eta_{\pm, n}^\mathcal{B}(t) \in \mathbb{S}_v^1- \xi_{\pm, v}(T_v\mathcal{T}_\pm)\},
$$
where $\eta_{\pm, n}^\mathcal{B}: \R/\Z \rightarrow \mathbb{S}^1$ is the marking for the Blaschke products $\bp_\pm$.
Here we use superscript $\mathcal{B}$ to distinguish from the markings of $f_n$.
Under the Assumption \ref{assumption:1}, we have
\begin{lem}\label{lem:JF}
If $v \in \mathcal{V}_\pm^\mathcal{J}$ is a Julia point, then $I_v$ is finite.
\end{lem}
\begin{proof}
If $v$ is a periodic point of period $q$, then the rescaling limit $\mathcal{R}^q_v: \D_v \longrightarrow \D_v$ has degree $1$.
By Lemma \ref{lem:critc}, the holes are contained in $\xi_v(T_v \mathcal{T}_\pm)$. 
Thus, there are at most finitely many periodic points $t\in \R/\Z$ with $\lim_v \eta_{\pm, n}^\mathcal{B}(t) \notin \xi_v(T_v\mathcal{T}_\pm)$.
Therfore $I_v$ is finite.

If $v$ is aperiodic with pre-period $l$, then let $\mathcal{R}^l_v: \D_v \longrightarrow \D_{F_\pm^l(v)}$ be the corresponding rescaling limit.
By Assumption \ref{assumption:1}, 
$$
\xi_{\pm,v}(T_v\mathcal{T}_\pm) = (\mathcal{R}^l_v)^{-1} (\xi_{\pm,F_\pm^l(v)}(T_{F_\pm^l(v)}\mathcal{T}_\pm)) \subseteq \mathbb{S}^1_v.
$$
Since $I_{F_\pm^l(v)}$ is finite, by Lemma \ref{lem:lacs}, $I_v$ is finite as well.
\end{proof}

\begin{figure}[ht]
  \centering
  \resizebox{1\linewidth}{!}{
    \def\svgwidth{\columnwidth}
    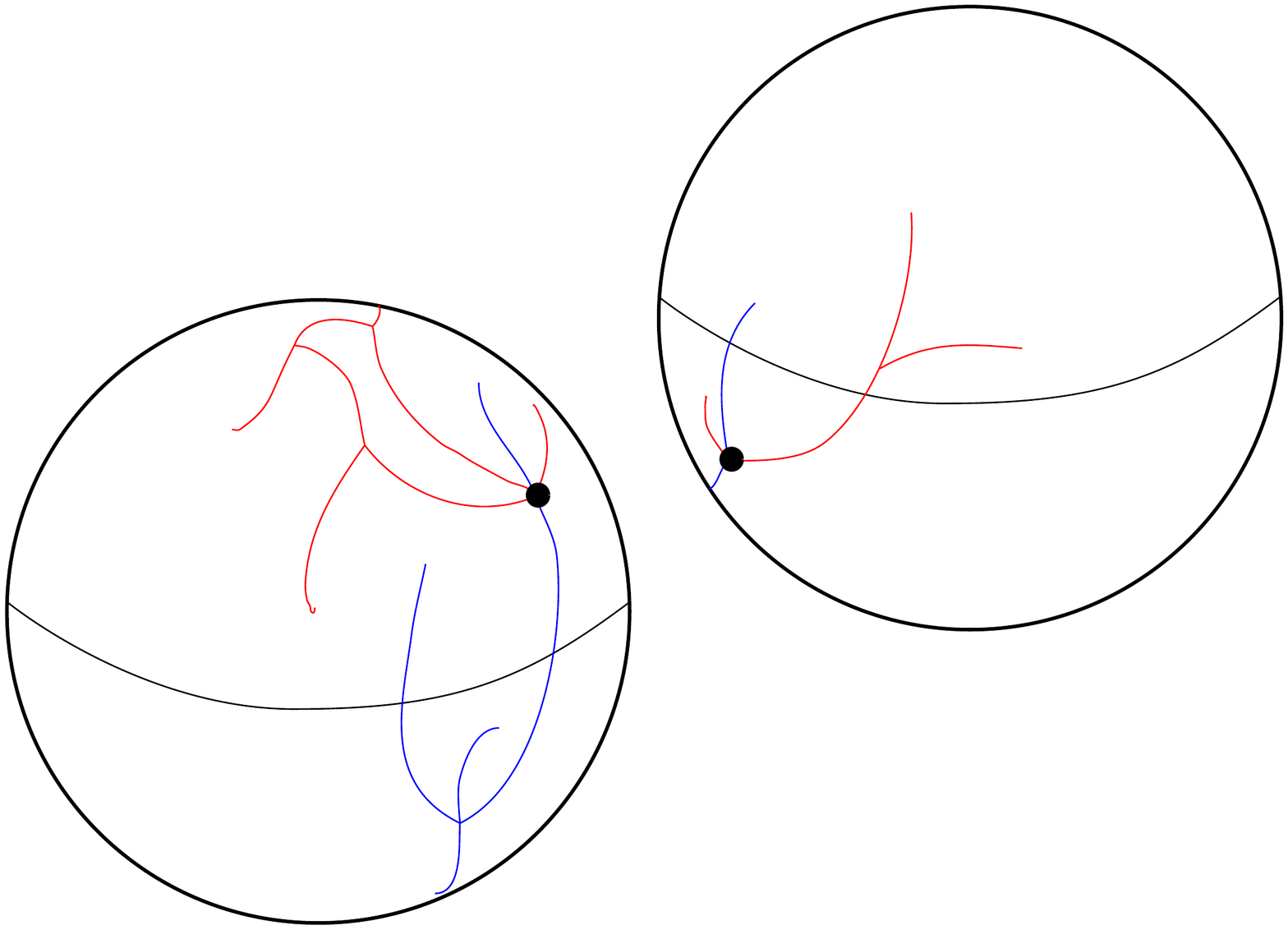

    \def\svgwidth{\columnwidth}
\begingroup%
  \makeatletter%
  \providecommand\color[2][]{%
    \errmessage{(Inkscape) Color is used for the text in Inkscape, but the package 'color.sty' is not loaded}%
    \renewcommand\color[2][]{}%
  }%
  \providecommand\transparent[1]{%
    \errmessage{(Inkscape) Transparency is used (non-zero) for the text in Inkscape, but the package 'transparent.sty' is not loaded}%
    \renewcommand\transparent[1]{}%
  }%
  \providecommand\rotatebox[2]{#2}%
  \newcommand*\fsize{\dimexpr\f@size pt\relax}%
  \newcommand*\lineheight[1]{\fontsize{\fsize}{#1\fsize}\selectfont}%
  \ifx\svgwidth\undefined%
    \setlength{\unitlength}{841.88976378bp}%
    \ifx\svgscale\undefined%
      \relax%
    \else%
      \setlength{\unitlength}{\unitlength * \real{\svgscale}}%
    \fi%
  \else%
    \setlength{\unitlength}{\svgwidth}%
  \fi%
  \global\let\svgwidth\undefined%
  \global\let\svgscale\undefined%
  \makeatother%
  \begin{picture}(1,0.70707071)%
    \lineheight{1}%
    \setlength\tabcolsep{0pt}%
    \put(0,0){\includegraphics[width=\unitlength,page=1]{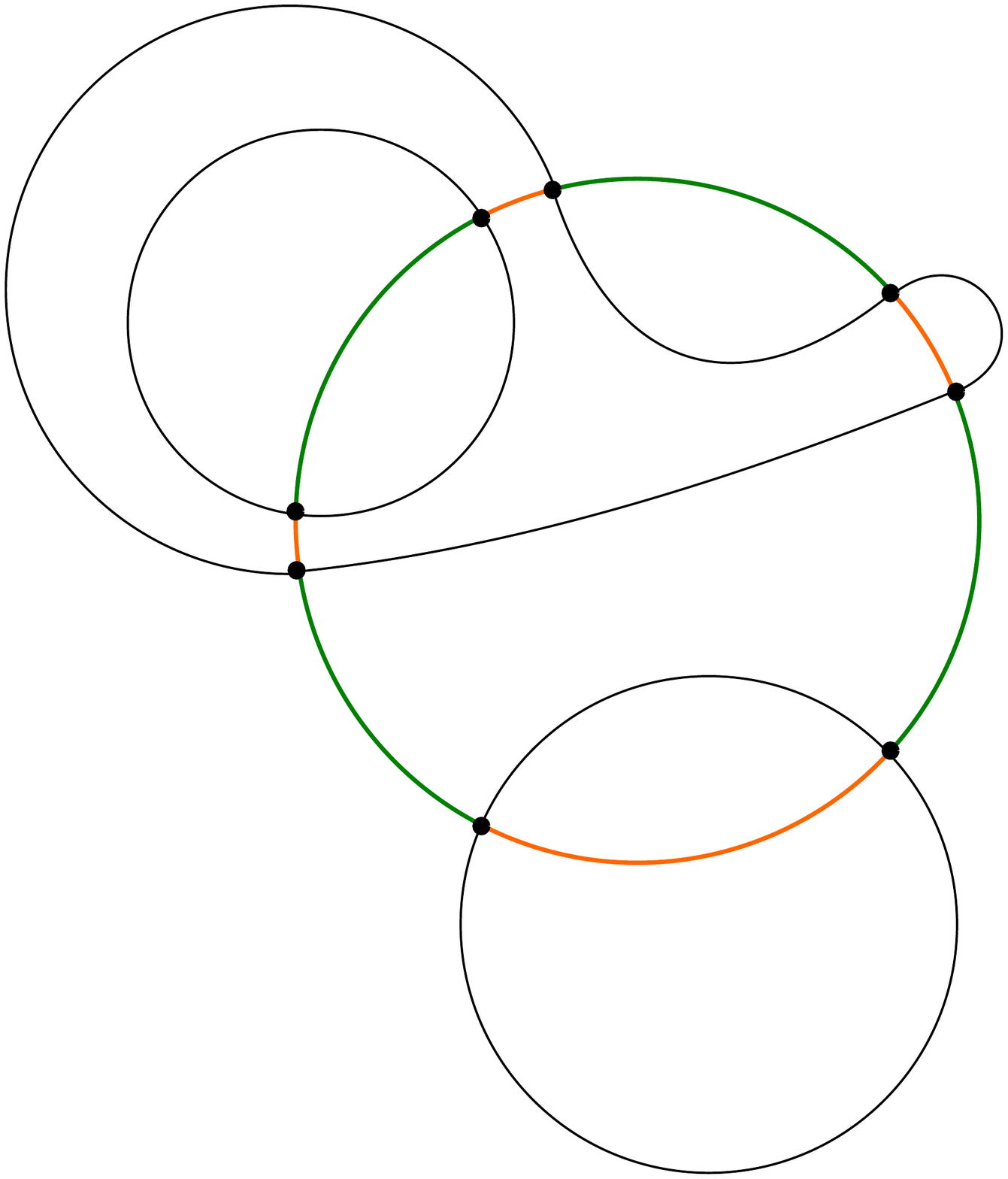}}%
    \put(0.5502924,0.60933142){\color[rgb]{0,0,0}\makebox(0,0)[lt]{\lineheight{1.25}\smash{\begin{tabular}[t]{l}$I_{+,a} = -I_{-,a}$\end{tabular}}}}%
    \put(0.49181281,0.1546226){\color[rgb]{0,0,0}\makebox(0,0)[lt]{\lineheight{1.25}\smash{\begin{tabular}[t]{l}$I_{+,b} = -I_{-,b}$\end{tabular}}}}%
  \end{picture}%
\endgroup%

  }
  \caption{On the left, a schematic picture of the graphs $\T_{+}$ and $\T_{-}$ are drawn in blue and red. The singular points $x_a, x_b$ cut the graphs into finitely many components giving partitions of the circle $\mathbb{S}^1$ in orange and green as appeared on the right.}
  \label{fig:Par}
\end{figure}

Recall that each vertex $v\in \mathcal{V}_\pm$ gives a rescaling $A_{v,n} \in \PSL_2(\C)$.
Two rescalings are equivalent if $A_{v,n}^{-1} \circ A_{w,n}$ stay in a bounded subset of $\PSL_2(\C)$, and $\RP = \mathcal{V}_\pm / \sim$ (see Definition \ref{defn:rld}).
Let $a\in \RP^\mathcal{F}$.
We define the angles at $a$ as
$$
I_{\pm, a} := \bigcup_{v\in \mathcal{V}_\pm, \; [v]=a} \overline{ \Int(I_{\pm, v})} \subseteq \R/\Z.
$$
Since each $ \overline{ \Int( I_{\pm, v})}$ is a finite union of closed intervals where the boundaries of any complementary component form a leaf of $\mathcal{L}_\pm$,
for each $x\in \partial I_{\pm, a}$, there exists $y\neq x\in \partial I_{\pm, a}$ connected to $x$ by a finite chain of leaves in $\mathcal{L}_\pm$ (see Figure \ref{fig:Par}).

We first show this $I_{\pm, a}$ induces partions of $\R/\Z$:
\begin{lem}
If $a\neq b\in \RP^\mathcal{F}$, then $\Int(I_{\pm, a}) \cap \Int(I_{\pm, b}) = \emptyset$.
Moreover, 
$\bigcup_{a\in \RP^\mathcal{F}} I_{\pm, a} = \R/\Z$.
\end{lem}
\begin{proof}
Since $I_{\pm, v}$ are disjoint, $\Int(I_{\pm, a}) \cap \Int(I_{\pm, b}) = \emptyset$.
Since $\bigcup_{v\in \mathcal{V}_\pm} I_{\pm, v}$ is dense in $\R/\Z$, by Lemma \ref{lem:JF}, $\bigcup_{v\in \mathcal{V}_\pm^\mathcal{F}} I_{\pm, v}$ is dense as well.
If $v\in \mathcal{V}_\pm^\mathcal{F}$, then $[v] \in \RP^\mathcal{F}$.
Thus, $\bigcup_{a\in \RP^\mathcal{F}} I_{\pm, a} = \R/\Z$.
\end{proof}

After passing to a subsequence, and a diagonal argument, we assume the limit 
$$
x_{t,a} := \lim_a \eta_n(t) = \lim_{n\to\infty} A_{a,n}^{-1}(\eta_n(t))
$$ 
exists for any rational angles $t \in \Q/\Z$ and any $a\in \RP^\mathcal{F}$.
Note that $x_{t,a} = \lim_a \eta_{\pm,n}(\pm t)$.

The following is the key of the argument, showing that the two partitions with respect to $+$ and $-$ are conjugate:
\begin{prop}\label{prop:pc}
Let $a\in \RP^\mathcal{F}$.
For all but finitely many rational angles $t\in \pm I_{\pm, a}$, $x_{t, a} \in \hat \C_a - \Xi_a$.
Thus, $I_{+,a} = -I_{-,a}$ for all $a\in \RP^\mathcal{F}$.
\end{prop}

Once the above is established, Theorem \ref{thm:dlm} follows immediately:
\begin{proof}[Proof of Theorem \ref{thm:dlm} assuming Proposition \ref{prop:pc}]
Suppose $[f_n]$ is diverging. 
Then we get non-trivial partition of the angles 
$$
\R/\Z = \bigcup_{a\in \RP^\mathcal{F}} I_{\pm, a},
$$
where each $I_{\pm, a}$ is the complement of a finite union of open intervals.
By Proposition \ref{prop:pc}, $I_{+,a} = -I_{-,a}$ for all $a\in \RP^\mathcal{F}$.
Thus, the two laminations $\mathcal{L}_\pm$ are parallel.
\end{proof}

\subsection{Conjugate partitions}
The remaining of the section is devoted to the proof of Proposition \ref{prop:pc}.
We first show, in contrast to disjoint limiting quasi-invariant graphs for degenerations of Sierpinski carpet rational maps, any singular point must intersect the limiting quasi-invariant graphs.
\begin{lem}\label{lem:spilqig}
Let $a\in \RP^\mathcal{F}$, and $x\in \Xi_a$.
Then there exists an edge of $\T_a = \T_{+, a} \cup \T_{-, a}$ incident at $x$.
\end{lem}
\begin{proof}
Suppose not. 
Then we can find a small closed curve $C'$ around $x$ that is disjoint from $\T_a$.
Since $a\in \RP^\mathcal{F}$, either $\T_{+, a}$ or $\T_{-, a}$ is non-trivial.
Without loss of generality, we assume $\T_{+, a}$ is non-trivial.
Since $x\in \Xi_a$ and $x\notin \T_{+, a}$, $\T_{-, a} = \{x\}$.

Thus there exists an annulus $\mathcal{A}_n \subseteq \hat\C$ separating the critical values in $\U_{\pm, n}$ whose modulus $m(\mathcal{A}_n) \to \infty$.
By enlarging the annulus, we assume the boundary components of $\mathcal{A}_n$ are in $\U_{\pm, n}$ respectively.
Thus, its preimage $\mathcal{A}_n'$ is compactly contained in $\mathcal{A}_n$ as an essential annulus with $f_n : \mathcal{A}_n' \longrightarrow \mathcal{A}_n$ as a degree $d$ covering map.
We can find $\mathcal{B}_n' \subseteq \mathcal{B}_n \subseteq \mathcal{A}_n$ with modulus $K \leq m(\mathcal{B}_n) < dK$ for some $K$, and $f_n : \mathcal{B}_n' \longrightarrow \mathcal{B}_n$ is a degree $d$ covering map.
By Theorem 2.1 in \cite{McM94}, $\mathcal{B}_n$ contains a round annulus $\widetilde{\mathcal{B}_n}$ of modulus $K/2$ for sufficiently large $K$.
Define $B_n \in \PSL_2(\C)$ so that $B_n (\mathcal{A}(K/2)) = \widetilde{\mathcal{B}_n}$, where $\mathcal{A}(K/2) \subseteq B(0,1)$ is the round annulus of modulus $K/2$ with one boundary component $\partial B(0,1)$.
Thus $B_n^{-1} \circ f_n \circ B_n$ converges to a degree $d$ map.
Since $a\in \RP^\mathcal{F}$, $M_n^{-1} \circ A_{a,n}$ is bounded.
So $R_a$ has degree $d$, and $x$ is a fixed point of degree $d$.
This is a contradiction to $x\in \Xi_a$.
\end{proof}
As an immediate corollary, we have
\begin{cor}\label{cor:ncpcs}
There are no critical periodic cycles of $R$ in the singular set $\Xi$.
\end{cor}
\begin{proof}
Let $x$ be a periodic singular point.
By Lemma \ref{lem:spilqig}, $x$ is a frontier vertex for some non-trivial limiting quasi-invariant graph. By Lemma \ref{lem:ncfh}, $x$ is not in a critical periodic cycle.
\end{proof}

A {\em ray} $\alpha$ in an infinite tree $T$ is a semi-infinite simple path, i.e., it is a sequence of vertices $v_0,v_1,v_2,...$ in which each vertex appears at most once and two consecutive vertices bounds an edge in $T$.
Two rays $\alpha_1, \alpha_2$ are {\em equivalent} if $\alpha_1 \cap \alpha_2$ is again a ray.
An equivalence class of rays is called an {\em end} of $T$, and we use $\epsilon(T)$ to denote the set of ends of $T$.
We shall use $\alpha$ to represent both the end and the ray it represents.

Let $F_\pm: \mathcal{T}_\pm^\infty \longrightarrow \mathcal{T}_\pm^\infty$.
Then any end can be represented by a ray disjoint from the critical points, whose image under $F$ is again a ray.
Therefore, we have an induced map $F_{\pm,*}: \epsilon(T) \longrightarrow \epsilon(T)$.
We call an end is periodic and pre-periodic if it is periodic or pre-periodic under this induced map.

Let $\rl_{\pm, v} : \D_v \longrightarrow \D_{F(v)}$ be the rescaling limit of $\bp_{\pm,n}$ at $v$.
Let $\rl_\pm := \bigcup_v \rl_{\pm, v}$ be the union of rescaling limits. We have
\begin{lem}\label{lem:pppe}
Given $t\in \Q/\Z$ that is not the end point of a leaf of $\mathcal{L}_\pm$, then $t$ corresponds to 
\begin{enumerate}
\item A pre-periodic point of $\mathbb{S}^1_v$ whose orbit avoids the holes of $\rl_\pm$; or
\item A pre-periodic end of $\mathcal{T}_\pm^\infty$.
\end{enumerate}
\end{lem}
\begin{proof}
Since $t$ is not the end point of a leaf of $\mathcal{L}_\pm$, $t$ lands on some vertex $v\in\mathcal{V}_\pm$.
If $x:=\lim_v \eta^\mathcal{B}_{\pm, n}(t)$ is a not a hole for any iterates of $\rl_{\pm,v}$, then we have the first case.
Otherwise, by the pull back construction, $t$ corresponds to a pre-periodic end for $\mathcal{T}_\pm^\infty$.
\end{proof}

First consider the case $t\in \Q/\Z$ is not the end point of a leaf of $\mathcal{L}_\pm$.
In this case, we call $t$ of Type (1) or Type (2) depending on the corresponding case in Lemma \ref{lem:pppe}.
\begin{lem}\label{lem:PV1}
Let $a\in \RP^\mathcal{F}$.
If $t\in \pm I_{\pm, a}$ is rational angle of Type (1), then with finitely many exceptions, $x_{t,a} \in \hat \C - \Xi_a$.
\end{lem}
\begin{proof}
We first assume $t$ is periodic under $m_d$ of period $q$.
Let $a=[v]$, then $A_{a,n}^{-1}((\U_{\pm,n}, \phi_{\pm,n}(v_n)))$ converges in Carath\'eodory limit to $(\U_{\pm, v}, v) \subseteq \hat \C_a$.
Suppose the period of $v$ is $k$, then $R_a ^k : \U_{\pm, v} \longrightarrow \U_{\pm, v}$ is conformally conjugate to the rescaling limit $\mathcal{R}_{\pm, v}^k: \D_v \longrightarrow \D_v$.

Let $u_\pm \in \mathbb{S}^1_v$ associated to $t$. 
Let $\hat \gamma_\pm \subseteq \D_v$ be an arc landing at $u_\pm$.
We assume it is invariant: $\mathcal{R}_{\pm, v}^q(\hat\gamma_\pm \cap U) = \hat\gamma_\pm \cap \mathcal{R}_{\pm, v}^q(U)$ for any sufficient small neighborhood of $u_\pm\in U$.
Let $\hat \gamma_{\pm, n}$ be a sequence of invariant rays landing at $\eta_{\pm,n}^\mathcal{B}(t)$ converging to $\hat\gamma_\pm$ in $v$-coordinate.

Let $\gamma_{\pm, n} \subseteq \U_{\pm, n}$ be the corresponding invariant ray for $f_n$.
Let $\gamma_\pm$ be the Hausdorff limit of $A_{a,n}^{-1}(\gamma_{\pm,n})$\footnote{Note that if $x\in \gamma_{\pm, n}$, then $d_{\U_{\pm, n}}(f_n^q(x),x)$ is uniformly bounded by $K$ for some $K$, so a similar argument in Theorem \ref{thm:hl} gives that $\gamma_\pm$ is a graph with finitely many vertices and finitely many non-loop edges.}, and $\gamma_{\pm,v} = \gamma_\pm \cap \U_{\pm,v}$.
Then $\gamma_{\pm,v}$ lands at a periodic point $y_\pm \in \partial \U_{\pm, v}$.
By Corollary \ref{cor:ncpcs}, the periodic singular points in $\Xi_a$ are not in critical periodic cycle.
So there are only finitely many such invariant rays landing on $\Xi_a$.
Similarly, there are only finitely many invariant rays landing at parabolic fixed points.
Thus, with finitely many exceptions, we assume $y_\pm \in\hat\C_a- \Xi_a$ is a repelling periodic point.
Since $\Xi$ and $\gamma_{\pm,v}$ are invariant, the iterates of $y_\pm$ avoid the singular set and thus avoid the holes of $R$.

\begin{figure}[ht]
  \centering
  \resizebox{0.6\linewidth}{!}{
    \def\svgwidth{\columnwidth}
\begingroup%
  \makeatletter%
  \providecommand\color[2][]{%
    \errmessage{(Inkscape) Color is used for the text in Inkscape, but the package 'color.sty' is not loaded}%
    \renewcommand\color[2][]{}%
  }%
  \providecommand\transparent[1]{%
    \errmessage{(Inkscape) Transparency is used (non-zero) for the text in Inkscape, but the package 'transparent.sty' is not loaded}%
    \renewcommand\transparent[1]{}%
  }%
  \providecommand\rotatebox[2]{#2}%
  \newcommand*\fsize{\dimexpr\f@size pt\relax}%
  \newcommand*\lineheight[1]{\fontsize{\fsize}{#1\fsize}\selectfont}%
  \ifx\svgwidth\undefined%
    \setlength{\unitlength}{1020.47244094bp}%
    \ifx\svgscale\undefined%
      \relax%
    \else%
      \setlength{\unitlength}{\unitlength * \real{\svgscale}}%
    \fi%
  \else%
    \setlength{\unitlength}{\svgwidth}%
  \fi%
  \global\let\svgwidth\undefined%
  \global\let\svgscale\undefined%
  \makeatother%
  \begin{picture}(1,0.5)%
    \lineheight{1}%
    \setlength\tabcolsep{0pt}%
    \put(0,0){\includegraphics[width=\unitlength,page=1]{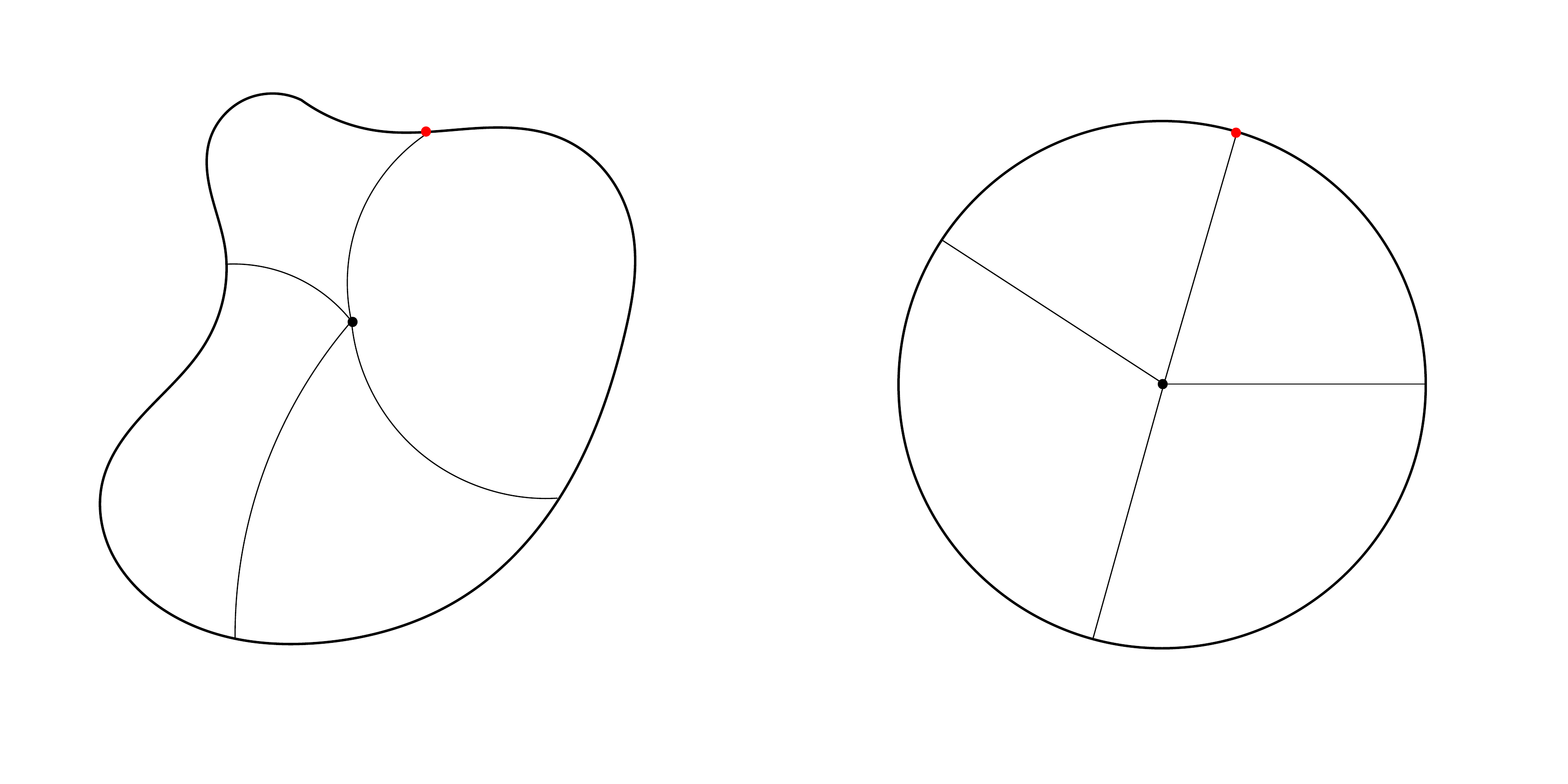}}%
    \put(0.83818698,0.41610707){\color[rgb]{0,0,0}\makebox(0,0)[lt]{\lineheight{1.25}\smash{\begin{tabular}[t]{l}{\large$\D_v$}\end{tabular}}}}%
    \put(0.35338972,0.42353122){\color[rgb]{0,0,0}\makebox(0,0)[lt]{\lineheight{1.25}\smash{\begin{tabular}[t]{l}{\large$\U_{\pm, v}$}\end{tabular}}}}%
    \put(0,0){\includegraphics[width=\unitlength,page=2]{LandingL.pdf}}%
    \put(0.73276393,0.39606185){\color[rgb]{0,0,0}\makebox(0,0)[lt]{\lineheight{1.25}\smash{\begin{tabular}[t]{l}$\hat\gamma_\pm$\end{tabular}}}}%
    \put(0.20045213,0.40422842){\color[rgb]{0,0,0}\makebox(0,0)[lt]{\lineheight{1.25}\smash{\begin{tabular}[t]{l}$\gamma_\pm$\end{tabular}}}}%
  \end{picture}%
\endgroup%

  }
  \caption{The red points on the left and right correspond to $y_\pm$ and $u_\pm$ in the proof of Lemma \ref{lem:PV1}.}
  \label{fig:LandingL}
\end{figure}

We claim that $y_\pm = x_{t,a}$. 
Let $U' \subseteq U$ be sufficiently small neighborhoods of $y$ with $R_a^{q}:U' \longrightarrow U$ is a homeomorphism.
Thus, for sufficiently large $n$, there exists $U_n' \subset U_n:= A_{a,n}(U)$ so that $f_n^{q}: U'_n \longrightarrow U_n$ is a homeomorphism.
Since $U'_n$ intersects the corresponding invariant ray $\gamma_{\pm, n}$ landing at $\eta_{\pm, n}(t)$ for sufficiently large $n$, by pulling back using $f_n^{q}: U'_n \longrightarrow U_n$, $\eta_{\pm, n}(t) \in U_n$.
Hence $A_{a,n}^{-1} (\eta_{\pm, n}(t)) \in U$.
Since $U$ is arbitrary, $y = x_{t,a}$.

If $t$ is aperiodic with pre-period $l$, then $s:=m_d^l(t)$ is periodic.
If $x_{s,a} \notin \Xi_a$, then by pulling back the Hausdorff limit $\gamma_\pm$, we conclude that $x_{t,a} \notin \Xi_a$.
If $x_{s,a} \in \Xi_a$, since the periodic singular points in $\Xi_a$ are not in critical periodic cycle, by pulling back $\gamma_\pm$, with finitely many exceptions, $x_{t,a} \notin \Xi_a$.
\end{proof}

Now suppose $t\in\pm I_{\pm, a} \cap \Q$ is not the end point of a leaf of $\mathcal{L}_\pm$, and is of Type (2) with pre-period $l$ and period $q$.
The itinerary of $t$ with respect to the partition $\mathbb{S}^1 = \bigcup_{a\in \RV^\mathcal{F}} \overline{I_{a,\pm}}$ is $(a_0 a_1... \overline{a_{l}...a_{l+q-1}})$ with $a_0 = a$.
We remark that this itinerary may be {\bf different} from the itinerary of $a_0$ under the dynamics of $F:\RV^\mathcal{F} \longrightarrow \RV^\mathcal{F}$, i.e., $a_{i+1}$ may be different from $F(a_i)$.
This happens when the end runs into a backward orbit of the singular set $\Xi$ (see the proof of Lemma \ref{lem:PV2} below for more discussions).

Let $t_i = m_d^i(t) \in \R/\Z$ and $\alpha_i = \alpha_\pm(t_i)$ be the corresponding end of $\mathcal{T}_\pm^\infty$.
To avoid too many subindices, we shall drop $\pm$ for $\alpha_i$ in the following.
Let $v_i \in \mathcal{T}_\pm$ be the projection of this end $\alpha_i$ onto $\mathcal{T}_\pm$.
By Assumption \ref{assumption:1}, we have
\begin{lem}
$v_i \in \mathcal{V}$ is a Fatou point and $[v_i] = a_i \in \RV$.
\end{lem}
\begin{proof}
Let $E$ be an edge incident at $v_i$ in the direction of $\alpha_i$, and let $E'$ be an edge of $\mathcal{T}_\pm$ incident at $v_i$.
Suppose for contradiction that $v_i$ is a Julia point. 
Then $v_i$ is eventually mapped to a periodic point of degree $1$.
Thus $E$ and $E'$ have the same pre-period, so $E'$ is an edge in $\mathcal{T}_\pm$ by Assumption \ref{assumption:1}.
This is a contradiction to $\alpha_i$ projects to $v_i$.

Since $\alpha_i$ projects to $v_i$, $\lim_{v_i} \eta_{\pm, n}^{\BP}(t_i) \notin \xi_{v_i}(T_{v_i}\mathcal{T}_\pm)$.
Thus $t_i \in \pm I_{\pm, v_i}$, so $[v_i] = a_i$.
\end{proof}

We assume $\alpha_i$ is the ray representing the end that starts at $v_i$.
Denote $\alpha_{i,n} \subseteq \T^\infty_{\pm,n}$ and $\alpha_{i,a} \subseteq \T^\infty_{\pm, a}$ be the corresponding limiting arcs.
Recall that by definition $\T^\infty_{\pm,a} = \bigcup_k \T^k_{\pm,a}$, so $\alpha_{i,a}$ is neither open nor closed.
\begin{lem}\label{lem:nonint}
For each $i=0,1..., q-1$, $\alpha_{i,a_i} \subseteq \hat\C_{a_i} - \Xi_{a_i}$.
\end{lem}
\begin{proof}
Since $\alpha_i$ starts at $v_i$ and $v_i$ is a Fatou point with $a_i = [v_i]$, $\alpha_{i,a_i}$ starts at a point in $\hat\C_{a_i} - \Xi_{a_i}$.
Suppose for contradiction $\alpha_{i,a_i} \cap \Xi_{a_i} \neq \emptyset$.
Then there exists $k$ so that $\alpha_{i,a_i} \cap \T^k_{\pm,a_i}\cap \Xi_{a_i} \neq \emptyset$.

By Lemma \ref{lem:spilqig}, each singular point $x\in \Xi_{a_i}$ is a vertex for some non-trivial edge of $\T_{\pm, {a_i}}$, thus the pre-period of $x$ is $\leq m$ where $m$ us the constant in the Assumption \ref{assumption:1}.
If $x$ is periodic, then a similar argument as Lemma \ref{lem:IEP} gives that any incident edge of $\T^k_{\pm, {a_i}}$ at $x$ corresponds to a periodic edge of $\mathcal{T}^k_\pm$.
Thus, in general, the incident edges of $\T^k_{\pm, a_i}$ at $x$ corresponds to edges of $\mathcal{T}^k_\pm$ of pre-period $\leq m$.
This is a contradiction as $\alpha_i$ contains no edges of pre-period $\leq m$ by Assumption \ref{assumption:1}.
\end{proof}

These rays $\alpha_{i,a_i} \subseteq \hat\C_{a_i} - \Xi_{a_i}$ play similar role as the invariant arcs in the proof of Lemma \ref{lem:PV1} for the Case (1).
We will show with finitely many exceptions, they land at $y_i \in \hat\C_{a_i} - \Xi_{a_i}$.
The landing point $y_i$ may correspond to a periodic point or a point in the backward orbit of the singular set $\Xi$.
We show that in either case, $y_i = x_{t_i, a_i}$, concluding the proof.

The rays $\alpha_i$ can be decomposed into
$$
\alpha_i = \widetilde{\alpha}_i\cup \bigcup_{k=0}^\infty \varepsilon_i^k,
$$
where $\widetilde{\alpha}_i$ is some initial segment, $F_\pm(\varepsilon_i^k) = \varepsilon_{i+1}^k$ if $i < q-1$ and $F_\pm(\varepsilon_{q-1}^k) = \varepsilon_0^{k-1}$.
In particular, $F_\pm^q(\varepsilon_i^k) = \varepsilon_i^{k-1}$.
We denote the corresponding arcs as $\varepsilon_{i,n}^k \subseteq \T_{\pm,n}^\infty$ and $\varepsilon_{i,a_i}^k \subseteq \T_{\pm,a_i}^\infty$.

\begin{lem}\label{lem:PV2}
Let $a\in \RP^\mathcal{F}$.
If $t\in \pm I_{\pm, a}$ is a rational angle of Type (2), then with finitely many exceptions, $x_{t,a} \in \hat \C - \Xi_a$.
\end{lem}
\begin{proof}
We assume $t$ is periodic under $m_d$. The aperiodic case can be proved with a similar modification as in Lemma \ref{lem:PV1}.

We consider two cases.

Case (a): $\varepsilon_{i,a_i}^k$ is never trivial, i.e, for all $i = 0, 1,..., q-1$ and $k\in \N$, $\varepsilon_{i,a_i}^k$ is not a singleton set.
Then $\varepsilon_{i,a_i}^k$ is a union of edges, and we denote the $U_{i}^k \subseteq \hat\C_{a_i}$ as the union of Carath\'eodory disks associated for the edges.
Then 
$$
R_{a_i}: U_i^k \longrightarrow U_{i+1}^k 
$$
for $k < q-1$ and 
$$
R_{a_{q-1}}: U_{q-1}^{k} \longrightarrow U_{0}^{k-1} 
$$
are homeomorphisms.
In particular, $F(a_i) = a_{i+1}$, and $R_{a_i}^q: U_i^k \longrightarrow U_{i}^{k-1}$ is a homoemorphism.
Note that $U_i^k$ and $U_i^{k-1}$ are connected at pre-periodic points, so we can choose a curve $\gamma_0^0 \subseteq U_i^0$ and get an invariant ray by pulling back by $R$
$$
\gamma_i = \bigcup_{k=0}^\infty \gamma_i^k;
$$
with $R_{a_i}(\gamma_i^k) = \gamma_{i+1}^k$ and $R_{a_{q-1}}(\gamma_{q-1}^k) = \gamma_{0}^{k-1}$.
Since $R_{a_i}^q: \hat\C\longrightarrow \hat\C$ is geometrically finite by Theorem \ref{thm:crm}, and $\gamma_0^0$ does not contain post-critical set, a standard argument using expansion in the complement of the post-critical set gives that each $\gamma_i$ converges to a point $y_i \in \hat\C_{a_i}$ (see Corollary A.6 in \cite{TY96}).
Now a similar argument as in Lemma \ref{lem:PV1} gives that with finitely many exceptions, $y_i \in \hat\C_{a_i} - \Xi_{a_i}$ and $y_i = x_{t_i, a_i}$, proving this case.

Case (b): $\varepsilon_{i,a_i}^k$ is trivial for some $i$ and $k$.

We first claim that if $\varepsilon_{i,a_i}^{k}$ is trivial and $\varepsilon_{i'}^{k'}$ is mapped to $\varepsilon_i^k$ by some iterate of $F_\pm$, then $\varepsilon_{i', a_{i'}}^{k'}$ is also trivial.
Without loss of generality, we assume $i \neq 0$, and it suffices to show $\varepsilon_{i-1, a_{i-1}}^{k}$ is trivial.

Suppose for contradiction that the claim is not true. 
Then $\varepsilon_{i, F_\pm(a_{i-1})}^{k} \subseteq \hat\C_{a_{i-1}}$ is non-trivial.
If $F_\pm(a_{i-1}) = a_i$, then we get a contradiction as $\varepsilon_{i,a_i}^{k} = \varepsilon_{i, F_\pm(a_{i-1})}^{k}$ is non-trivial.
Thus $F_\pm(a_{i-1}) \neq a_i$. Let $x \in \Xi_{a_i} \subseteq \hat\C_{a_i}$ be the point associated to $F_\pm(a_{i-1})$, then $\varepsilon_{i,a_i}^{k} = \{x\}$, contradicting to Lemma \ref{lem:nonint}.

\begin{figure}[ht]
  \centering
  \resizebox{0.7\linewidth}{!}{
    \def\svgwidth{\columnwidth}
    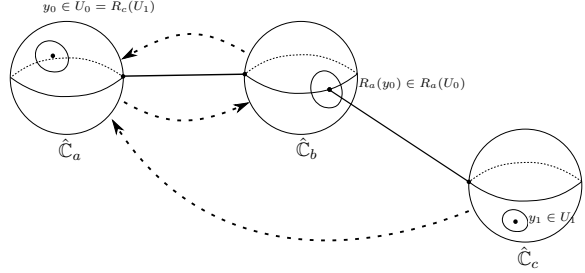

  }
  \caption{An illustration when $\varepsilon_{i,a_i}^k$ is trivial.}
  \label{fig:DS}
\end{figure}

Now assume that $\varepsilon_{i,a_i}^{k_0} = \{y_i\}$ for some $k_0$, then $\varepsilon_{i,a_i}^{k} = \{y_i\}$ for all $k\geq k_0$ as the end is connected.
By Lemma \ref{lem:nonint}, $y_i \notin \Xi_{a_i}$, so $y_i$ is not a hole of $R_{a_i}$ by Lemma \ref{lem:hcs}.
It is not a critical point, as any critical point that is not in $\Xi_{a_i}$ corresponds to some vertex in $\T_{\pm, a_i}$.
Let $U_i \subseteq \hat\C_{a_i}$ be a sufficiently small neighborhood of $y_i$ on which $R_{a_i}$ is univalent (see Figure \ref{fig:DS}).

If $a_{i+1} = F_\pm(a_i)$, then $R_{a_i}(U_i) \subseteq \hat\C_{F_\pm(a_i)}$ contains $y_{i+1}$.
Otherwise, $a_{i+1} \neq F_\pm(a_i)$, and $R_{a_i}(y_i)$ is the singular point in $\Xi_{F_\pm(a_i)}$ associated to $a_{i+1}$\footnote{We remark that since $\varepsilon_{i,a_i}^k$ is eventually trivial, this case must occur for some $i$. In particular, $F(a_i) \neq a_{i+1}$ for some $i$.}.
In any case, given $U_0$, we can inductively choose a neighborhood $U_i \subseteq\hat\C_{a_i}$ of $y_i$ so that (see Figure \ref{fig:DS})
\begin{itemize}
\item $R_{a_i}$ is univalent on $U_i$;
\item $U_{i+1} \subseteq_\RT R_{a_i}(U_i)$.\footnote{See Definition \ref{defn:subset} for the notation $\subseteq_\RT$.}
\end{itemize}
Thus, inductively, for sufficiently large $n$, there exists a sequence of $q$ disks $D_{i,n}$ with 
\begin{itemize}
\item $D_{i,n}$ is compactly contained in $A_{a_i,n}(U_i)$;
\item $R_n: D_{i,n} \longrightarrow D_{i+1,n}$ is a homeomorphism for $i =0,..., q-2$;
\item $R_n: D_{q-1,n} \longrightarrow A_{a_0,n}(U_0)$ is a homeomorphism.
\end{itemize}
A similar argument as in Lemma \ref{lem:PV1} gives that $x_{t_i, a_i} \in U_i$.
Since $U_i$ is abitrary, $x_{t_i, a_i} = y_i$ proving this case.
\end{proof}

Since there are only finitely many periodic points as an endpoint of a leaf of $\mathcal{L}_\pm$, similar arguments as above for the aperiodic points give
\begin{lem}\label{lem:PV3}
Let $a\in \RP^\mathcal{F}$.
If $t\in \pm I_{\pm, a}$ is an end point of a leaf for $\mathcal{L}_\pm$, then with finitely many exceptions, $x_{t,a} \in \hat \C - \Xi_a$.
\end{lem}

\begin{proof}[Proof of Proposition \ref{prop:pc}]
Combining Lemma \ref{lem:PV1}, Lemma \ref{lem:PV2} and Lemma \ref{lem:PV3}, we get the proposition.
\end{proof}

\appendix
\section{Applications on polynomial matings}
In this appendix, we show how these convergence results can be used to show the existence of polynomial matings.

A polynomial $P(z) = a_dz^d+...+a_0$ is said to be {\em monic} if $a_d = 1$ and {\em centered} if $a_{d-1} = 0$.
Let $\MP_d$ be the space of all monic and centered polynomials.
A monic and centered polynomial with connected Julia set has a unique choice of the B\"ottcher map normalized so that the derivative at infinity is $1$.
Thus, if $P$ has locally connected Julia set, we have a {\em marking} on the Julia set, i.e., a particular choice of the conjugacy $\eta_P: \R/\Z \longrightarrow \mathcal{J}(P)$ with $\eta_P \circ m_d = P \circ \eta_P$.
In this case, the external rays also defines a {\em polynomial lamination} $\mathcal{L}_P$ for $P$ (see \cite{Thurston19}).

A rational map $R$ is a {\em geometric mating} (or briefly a {\em mating}) of two polynomials $P_\pm \in \MP_d$ with connected and locally connected Julia sets if we have
\begin{itemize}
\item a decomposition of the Fatou set $\mathcal{F}(R) = \mathcal{F}_+ \sqcup \mathcal{F}_-$ with $\mathcal{J}(R) = \partial \mathcal{F}_+ = \partial \mathcal{F}_-$; and
\item two continuous surjections on the filled Julia sets $\psi_\pm: \mathcal{K}(P_\pm) \rightarrow \overline{\mathcal{F}_\pm}$ that are conformal between $\Int(\mathcal{K}(P_\pm))$ and $\mathcal{F}_\pm$
\end{itemize}
so that
\begin{align*}
\psi_\pm \circ P_\pm &= R \circ \psi_\pm\\
\psi^+ \circ \eta^+ &= \psi^- \circ \eta^- \circ m_{-1}
\end{align*}
where $\eta^{\pm}: \R/\Z \rightarrow \mathcal{J}(P_\pm)$ are the markings and $m_{-1}$ is the multiplication by $-1$ map.
Two polynomials $P_\pm \in \MP_d$ are said to be {\em mateable} if there exists a rational map $R$ that is a mating of them.

Given a geometrically finite polynomial $P \in \MP_d$ with connected Julia set, there exists a unique post-critically finite polynomial $\hat P \in \MP_d$ with conjugate dynamics on the Julia sets with compatible markings \cite{Haissinsky00}.

The dynamics of $\hat P$ can be described combinatorially by the dynamics on its Hubbard tree $\hat P: H \rightarrow H$.
We say $P$ has simplicial Hubbard tree if the corresponding dynamics $\hat P: H \rightarrow H$ is simplicial: there exists a finite simplicial structure on $H$ so that $\hat P$ sends an edge of $H$ to an edge of $H$.
We will show:
\begin{theorem}\label{thm:mating}
Let $P_\pm \in \MP_d$ be two geometrically finite polynomials, where $P_-$ has simplicial Hubbard tree.
Then $P_\pm$ are mateable if and only if their laminations $\mathcal{L}_\pm$ are non-parallel.
\end{theorem}

The only if direction of the Theorem is immediate: if $\mathcal{L}_\pm$ are parallel, then there is a topological obstruction on mating coming from the Levy cycles.

To prove the if direction, we sketch the idea as follows.
We first define an equivalence relation on geometrically finite polynomial in $\MP_d$ with connected Julia set by
$$
P \sim_\mathcal{J} Q
$$
if $P$ and $Q$ have conjugate dynamics on the Julia sets with compatible markings.
By a Theorem of Ha\"issinsky and Tan Lei \cite{HT04}, $P_\pm$ are mateable if and only if some (and hence all) geometrically finite polynomials $Q_\pm \in [P_\pm]_\mathcal{J}$ are mateable.
Thus, it suffices to construct a mating between some $Q_\pm \in [P_\pm]_\mathcal{J}$.

We start with a post-critically finite polynomial $\hat P_+ \in [P_+]_\mathcal{J}$.
By Theorem 1.3 and Proposition 6.2 in \cite{Luo21}, there exists quasi post-critically finite $(\bp_n)_n \in \BP_d$ whose dual lamination $\mathcal{L}_\bp$ generates the same equivalence relation as the polynomial lamination $\mathcal{L}_-$.
We perform quasi post-critically finite on the unbounded Fatou component using this sequence $\bp_n$. 
In other words, we consider the matings $\hat P_+ \sqcup \bp_n$.
We show the sequence has a convergent subsequence if the laminations $\mathcal{L}_+$ and $\mathcal{L}_{\bp}$ are non-parallel.
Let $f$ be the limiting rational map.
We show $f$ is a mating of $\hat P_+$ and $Q_-$ for some $Q_- \in [P_-]_\mathcal{J}$ and conclude the proposition.

\subsection*{Degenerations on the unbounded Fatou component}
Let $\hat P_+ \in \MP_d$ be a post-critically finite polynomial, and $(\bp_n)_n \in \BP_d$ be quasi post-critically finite.
Let $f_n = \hat P_+ \sqcup \bp_n$ be the mating of $\hat P_+$ and $\bp_n$.

Although in \S \ref{sec:qpcfd}, we only deal with quasi post-critically finite degenerations for hyperbolic rational maps, the construction of rescaling limits also works for degenerations of geometrically finite rational maps $f_n$.

More precisely, let $F: (\mathcal{T}, \p) \longrightarrow (\mathcal{T}, \p)$ be the corresponding quasi-invariant tree for $\bp_n$, with vertex set $\mathcal{V}$.
The equivalence classes of rescalings $A_{v,n} \in \PSL_2(\C)$ for $v\in \mathcal{V}$ gives a set $\RP = \mathcal{V}/\sim$ (see \S \ref{sec:qpcfd} for more details), and an induced map
$F: \RP \longrightarrow \RP$ 
with rescaling rational maps 
$$
R_a: \hat\C_a \longrightarrow \hat\C_{F(a)}.
$$
To get to the result in the shortest path, although one can, we don't need to define the rescaling trees and rescaling tree maps.

Let $\Xi_a \subseteq \C_a$ be the singular set at $a\in \RP$.
The Fatou points and Julia points are defined similarly $\RP= \RP^\mathcal{F} \sqcup \RP^\mathcal{J}$.
The same proof of Theorem \ref{thm:hl} gives that after passing to a subsequence, for any $a\in \RP^\mathcal{F}$, the quasi-invariant tree $\T_n$ converges to a graph $\T_a$ with finitely many vertices and finitely many non-loop edges.

By considering a subtree of the pullback, we assume that $\mathcal{T}$ satisfies the Assumption \ref{assumption:1}.
The proof of Lemma \ref{lem:hcs} needs a modification:
 \begin{lem}\label{lem:hcs2}
Let $a\in \RP^\mathcal{F}$. Then the holes of $R_a$ are contained in $\Xi_a$.
\end{lem}
\begin{proof}
Let $x\in \hat\C_a$ be a hole of $R_a$. We consider two cases.

If $x\notin \T_a$, then there exists a small simple closed curve $C$ bounding $x$ so that the corresponding disk $D$ does not intersect $\T_a$.
Modify $C$ if necessary, we may assume $R_a(C)$ is a simple closed curve.
Shrink $C$ if necessary, we may assume for sufficiently large $n$, $A_{F(a),n}(R_a(D))$ does not contain critical values in the unbounded Fatou component $\U_n$ of $f_n$.
Therefore, each component of $f_n^{-1}(A_{F(a),n}(R_a(D)))$ is a disk.
Let $D_n$ be the component with $\partial D_n \to_a C$, then $f_n(D_n)\neq \hat \C$.
Thus $A_{F(a),n}^{-1}\circ f_n \circ A_{a,n}(D)\neq\hat\C$ for all large $n$. This gives a contradiction to $x$ is a hole (see Lemma 4.5 in \cite{DeM05}).

If $x\in \T_a$, since $a\in \RP^\mathcal{F}$, $\T_a$ is non-trivial, we can choose a small simple closed curve $C$ intersecting an edge of $\T_a$.
Let $C_n \to_a C$ with $f_n(C_n)$ is a simple closed curve.
Since $x$ is a hole, there exists $C_n'$ with 
\begin{itemize}
\item $C_n' \to_a x$;
\item $f_n(C_n) = f_n(C_n')$.
\end{itemize}
By Assumption \ref{assumption:1}, $C_n'$ intersects some edge $E_n$ of $\T_{n}$, so $E_n \to_a x$. Therefore $x \in \Xi_a$.
\end{proof}

The above plus a similar proof as in Proposition \ref{prop:gfl} gives
\begin{cor}
The union of rational maps $R: \hat\C^{\RP^\mathcal{F}} \longrightarrow \hat\C^{\RP^\mathcal{F}}$ is geometrically finite.
\end{cor}

As in \S \ref{sec:ptq}, we can define a partition by
$$
I_a := \bigcup_{[v] = a} \overline{\Int(I_v)} \subseteq \R/\Z,
$$
for $a\in \RP^\mathcal{F}$ and
$$
\R/\Z = \bigcup_{a\in \RP^\mathcal{F}} I_a.
$$
Let $\eta_n: \R/\Z \longrightarrow \mathcal{J}(f_n)$ be the {\em markings} of the Julia set for $f_n$, i.e., the unique semi-conjugacy between $m_d$ and $f_n$ that compatible with the markings of the unbounded Fatou component.
After passing to a subsequence, and a diagonal argument, we assume 
$$
x_{t,a} := \lim_a \eta_n(t) = \lim_{n\to\infty} A_{a,n}^{-1}(\eta_n(t))
$$ 
exist for all rational angles $t \in \R/\Z$ and all $a\in \RP^\mathcal{F}$.
Then the same proof of Proposition \ref{prop:pc} gives that 
\begin{lem}
Let $a\in \RP^\mathcal{F}$. For all but finitely many rational angles $t\in -I_a$, $x_{t,a} \in \hat \C_a - \Xi_a$.
\end{lem}

\begin{prop}\label{prop:cnp}
If $[f_n] = [\hat P_+ \sqcup \bp_n]$ diverges, then $\mathcal{L}_+$ and $\mathcal{L}_\bp$ are parallel.
\end{prop}
\begin{proof}
Since $[f_n] = [\hat P_+ \sqcup \bp_n]$, by the above construction, we get a non-trivial partition of angles $\R/\Z = \bigcup_{a\in \RP^\mathcal{F}} I_a$.

We claim that if $t\in -\partial I_a$, $t$ lands at a cut point for $\hat P_+$.
There are two cases to consider.
Case 1: If $t$ lands on the boundary of a Fatou component $U$, then there exists $t_n < t < t_n'$ with $t_n, t_n' \to t$ such that $t_n, t_n'$ all land at a Fatou component $U$.
Case 2: Otherwise, there exists $t_n < t < t_n'$ with $t_n, t_n' \to t$ such that  $(t_n, t_n')$ is a leaf of $\mathcal{L}_+$.

The first case is not possible, as the partition of $\hat\C^{\RP^\mathcal{F}}$ does not separate the corresponding limit of Fatou component of $U$, so rational external angles landing on $\partial U$ are all contained in the same $I_a$.
The second case is not possible either as rational angles $t, t'$ landing at the same point for $\hat P_+$ would belong to the same $I_a$.
 
Similar argument also shows that there exists $t' \neq t \in -\partial I_a$ landing at the same point as $t$ for $\hat P_+$. 
Thus, two laminations are parallel.
\end{proof}

\subsection*{Construction of a mating}
Let $P_\pm \in \MP_d$ be two geometrically finite polynomials, where $P_-$ has simplicial Hubbard tree such that their polynomial laminations $\mathcal{L}_\pm$ are non-parallel.
By Theorem 1.3 and Proposition 6.2 in \cite{Luo21}, there exists quasi post-critically finite $(\bp_n)_n \in \BP_d$ whose dual lamination $\mathcal{L}_\bp$ generates the same equivalence relation as the polynomial lamination $\mathcal{L}_-$. 
In particular, $\mathcal{L}_+$ and $\mathcal{L}_\bp$ are non-parallel.

Let $f_n = \hat P_+ \sqcup \bp_n$ where $\hat P_+ \in [P_+]_\mathcal{J}$ is post-critically finite.
By Proposition \ref{prop:cnp}, after passing to a subsequence, $[f_n]$ converges to a rational map $[f]$ of the same degree.
We shall now verify that $f$ is indeed a mating of $\hat P_+$ with $Q_- \in [P_-]_\mathcal{J}$.

Let $v\in \mathcal{V}^\mathcal{F}$, then the corresponding Carath\'eodory disk $\U_v \subseteq \Omega_v$ for some Fatou component of $f$. 
Using the limiting quasi-invariant graphs (see Proposition 6.12 and Proposition 6.13 in \cite{Luo21}), one can prove 
\begin{prop}\label{prop:dfcr}
If $v\neq w\in \mathcal{V}^\mathcal{F}$, then
$\Omega_v \neq \Omega_w$.
Moreover, if $t, t'$ form a leaf of $\mathcal{L}_\bp$, then the corresponding pre-periodic points $x_{t,n}$, $x_{t',n}$ for $f_n$ converge to the same pre-periodic point $x$ for $f$.
\end{prop}

\begin{proof}[Proof of Theorem \ref{thm:mating}]
It suffices to prove the if direction.
Let $f$ be the limit as $[f_n] = [\hat P_+ \sqcup \bp_n]$.
Then $f$ is geometrically finite.
By proposition \ref{prop:dfcr}, each Fatou component is simply connected so $f$ has connected Julia set, so it has locally connected Julia set \cite{TY96}.
Let $\mathcal{K}_-$ be the closure of unions of Fatou components $\Omega_v$ where $v$ is a Fatou point of $\mathcal{T}^\infty$, then $\partial \mathcal{K}_- = \mathcal{J}(f)$.
Let $\mathcal{K}_+ = \mathcal{J}(f) \cup (\hat\C- \mathcal{K}_-)$, then $\partial \mathcal{K}_+ = \mathcal{J}(f)$.
By Proposition \ref{prop:dfcr}, we have a continuous surjection
$$
\psi_- : \mathcal{K}(\hat P_-) \longrightarrow \mathcal{K}_-
$$
which semi-conjugates the dynamics of $\hat P_-$ with $f$ on the Julia sets.
Indeed, by Proposition \ref{prop:dfcr}, we can construct a semi-conjugacy between the closures of unions of finitely many components for Fatou points in $\mathcal{T}_-$, and then use the dynamics to pull back the semi-conjugacy.
Since $f$ is geometrically finite, the diameters of the pullback Fatou components are shrinking to $0$\cite{TY96}, giving a continuous map in the limit.

If $t, t' \in \Q/\Z$ land at the same point for $\hat P_+$, then the corresponding point $x_{t,n} = x_{t',n}$ for all $n$. 
Thus, the equation $\psi^+ \circ \eta^+ = \psi^- \circ \eta^- \circ m_{-1}$ gives a continuous map on pre-periodic points on pre-periodic points on $\mathcal{J}(\hat P_+)$, which extends to a continuous surjection $\psi_+: \mathcal{K}(\hat P_+) \longrightarrow \mathcal{K}_+$ semi-conjugating the dynamics of $\hat P_+$ with $f$ on the Julia sets.

Changing $\hat P_-$ to $Q_- \in [P_-]_\mathcal{J}$, we may modify the conjugacy $\psi_\pm$ to be conformal in the interior of the filled Julia set.
Thus $f$ is a mating of $\hat P_+$ and $Q_-$.
By Theorem D in \cite{HT04}, $P_\pm$ are mateable.
\end{proof}


\newcommand{\etalchar}[1]{$^{#1}$}

\end{document}